\documentclass[11pt]{article} 
\usepackage[utf8]{inputenc}

\usepackage{float}
\usepackage{tikz-cd}
\usepackage{amsxtra}
\usepackage{amsmath}
\usepackage{amssymb}
\usepackage{amsfonts}
\usepackage{amscd}
\usepackage{extarrows}
\usepackage{mathtools}
\usepackage[all,2cell]{xy}
\usepackage{amsthm}
\usepackage[shortlabels]{enumitem}
\usepackage{hyperref}
\usepackage{cleveref}
\usepackage{subfigure}
\usepackage{euscript}
\usepackage[margin=3cm]{geometry}
\usepackage[bbgreekl]{mathbbol}
\usepackage{mathabx}
\usepackage{graphicx}
\graphicspath{{pics/}}
\usepackage{booktabs}

\usepackage{chngcntr}
\counterwithin{figure}{subsection}

\usepackage{xstring} 

\numberwithin{equation}{subsection}

\newtheorem{thm}[equation]{Theorem}
\newtheorem*{thm*}{Theorem}
\newtheorem{cor}[equation]{Corollary}
\newtheorem{lem}[equation]{Lemma}
\newtheorem{prop}[equation]{Proposition}

\theoremstyle{definition}

\newtheorem{defi}[equation]{Definition}
\newtheorem{rem}[equation]{Remark}

\newtheorem{exa}[equation]{Example}
\newtheorem*{exa*}{Example}
\newtheorem{con}[equation]{Construction}
\newtheorem{conj}[equation]{Conjecture}

\usepackage{xcolor}%
\hypersetup{%
    colorlinks,%
    linkcolor={red!50!black},%
    citecolor={blue!50!black},%
    urlcolor={blue!80!black}%
}

\usepackage{macros/utility-macros}
\usepackage{macros/diagram-macros}
\usepackage{macros/macros}

\setlist[enumerate,1]{label=(\arabic{*})}
\setlist[enumerate,2]{label=(\alph{*})}
\setlist[enumerate,3]{label=(\roman{*})}

\usepackage{subfiles}

\title{Complexes of stable $\infty$-categories} 
\author{Merlin Christ\footnote{Université Paris Cité and Sorbonne Université, CNRS, IMJ-PRG, F-75013 Paris, France, email: {\tt merlin.christ@imj-prg.fr}}, Tobias
Dyckerhoff\footnote{Universität Hamburg, Fachbereich Mathematik, Bundesstraße 55, 20146 Hamburg, Germany
	email: {\tt tobias.dyckerhoff@uni-hamburg.de}}, Tashi
	Walde\footnote{Technische Universität München, Zentrum Mathematik, Boltzmannstr. 3, 85748 Garching, Germany, email: {\tt tashi.walde@ma.tum.de}}}

\begin{document}

\maketitle

\begin{abstract}
	We study complexes of stable $\infty$-categories, referred to as categorical complexes.  As
	we demonstrate, examples of such complexes arise in a variety of subjects including
	representation theory, algebraic geometry, symplectic geometry, and differential topology.
	One of the key techniques we introduce is a totalization construction for categorical cubes
	which is particularly well-behaved in the presence of Beck-Chevalley conditions. As a direct
	application we establish a categorical Koszul duality result which generalizes previously
	known derived Morita equivalences among higher Auslander algebras and puts them into a
	conceptual context. We explain how spherical categorical complexes can be interpreted as
	higher-dimensional perverse schobers, and introduce Calabi-Yau structures on
	categorical complexes to capture noncommutative orientation data. A variant of homological
	mirror symmetry for categorical complexes is proposed and verified for $\CP^2$. 
\end{abstract}

\tableofcontents



\section{Introduction}
\label{sec:introduction}

While originally intended as a habitat to formulate and prove homological duality results, derived
categories have since evolved into invariants in their own right. One modern perspective is to treat
derived categories as {\em stable $\infty$-categories}, as introduced by J. Lurie \cite{lurie:ha},
which allows for the application of powerful tools from homotopy theory.
In this work, we investigate {\em complexes} of stable $\infty$-categories, referred to as {\em categorical
complexes}. These are sequences
\[
	\begin{tikzcd}
		\cdots \ar{r}{d} & \A_{2} \ar{r}{d} & \A_1 \ar{r}{d} & \A_{0} \ar{r}{d} & \cdots
	\end{tikzcd}
\]
of stable $\infty$-categories $\A_i$ with exact functors $d$ such that $d^2 \simeq 0$. In the form
of {\em localization sequences}, special examples of such complexes have already been crucial in the
study of algebraic $K$-theory and related additive or localizing invariants
\cite{waldhausen,thomason,BGT13}. 
Several developments, mainly within the vicinity of homological mirror symmetry, suggest the
existence of interesting more general classes of examples of categorical complexes and arouse the
desire to investigate them within an effective axiomatic framework of {\em categorified homological
algebra}.
Here, the primary interest lies in the stable $\infty$-categories {\em themselves} and not just
their additive invariants, so that the familiar axiomatic frameworks for homological algebra can not be
applied but needs to be replaced by a suitable $(\infty,2)$-categorical counterpart. 

In this article, using this circle of ideas as a starting point, we
\begin{itemize}
	\item construct classes of examples of categorical complexes arising in several
		different areas, 
	\item study various totalization constructions and their interplay with Beck--Chevalley
		conditions,
	\item demonstrate the usefulness of these constructions by proving concrete results, such as
		categorical Koszul duality, 
	\item formulate conjectures in the context of homological mirror symmetry
	\item and discuss properties (sphericalness) and structures (Calabi--Yau) of categorical
		complexes.
\end{itemize}
In the companion paper \cite{CDW24}, we furthermore introduce a notion of lax $2$-additivity for $(\infty,2)$-categories which we believe to be a reasonable axiomatic framework for the study of the purely additive aspects of categorical complexes. We emphasize, that these are merely preliminary steps: In this work, we do not offer well--behaved
$(\infty,2)$-categorical counterparts of homotopical or homological aspects, in particular, we do not
introduce an $(\infty,2)$-categorical version of the derived category --- we hope to return to this question in future work.

We proceed with a more detailed overview of the contents of this work. 

\subsection{Rules of categorification}
\label{subsec:rules-of-categorification}

Most examples and results of this work are based on the insight that some of the building
blocks of the construction of complexes of abelian groups can be ``categorified'', replacing abelian
groups by stable $\infty$-categories. This procedure of categorification typically follows the
following set of rules:
{
\renewcommand{\arraystretch}{3}
\begin{center}
\begin{tabular}{lll}
	 & classical & categorified\\ \hline
	1) & abelian group $A$ & stable $\infty$-category $\A$\\
	2) & element $x \in A$ & object $X \in \A$\\
	3) &  $y - x$ & $\cone(X \overset{f}{\to} Y)$\\
	4) &  $\sum (-1)^i x_i$ & $\tot(\begin{tikzcd}[column sep=small] X_0  \ar{r}{d} &  X_1 \ar{r}{d} & X_2 \ar{r}{d} & \cdots
	 \ar{r}{d}&  X_n\end{tikzcd})$\\
		5) & \parbox{6cm}{direct sum decomposition\\$C \cong A\oplus B$} &
		\parbox{6cm}{semiorthogonal decomposition\\$\C \simeq \langle \A, \B \rangle$}\\
		6) & external direct sum $A\oplus B$ &
		lax sum $\laxsum{\A}{F}{\B}$
\end{tabular}
\end{center}
}
\vspace{2ex}

While the first two rules should be apparent, we start commenting on rule 3). This is a first
crucial difference between the classical and the categorified context: In order to take a
``difference'' between objects $X,Y$ of a stable $\infty$-category $\A$, we need to be given the
additional datum of a morphism $f\colon X \to Y$ -- the difference will then be the cone of $f$.
Compliance with this rule will force us to include certain $2$-categorical data which becomes
invisible upon passing to the Grothendieck group $K_0$. As we will see in examples, this typically results in rather
natural {\em lax variants} of $1$-categorical constructions.  

Rule 4) is a natural generalization of Rule 3): An alternating sum over $n$ elements will
be categorified by the totalization of an $n$-term complex in $\A$. Here we do not only need to
specify the differentials of this complex, but also a coherent system of null homotopies -- this is
necessary to make sense of the totalization in the $\infty$-categorical context.

Rule 5) is almost evident after having accepted Rule 3): While in a direct sum $A \oplus B$, every
element is uniquely the sum of elements from the components $A$ and $B$, respectively, in a
semiorthogonal decomposition $\langle \A, \B \rangle$, every object is uniquely an {\em extension} of
an object $A \in \A$ by an object $B \in \B$. Put differently, by shifting the exact triangle
of the extension, every object is uniquely the cone of a morphism $A[-1] \to B$, thus
connecting back to Rule 3). 

Conceptually distinct to a direct sum {\em decomposition} of a given abelian group are the {\em
universal properties} satisfied by the {\em external} direct sum of a pair of abelian groups.  When
formulating a categorified analog of these universal properties, it is not sufficient to just
provide a pair of stable $\infty$-categories: As an additional datum, we need to specify a functor
$F:\A \rightarrow \B$ (similar to the additional choice of a morphism $f\colon X \to Y$ needed in Rule
3). The categorified ``direct sum'' is then the {\em lax sum} 
\[
	\laxsum{\A}{F}{\B}
\]
of the diagram of stable $\infty$-categories described by $F$ (see Rule 6). As will be explained in \S \ref{intro:laxadd} and \S \ref{subsec:laxsum}, the lax sum is simultaneously
a lax limit and lax colimit ({\em lax semiadditivity}), as well as an oplax limit and oplax
colimit ({\em lax additivity}) of the diagram $\Delta^1 \to \Stk$ given by $F$.  The lax sum admits
a semiorthogonal decomposition with components $\A$ and $\B$ (which according to Rule 5 should be
interpreted as a categorified direct sum {\em decomposition}). Vice versa, a stable
$\infty$-category with a semiorthogonal decomposition can be described as a lax sum if and only
if it admits a {\em gluing functor}. 

\subsection{Examples of categorical complexes}\label{intro:ex}

We survey the main examples of {\em categorical complexes} (i.e.~chain complexes of stable
$\infty$-categories) that will be constructed in this article. 

\subsubsection{The categorified Dold--Kan correspondence}

As explained in \cite{dyck:dk}, the classical Dold--Kan correspondence can be categorified along the
lines of the rules explained in \S \ref{subsec:rules-of-categorification}. This implies that every
connective categorical complex arises as the {\em totalization of a $2$-simplicial stable
$\infty$-category}. The relevant constructions will be reviewed in some more detail in \S \ref{sec:dk}.
The appearance of categorical complexes in this context can be viewed as a first evidence that one may
imitate somewhat more elaborate constructions from homological algebra in an $(\infty,2)$-categorical context.

\subsubsection{Categorical Koszul complexes and representation theory}

The classical construction of Koszul complexes by means of totalizing tensor products of two-term
complexes can be categorified resulting in a notion of {\em categorical Koszul complex}. The details
will be explained in \S \ref{sec:koszul}. In this context, we prove a Koszul duality statement
(Theorem \ref{thm:catkoszul}) which amounts to a categorification of the classical self-duality of
the Koszul complex. As a consequence, we obtain a conceptual proof of a derived Morita duality,
originally discovered in \cite{Bec18}, among the {\em higher Auslander algebras} as introduced in
\cite{iyama:higher}. 

\subsubsection{Categorical intersection complexes and lax algebraic geometry}
\label{intro:complete}

Given a normal crossings divisor in a smooth projective variety, we consider its \v{C}ech nerve, formed by
the smooth components of the divisor along with their intersections. Passing to derived categories
of coherent sheaves, we obtain a cubical diagram of stable $\infty$-categories which may be
totalized to yield a categorical complex. We refer to this complex as the {\em categorical intersection
complex} of the divisor, introduced in \S \ref{subsec:completeintersection}. These complexes
naturally appear on the $B$-side of homological mirror symmetry and are part of a {\em homological
mirror symmetry conjecture for categorical complexes}. The terms of the categorical intersection
complex can be described somewhat appealingly as a hybrid of ``commutative data'' (coherent sheaves
on a smooth variety) and ``non--commutative data'' (lax gluing data). We feel that these types of
stable $\infty$-categories could be of independent interest for a theory of ``lax algebraic geometry''.

In this context, it is of particular relevance that categorical complexes are capable of capturing
noncommutative orientation data in the form of Calabi--Yau structures and have geometric
interpretations in terms of higher--dimensional perverse schobers. This will be explained in \S
\ref{intro:cy} and \S \ref{intro:spherical}. 

The construction of the categorical intersection complex arises from a general totalization
procedure which can be applied to any cubical diagram of stable $\infty$-categories (see \S
\ref{intro:tot}). In particular, we obtain a categorical complex associated to the cubical
resolutions of singularities familiar from mixed Hodge theory. One can expect that these
categorical complexes yield some variant of noncommutative resolutions of
singularities, but this will not be explored here. 

\subsubsection{Fukaya-Seidel complexes in symplectic geometry}
\label{intro:fs}

Let $X \subset \CC^{N}$ be an $n$-dimensional smooth affine subvariety and let $X_1 \subset X$ be a
generic hyperplane section of $X$ cut out by the equation $\pi_1(x) = 1$ for some linear function
$\pi_1 \colon \CC^N \to \CC$. In \cite{Sei08}, using further assumptions and auxiliary data, a left exact
sequence of stable $\infty$-categories (in loc. cit. these are modeled as $A_{\infty}$-categories)
\begin{equation}\label{eq:fs}
	\Fuk(X) \hra \FS(\pi_1) \overset{\partial}{\lra} \Fuk(X_1)
\end{equation}
is describes, which relates the Fukaya categories of $X$ and $X_1$ with the Fukaya--Seidel category $\FS(\pi_1)$. We may
iterate this construction by choosing a generic hyperplane section $X_2$ of $X_1$, corresponding to
a linear function $\pi_2$, and so on. Splicing together the various left exact sequences
\eqref{eq:fs} yields a categorical complex 
\[
	\FS(\pi_1) \lra \FS(\pi_2) \lra ... \lra \Fuk(X_n),
\]
which we call the {\em Fukaya--Seidel complex} of the family $\{\pi_{\bullet}\}$ (see \S
\ref{subsec:FukayaSeidelcplx}).

It is a consequence of classical Picard--Lefschetz theory, that the passage to Grothendieck groups
yields a complex of abelian groups with homology $\on{H}_{*}(X)$. This suggests that, in comparison
with Fukaya categories, we may think of Fukaya--Seidel complexes as categorical invariants that
attempt to capture information about the full homology as opposed to just the middle--dimensional
one. In the given context, there is one caveat though: Currently it is not really justified to use
the term ``invariant'' here, since it is unclear to us at the moment, in which sense the
Fukaya--Seidel complex is independent of the choice of the family $\{\pi_{\bullet}\}$.

Fukaya--Seidel complexes are (conjecturally) the symplectic counterparts of categorical intersection
complexes in the context of homological mirror symmetry. Thus, just like categorical intersection
complexes of anticanonical normal crossing divisors, they typically carry Calabi--Yau structures and
naturally fit within a higher--dimensional framework of perverse schobers (cf. \S \ref{intro:cy} and
\S \ref{intro:spherical}).

\subsubsection{Complexes of local systems on manifolds with corners}
\label{intro:corners}

We consider the $n$-simplex $\Delta^n$ as an $n$-category, e.g.~via the $n$th oriental. A functor
$\M$ from $\Delta^n$ into the $(\infty,n)$-category $\on{Bord}_n$ of bordisms, as defined in
\cite{Lur:tft,CS19}, amounts to an $n$-manifold $M$ with corners together with an iterative decomposition of its boundary into smaller dimensional bordisms. The appearing submanifolds of $M$ are parametrized by the barycentric subdivision of the $n$-simplex, so that we may interpret this data, augmented by the empty space, as a
cubical diagram of topological spaces.  Similarly to \S \ref{intro:complete}, we may then pass to
the corresponding cube of stable $\infty$-categories of $\mathcal{D}(k)$-valued local systems and
totalize to obtain the {\em categorical complex of local systems on $\M$}, introduced in \S
\ref{subsec:localsystems}. We will see that an orientation of the $n$-manifold $M$ yields a Calabi--Yau structure on this complex. 

\subsection{Totalization constructions}\label{intro:tot}

Many of the examples of categorical complexes described in \S \ref{intro:ex} arise as totalizations
of categorical multi--complexes. In \S \ref{sec:tot}, we will discuss such totalization
constructions systematically. In particular, we explain that there are (at least) two natural
choices of totalizations, the {\em product} and {\em coproduct totalizations}. While, in general,
these totalizations lead to different categorical complexes, we show that they do agree in
the presence of suitable Beck--Chevalley conditions. 

\subsection{Mirror symmetry for categorical complexes}

We formulate homological mirror symmetry conjectures relating the Fukaya-Seidel
complexes from \S \ref{intro:fs} on the symplectic side to the categorical intersection
complexes from \S \ref{intro:complete} on the algebro--geometric side. We work out some details for
this {\em homological mirror symmetry for categorical complexes} in the simplest  illustrative case of Hori-Vafa mirror symmetry for $\CP^2$. The claimed equivalence for $\CP^n$, $n
\ge 3$, is already conjectural.  

\subsection{Geometric structures associated with categorical complexes}

\subsubsection{Calabi-Yau structures}
\label{intro:cy}

The concept of a left Calabi--Yau structure on a functor between $k$-linear stable
$\infty$-categories, as introduced for $k$ a field in \cite{BD19} and described for $k$ an $\mathbb{E}_\infty$-ring spectrum in \cite{Chr23}, can be regarded as a noncommutative
generalization of an oriented manifold with boundary. In fact, one example is precisely obtained
from such a manifold by considering the left Kan extension of derived local systems along the
boundary inclusion. 

In \S \ref{sec:CY}, we introduce the notion of a left Calabi--Yau structure on a categorical
complex, which can be thought of as a noncommutative analog of an oriented manifold with corners. It
generalizes the left Calabi--Yau structures from \cite{BD19} when interpreting the given functor as
a categorical $2$-term complex. The constructions of Calabi--Yau structures given in \cite{BD19} can
be generalized to provide examples of Calabi--Yau structures on categorical complexes:
\begin{enumerate}
	\item the categorical intersection complexes from \S \ref{intro:complete} associated with anticanonical divisors in projective varieties, see \Cref{cor:anticantotCY},
	\item conjecturally, the Fukaya--Seidel complexes from \S \ref{intro:fs},
	\item the categorical complexes of local systems associated to $n$-bordisms in \S
		\ref{intro:corners}, see \Cref{cor:bord},
\end{enumerate}

To produce Calabi--Yau structures on categorical complexes, we also introduce a notion of a Calabi--Yau structure on categorical cubes, with the property that the coproduct totalization of a Calabi--Yau cube inherits a Calabi-Yau structure. The Calabi--Yau structures in the examples (1) and (3) above arise in this way. 

\subsubsection{Spherical complexes and perverse schobers} 
\label{intro:spherical}

Spherical functors arise in several areas, in particular in homological mirror symmetry. As
proposed in \cite{KS14}, one may interpret them geometrically as categorical analogs of perverse
sheaves, called perverse schobers, on the complex plane $\CC$, stratified as $\CC \setminus \{0\}
\cup \{0\}$. 

As discussed in \S \ref{sec:spherical}, it turns out that the natural generalization of the concept
of a spherical functor to categorical complexes is simply given by a {\em spherical complex}: a
complex where all differentials are spherical functors. Spherical $n$-term complexes correspond geometrically to
categorified perverse sheaves on $\CC^n$ stratified by the linear subspaces
\[
	\{0\} \subset \mathbb{C} \subset \mathbb{C}^2 \subset \dots \subset \CC^n
\]
and their complements.

Just as for Calabi--Yau structures, many of the examples of spherical functors admit natural generalizations
to spherical complexes. These spherical complexes often arise as totalizations of spherical cubes
satisfying a Beck--Chevalley condition (Proposition \ref{prop:sphtot}). Remarkably, these cubical
diagrams can also be interpreted geometrically as perverse schobers on $\CC^n$ but with
stratification given by the coordinate hyperplanes along with their intersections and complements
(\S\ref{subsec:sphcub}). In this context, the totalization procedure can be regarded as a
higher--dimensional analog of Seidel's formation of the ``directed subcategory on the vanishing
cycles'' as a description of the Fukaya--Seidel category. We hope to explore the geometric relevance
of this in future work.

\subsection{Lax additivity}\label{intro:laxadd}

Plenty of thought has been invested in the question into what kind of categorical structure the
collection of derived categories should be organized, or, as Tamarkin put it:
``What do DG-categories form?'' \cite{tamarkin:dg}.
Several precise answers to this question have been given within the
various frameworks of enhanced triangulated categories that can be used to study derived categories,
such as differential graded categories \cite{toen:lectures}, $A_{\infty}$-categories
\cite{stasheff,keller:ainfty}, or stable $\infty$-categories \cite{lurie:ha}. 

Relevant for this work, in which we take the approach via stable $\infty$-categories, will be that
the \((\infty,2)\)-category of stable $\infty$-categories is {\em additive} in a suitable lax
$2$-categorical sense. We give an informal introduction to the circle of ideas surrounding the
concept of lax additivity in \S\ref{sub:marked_limits} and \S\ref{subsec:laxsum} adapted specifically
to the case of stable \(\infty\)-categories. In the companion paper \cite{CDW24}, we systematically develop the
general theory of lax additivity and lax matrices. The following table summarizes the main features
and should be seen as a natural continuation of the table of
\S\ref{subsec:rules-of-categorification}:
\begin{center}
\renewcommand{\arraystretch}{3}
  \begin{tabular}{lll}
    & classical & categorified\\ \hline
    7) & additive (\(\infty\)-)category \(\A\) & lax additive \((\infty,2)\)-category \(\mathbb{A}\)\\
    8) & hom-sets \(\A(X,Y)\) have addition & hom-categories \(\mathbb{A}(\X,\Y)\) have colimits\\
    9) & hom-sets \(\A(X,Y)\) are abelian groups & hom-categories \(\mathbb{A}(\X,\Y)\) are stable\\
    10) & \parbox{8cm}{finite direct sums\\\(\oplus_{s=1}^kx_s=\prod_{s=1}^kx_s=\coprod_{s=1}^kx_s\)} &
         \parbox{8cm}{general lax bilimits\\\(\laxbilim_{s:S}\X_s=\laxlim_{s:S}\X_s=\laxcolim_{s:S}\X_s\)}\\
    11) & \parbox{8cm}{matrices
         \(
         \begin{psmallmatrix}
           m_{11}& m_{12}\\
           m_{21}& m_{22}
         \end{psmallmatrix}
         \)\\ and matrix multiplication} &
               \parbox{8cm}{lax matrices
               \(
               \Donematrix{M_{11}}{M_{12}}{M_{21}}{M_{22}}
               \)\\ and lax matrix multiplication
               }
  \end{tabular}
\end{center}
\vspace{2ex}

Accepting our basic premise that abelian groups are to be categorified by stable
\(\infty\)-categories, Rule 9) requires no further comment.  Rule 8) is a convenient intermediate
step, categorifying the situation where the addition on hom-sets does not necessarily have inverses;
just like the uncategorified case, many basic lemmas are most naturally expressed in this generality
leading to the notion of lax \emph{semi}-additive \((\infty,2)\)-category.

The direct sum of abelian groups is both a categorical product and a categorical coproduct, a
universal property that is taken as the definition in general additive categories.  Rule 10) states
that the same definition can be categorified, as long as we replace finite products and coproducts,
i.e., limits and colimits indexed by finite discrete categories
\(S=\set{1,\dots,k}\),
with \emph{lax} limits and colimits indexed by arbitrary \(\infty\)-categories \(S\).  Apart from
this change, the theory is exactly analogous: if the hom-categories have colimits (categorifying
addition) then such lax limits and colimits always agree if they exits, yielding the notion of lax
bilimits.  When the hom-categories are stable, we can say even more for certain special shapes, such
as \(S=\Delta^1\): in this case all four possible universal 2-categorical constructions (lax/oplax,
limit/colimit) associated to an \(S\)-diagram \(\X\xrightarrow{F}\Y\) agree with each other.

A very convenient feature of additive categories is that maps \(m\colon x_1\oplus x_2\to
y_1\oplus y_2\) between direct sums can be represented as matrices of the form
\[
m = \begin{pmatrix}
  m_{11}\colon x_1\to y_1, & m_{12}\colon x_2\to y_1\\
  m_{21}\colon x_1\to y_2, & m_{22}\colon x_2\to y_2
\end{pmatrix}.
\]
Composing such maps then just amounts to the usual matrix multiplication.  Rule 11) states that these matrices and the matrix multiplication formula can be categorified, yielding an analogous theory of matrices
indexed in each coordinate not by a finite set but by arbitrary \(\infty\)-categories.  These
matrices are just a dependent version of bimodules, which by Morita theory encode functors between
module categories.

\subsection{Relations to previous work}

As already indicated in the very beginning, the idea to consider categorical complexes is, in
principle, by no means new: short exact categorical sequences in the form of localization sequences
already play a central role in algebraic $K$-theory. Also the directed pullback construction has
already been used in this context, cf.~\cite{LT19}, where (implicitly) also the product totalization
of a square is being considered.  As the authors have told us, they have also observed a type of
Beck--Chevalley phenomenon, that will be relevant for applications to $K$-theory. A version of
categorical complexes for non-stable symmetric monoidal $\infty$-categories was introduced in
\cite{Lur:tft}, but the motivation in this context is quite different from ours. 
The Calabi-Yau structures introduced in \S \ref{sec:CY} can be interpreted in the framework of
\cite{ST:cy} but our examples have not been considered there. 
The work in progress \cite{auroux:multi} studies Fukaya categories for Landau-Ginzburg models with
{\em multi-potentials}. As indicated to us by the authors, this theory will provide a means to
introduce cubical categorical diagrams in terms of partially wrapped Fukaya categories which one
would expect to mirror the cubical diagrams in \S \ref{subsec:completeintersection}. Variants of
this homological mirror symmetry conjecture for categorical cubes were already recorded in
\cite{Lee22}. The expected relation to our formulation will be explained below.

\subsection{Acknowledgements}

We are very grateful to Jeffrey Hicks, Sukjoo Lee, Nick Sheridan, and the participants of the Edinburgh
Hodge Seminar in general for interesting conversations about the topic of this work. In particular,
we thank them for making us aware of the work in progress on Fukaya categories of multi-potentials.
Many thanks to Denis Auroux, Sheel Ganatra, Andrew Hanlon, and Maxim Jeffs for taking the time and
effort to provide us with a report on their work in progress \cite{auroux:multi} (and to Nick
Sheridan for organizing the lecture series). 
We further thank Jacob Lurie for interesting remarks on categorical complexes. We also thank
Federico Barbacovi for discussions about spherical functors.
T.D.\ thanks M.\ Kapranov and V.\ Schechtman for many inspiring discussions on perverse sheaves and
schobers, in particular, this work draws substantial inspiration from their perspectives proposed in
\cite{KS14} and related subsequent work.
M.C.\ thanks the Hausdorff Research Institute for the hospitality during his stay, during which part of this work was written. M.C.\ was funded by the Deutsche Forschungsgemeinschaft (DFG,
German Research Foundation) under Germany's Excellence Strategy -- EXC-2047/1 -- 390685813. 
M.C.\ and T.D.\ acknowledge support by the Deutsche Forschungsgemeinschaft under Germany’s Excellence
Strategy -- EXC 2121 “Quantum Universe” -- 390833306. 
T.D.\ acknowledges support of the VolkswagenStiftung through the Lichtenberg Professorship Programme.
T.W.\ acknowledges support by the SFB 1085 -- Higher Invariants, funded by the DFG.

\section{Basic concepts}
\label{sec:basic}

\subsection{Stable \texorpdfstring{$\infty$}{infinity}-categories}
\label{subsec:basic_stable}

As a model for ``enhanced triangulated categories'', we will use stable $\infty$-categories as
introduced in \cite{lurie:ha} which, along with \cite{lurie:htt} and \cite{GR17}, are our standard
references. It will be most convenient to work with presentable stable $\infty$-categories -- the
$\infty$-category formed by these with colimit--preserving functors will be denoted by $\St$. 

For applications, it is also important to consider presentable stable $\infty$-categories linear
over a given field, or more generally over a commutative ring spectrum, denoted $k$. Following \cite{lurie:ha},
we define these as modules in $\St$ over the symmetric monoidal stable $\infty$-category $\Modk$ of
$k$-module spectra. The resulting $\infty$-category of {\em $k$-linear stable $\infty$-categories}
will be denoted by $\Stk$. Given objects $\C,\D \in \Stk$, we denote by $\Stk(\C,\D)$ the
$\infty$-category of morphisms in $\Stk$, i.e., colimit--preserving $k$-linear functors. The
$\infty$-category $\Stk(\C,\D)$ can be regarded as an object of $\Stk$ and thus defines an {\em
internal hom}. Further, $\Stk$ admits a symmetric monoidal structure (adjoint to the internal hom)
and we denote the corresponding binary tensor product by $\otimes$.

Any compactly generated $\infty$-category $\C \in \Stk$ is dualizable with respect to this monoidal
structure and we denote its dual by $\C^{\vee}$. For the subcategory $\C_0 \subset \C$ of compact
objects, the $\infty$-category $\C$ is equivalent to the $\infty$-category $\Ind(\C_0)$ of
$\Ind$-objects of $\C_0$, and we have
\[
	\C^{\vee} \simeq \Stk(\C, \Modk) \simeq \Ind(\C_0^{\op}).
\]
For a morphism $F: \C \to \D$ in $\Stk$ between compactly generated $\infty$-categories, we denote its
dual by
\[
	F^{\vee}: \D^{\vee} \lra \C^{\vee}.
\]

By the adjoint functor theorem, any functor $F$ in $\Stk$ admits a right adjoint, denoted by
$\ra{F}$.  The adjoint $\ra{F}$ is not necessarily a morphism in $\Stk$, since it may not preserve
colimits. The condition that $\ra{F}$ preserves colimits and admits a $k$-linear
structure is equivalent to the condition that $F$ preserves compact objects. When we speak of a $k$-linear adjoint, we mean that the adjoint also lies in $\Stk$, and the adjunction is in this case given in a $2$-categorical sense in $\Stk$. 

\subsection{Complexes of stable \texorpdfstring{$\infty$}{infinity}-categories}\label{subsec:complexesofstables}

We define a {\em chain complex of stable $\infty$-categories} to consist of a sequence
\[
	\begin{tikzcd}
		\dots \ar{r}{d} &  \A_2 \ar{r}{d} & \A_{1} \ar{r}{d} & \A_0 \ar{r}{d} &   \dots 
	\end{tikzcd}
\]
of morphisms $d\colon \A_i \to \A_{i-1}$ in $\Stk$ satisfying $d^2 \simeq 0$. Note, that the existence
of an equivalence $d^2 \simeq 0$ simply amounts to the statement that $d^2$ is a zero object in the
stable $\infty$-category $\Stk(\A_i, \A_{i-2})$ and the space of zero objects is contractible.
Thus, in order to specify a complex of stable $\infty$-categories, one only needs to provide the
list of differentials, and verify the condition $d^2 \simeq 0$. In particular, in contrast to the
notion of a complex of {\em objects in a fixed stable $\infty$-category}, no higher coherence data
is involved. We justify this a bit more formally: 

\begin{defi} Following \cite{tashi:coherent}, we define the category
\[
	\Ch \coloneqq \ZZ_* / \sim
\]
where $\ZZ_{*} $ denotes the category obtained from $\ZZ$ by adjoining a zero object $*$ and $\Ch$
is the quotient obtained by identifying each composite $d^2\colon i \to i+2$ with the $0$ morphism (i.e.
the unique composite $i \to * \to i+2$). We then define a {\em coherent chain complex} in a
pointed $\infty$-category $\C$ to be a functor $\Ch^{\op} \to \C$ that preserve zero objects.
\end{defi}

Furthermore, we denote by
\begin{equation}
  G\coloneqq \dots \amalg_{\set{-1}}\set{-1\to 0}
  \amalg_{\set{0}}\set{0\to 1} \amalg_{\set{1}}\set{1\to2}\amalg_{\set{2}}\dots
  \subset \ZZ
\end{equation}
the \(1\)-skeleton of (the nerve of) \(\ZZ\).
By definition, a map of simplicial sets \(G\to K\)
just consists of a sequence of composable edges in the simplicial set \(K\).

\begin{lem}\label{lem:coherent} We regard $\Stk$ as an $\infty$-category (by discarding noninvertible
	$2$-morphisms) and consider further the $\infty$-categories
	\begin{itemize}
		\item $\Fun_*(\Ch^{\op},\Stk)$ of functors preserving zero objects, 
		\item $\Fun_{(d^2)}(\ZZ^{\op},\Stk)$ of functors such that each $d^2$ is equivalent to the
			$0$ functor,
		\item $\Fun_{(d^2)}(G^\op,\Stk)$ of functors such that each $d^2$ is equivalent to the
			$0$ functor.
	\end{itemize}
	The restriction functors
	\[
		\Fun_*(\Ch^{\op},\Stk) \lra \Fun_{(d^2)}(\ZZ^{\op},\Stk) \lra \Fun_{(d^2)}(G^\op,\Stk)
	\]
	are equivalences of $\infty$-categories.
\end{lem}
\begin{proof}[Proof sketch]
  The rightmost restriction functor is an equivalence
  because the inclusion \(G\subset \ZZ\) is an inner anodyne extension
  of simplicial sets.

  Extending a functor \(F\colon \ZZ^\op\to\Stk\) to \(\Ch^\op\)
  involves various lifting problems in the hom-spaces \(\Stk(F(i),F(j))\).
  If we assume that for each \(k\geq 2\) the edge \(F(d^k)\) is equivalent to zero,
  then all of these lifting problems take place in the zero component
  of the mapping space \(\Stk(F(i),F(j))\).
  Hence the statement follows from the key fact that for each \(\A,\B\) in \(\Stk\),
  the zero connected component of \(\Stk(\A,\B)\)
  is just the groupoid of zero functors \(\A\to\B\), which is contractible.
  We omit the remaining details.
\end{proof}
	
In virtue of Lemma \ref{lem:coherent}, we may always extend a given complex in $\Stk$ presented as a
sequence of differentials $d$ satisfying $d^2 \simeq 0$, to a fully coherent complex in an
essentially unique way. Further, we will sometimes want to turn a coherent complex concentrated in degrees $n \ge i \ge
0$ into a cubical diagram: To this end, we may define the functor
\[
	q\colon \P(\{1,...,n\}) \lra \Ch, M \mapsto \begin{cases} 0 & \text{if $M = \emptyset$,}\\
								i & \text{if $M = \{1,...,i\}$,}\\
								* & \text{else}
							\end{cases}
\]
and consider the pullback
\begin{equation}\label{eq:coherent_cube}
	q^*\colon \Fun_*(\Ch^{\op}, \Stk) \lra \Fun_*(\P(\{1,...,n\})^{\op}, \Stk).
\end{equation}

For the sake of brevity, we will refer to a complex in $\St_k$ as a {\em categorical complex}.

\begin{defi}
We denote by $\Ch(\St_k)$ the full subcategory of $\Fun(\ZZ^\op,\St_k)$ spanned by categorical
complexes.
\end{defi}

An categorical $n$-complex is a functor $\A_{*}\colon \ZZ^{\op,n}\rightarrow \St_k$, meaning
collections of commuting $n$ differentials
\[ d_i\colon\A_{a_1,\dots,a_i,\dots,a_n}\longrightarrow \A_{a_1,\dots,a_i-1,\dots,a_n}\,,\]
satisfying that $d_i^2\simeq 0$ for all $1\leq i\leq n$. 

\begin{defi}\label{def:ncplx}
Let $n\geq 1$.
\begin{enumerate}[(1)]
\item We denote by $\Ch_n(\St_k)$ the full subcategory of $\Fun((\ZZ^{\op})^n,\St_k)$ spanned by categorical $n$-complexes. 
\item A categorical $n$-complex $\A_*$ is called a categorical $n$-cube if $\A_{{\bf j}}\simeq 0$
for all ${\bf j}\in (\ZZ^{\op})^n$ and ${\bf j} \notin \I^n$, where $\I=[1]^{\op}=\{1\to 0\}\subset
	\ZZ^{\op}$. We denote by $\Cube_n(\St_k)\subset \Ch_n(\St_k)$ the full subcategory spanned by categorical $n$-cubes.
\end{enumerate}
\end{defi}

\subsection{Partially lax limits}%
\label{sub:marked_limits}

There are different types of limits of a given diagram $\{\A_i\}_{i : I}$ in an $(\infty,2)$-category. Namely, for each
triangle 
\[
  \begin{tikzcd}[column sep=small]
    & \ar[dr, ""{name=U}]\ar[dl] \lim   & \\
    \A_i \ar[rr] & & \A_j
  \end{tikzcd}
\]
in the limit cone, one can require either of the following:
\begin{enumerate}[label=(\roman*)]
  \item {\em strict}: it commutes up to natural equivalence:
	\[
		  \begin{tikzcd}[column sep=small]
		    & \ar[dr, ""{name=U}]\ar[dl] \lim   & \\
		    \A_i \ar[rr]\ar[to=U, shorten <=20pt, shorten >=10pt, equal] & & \A_j
		  \end{tikzcd}
	  \]
  \item {\em lax}: it commutes up to a possibly noninvertible $2$-morphism as in 
	  \[
		  \begin{tikzcd}[column sep=small]
		    & \ar[dr, ""{name=U}]\ar[dl] \lim   & \\
		    \A_i \ar[rr]\ar[to=U, shorten <=20pt, shorten >=10pt, Rightarrow] & & \A_j
		  \end{tikzcd}
	  \]
  \item {\em oplax}: it commutes up to a possibly noninvertible $2$-morphism as in 
	  \[
		  \begin{tikzcd}[column sep=small]
		    & \ar[dr, ""{name=U}]\ar[dl] \lim   & \\
		    \A_i \ar[rr]\ar[to=U, shorten <=20pt, shorten >=10pt, Leftarrow] & & \A_j
		  \end{tikzcd}
	  \]
\end{enumerate}
The oplax triangle may be modelled as an amalgamate of a strict and a lax triangle as in 
	  \[
		  \begin{tikzcd}[column sep=small]
		    & & \ar[dll]\ar[d, ""{name=U}] \lim \ar[drr, ""{name=V}] &  & \\
		    \A_i \ar[rr]\ar[to=U, shorten <=20pt, shorten >=10pt, equal] & & \A_j & &
		    \ar[ll,"\id"] \ar[to=U, shorten <=15pt, shorten >=5pt, Rightarrow] \A_j
		  \end{tikzcd}
	  \]
so that it suffices to distinguish between strict and lax triangles. One way to encode this ``partial
laxness'' is to simply mark those edges in the diagram $\{\A_i\}_{i:I}$ over which we require the
cone to strictly commute. The resulting notion of {\em marked limit} has been developed and studied
by several groups of authors, see, e.g., \cite{abellan-stern:cofinality} for a formal definition.

Throughout this work, we will use partially lax universal constructions to describe stable $\infty$-categories and functors among them (see \S \ref{sub:marked_limits}). In \S \ref{sec:dk} through \S
\ref{sec:examples}, the focus lies on concreteness and explicit examples, so that we take a
hands-on approach and use explicit models for the relevant constructions. The framework of
$\infty$-categories is sufficient to describe the needed universal properties, even though we often
formulate them in $(\infty,2)$-categorical terms. We also wish to point out that we do not make use of a specific model for $(\infty,2)$-categories, but only use typical features of any such theory.

In \cite{CDW24}, we discuss (some of) these universal constructions from a more foundational
point of view within a framework of {\em lax additive $(\infty,2)$-categories} (also see \S
\ref{subsec:laxsum}). 

\subsection{Lax additivity}
\label{subsec:laxsum}

For given abelian groups $A$ and $B$, their product and coproduct are characterized by the universal
cones
\begin{equation}\label{eq:universal_sum}
	\begin{tikzcd}
	&A \times B\ar{d}\ar{r} &  B & 			      	 & A \amalg B  & B  \ar{l}         \\
	&     A  &                  & \quad\text{and}\quad            & 	A \ar{u}, & 
	\end{tikzcd}
\end{equation}
respectively. The fact that these constructions agree up to canonical isomorphism is typically
indicated by using the terminology {\em direct sum}, and the notation $A \oplus B$. This phenomenon
is referred to as the {\em semiadditivity} of the category of abelian groups. {\em Additivity} then
amounts to the extra condition of the resulting abelian monoid structure on $\Hom(A,B)$ being a
group.

As already indicated in \S \ref{subsec:rules-of-categorification}, the analogous categorified
universal constructions typically involve more data, invisible upon passing to $K_0$. 
Namely, beyond a pair of stable $\infty$-categories $\A$, $\B$, we are also given an exact functor
$F\colon \A \to \B$. The relevant universal constructions that we may build from this data are
characterized by the universal cones depicted in Figure \ref{fig:laxcones}.
\begin{figure}[h]
	\begin{center}
\begin{tabular}{c|c}
	lax limit: 
  \laxcone{\A}{\laxlima{\A}{F}{\B}}{B}{}{}{F}
	&oplax limit: 
    \oplaxcone[]{\A}{\oplaxlima{\A}{F}{\B}}{B}{}{}{F}
	\\[1ex]\hline\\
	lax colimit: 
  \laxcocone{\A}{\laxcolima{\A}{F}{\B}}{\B}{}{}{F}
	&
	oplax colimit: 
    \oplaxcocone{\A}{\oplaxcolima{\A}{F}{\B}}{\B}{}{}{F}
\end{tabular}
\caption{The four lax cones with base given by a functor $F\colon \A \to \B$.}
\label{fig:laxcones}
\end{center}
\end{figure}

\begin{thm}[Lax Additivity, \cite{CDW24}] In the $(\infty,2)$-category $\Stk$ of presentable stable
	$\infty$-categories, the four lax universal constructions from Figure \ref{fig:laxcones}
	exist and are canonically equivalent. 
\end{thm}

As explained in \cite{CDW24}, the analogous statement holds in any {\em lax additive}
$(\infty,2)$-category and $\St_k$ is an example of such. In fact, \cite{CDW24} provides a more
refined systematic analysis also introducing the notion of a {\em lax semiadditive}
$(\infty,2)$-category as an intermediate step. 

From an axiomatic perspective, the canonical equivalence of these four universal constructions
captures the essence of what seems to make lax additive $(\infty,2)$-categories a suitable context
for categorified homological algebra (just like additive categories are used for classical
homological algebra).

While complete proofs can be found in \cite{CDW24}, we explain here, how to construct a
particular model for the lax limit in $\St$ and describe the four universal lax cones from Figure
\ref{fig:laxcones} in terms of this model (the analogous statement for $\Stk$ is somewhat more
involved as we need to keep track of the $\Modk$-module structures -- we don't discuss this here).
Let $\Sigma(F)$ denote the $\infty$-category of sections of the covariant Grothendieck construction
of $F$, considered as a diagram $\Delta^1 \to \St$. Due to the simplicity of the indexing category,
this can be described even more concretely as the pullback of simplicial sets
\[
	\begin{tikzcd}
		\Sigma(F) \ar[r]\ar[d] &  \Fun(\Delta^1, \B) \ar[d, "\ev_0"]\\
		\A  \ar[r, "F"]	& \B
	\end{tikzcd}
\]
Notationally, we simply write
\[
	\Sigma(F) = \{ (a,b,\eta)\; |\; \text{$a \in \A$, $b \in \B$, $\eta\colon F(a) \to b$} \}.
\]
We now describe the universal lax cones that characterize $\Sigma(F)$ as each of the four universal
constructions:
\begin{enumerate}[label=\arabic*.]
	\item The {\em lax limit cone} is evident: The functors to $\A$ and $\B$, respectively,
are simply given by projecting to the components $a \in \A$ and $b \in \B$ while the natural
transformation which is part of the cone is given by $F(a) \to b$. Note that, while we specify the
data of the cone on the level of objects, it is evident how to extend these formulas to actual
functors and natural transformations of $\infty$-categories. 

\item The {\em lax colimit cone} of $\Sigma(F)$ corresponds to the functors
\begin{align*}
	\A \to \Sigma(F),\quad & a \mapsto (a, F(a), F(a) \overset{\id}{\to} F(a))\\
	\B \to \Sigma(F),\quad & b \mapsto (0, b, 0 \to b)
\end{align*}
where the natural transformation assigns to $a \in \A$, the morphism in $\Sigma(F)$ given by
\[
	\begin{tikzcd}
		( & 0 \ar[d], & F(a) \ar[d,"\id"], &  0 \ar[d] \ar[r] & F(a) \ar[d,"\id"] & )\\
		( & a, & F(a), &  F(a) \ar[r,"\id"] & F(a) & )
	\end{tikzcd}
\]

\item The {\em oplax limit cone} is given by the functors
\begin{align*}
	\Sigma(F) \to \A,\quad & (a,b,\eta) \mapsto a,\\
	\Sigma(F) \to \B,\quad & (a,b,\eta) \mapsto \fib(\eta),
\end{align*}
and the natural transformation which assigns to $(a,b,\eta) \in \Sigma$, the morphism in $\B$ given by
\[
	\fib(\eta) \to F(a)
\]
which is part of the fiber square
\[
	\begin{tikzcd}
		\fib(\eta) \ar[r]\ar[d] & F(a) \ar[d]\\
		0 \ar[r] & b.
	\end{tikzcd}
\]

\item Finally, the {\em oplax colimit cone} is determined by the functors
\begin{align*}
	\A \to \Sigma(F),\quad & a \mapsto (a, 0, F(a) \to 0)\\
	\B \to \Sigma(F),\quad & b \mapsto (0, b[1], 0 \to b[1])
\end{align*}
while the natural transformation assigns to $a \in \A$, the morphism in $\Sigma(F)$ given by
\[
	\begin{tikzcd}
		( & a \ar[d], & 0 \ar[d], &  F(a) \ar[d] \ar[r] & 0 \ar[d] & )\\
		( & 0, & F(a)[1], &  0 \ar[r] & F(a)[1] & )
	\end{tikzcd}
\]
where the square is the biCartesian square exhibiting $F(a)[1]$ as the suspension of $F(a)$. 
\end{enumerate}

\begin{rem}\label{rem:semiadditive} Note that, for $F = \id\colon \A \to \A$, we have 
	\[
		\Sigma(F) = \Fun(\Delta^1, \A)
	\]
	so that $\Sigma(F)$ is simply a presheaf category, or put differently, the category of
	$\A$-valued representations of the quiver of Dynkin type $\AA_2$. In this situation, the fact that
	$\Sigma(F)$, which most evidently is a lax limit of $F$, is in fact also a lax colimit of
	$F$ amounts to the statement that any presheaf category is generated by the subcategory of
	representable presheaves (tensored with objects of $\A$). As explained in more detail in \cite{CDW24}, the equivalence (for general $F$)
	\[
		\laxlima{A}{F}{B} \simeq \laxcolima{A}{F}{B}
	\]
	generalizes to any so-called {\em lax semiadditive} $(\infty,2)$-category. In particular, it
	does not require stability and also holds, for example, in the $(\infty,2)$-category $\PrL$ of
	presentable $\infty$-categories.
\end{rem}

\begin{rem}\label{rem:additive}  In contrast to the equivalence between the lax limit and colimit
	discussed in Remark \ref{rem:semiadditive}, the phenomenon that the lax and oplax limit are
	equivalent does {\em not} hold in a general lax semiadditive $(\infty,2)$-category but requires
	{\em lax additivity}. For example, the $(\infty,2)$-category $\St$ of {\em stable} presentable
	$\infty$-categories has this property. 

	The resulting equivalence
	\begin{equation}\label{eq:lax_oplax}
			\laxlima{A}{F}{B} \simeq \oplaxlima{A}{F}{B}
	\end{equation}
	is exploited in \cite{DJW19-BGP} where it is described explicitly, in terms of models for lax
	and oplax limits, and used to provide a framework for generalized
	Bernstein-Gelfand-Ponomarev reflection functors. 

	This equivalence also has a natural interpretation in terms of {\em semiorthogonal
	decompositions}, as we now explain. Recall (cf. \cite{DKSS:spherical}) that a semiorthogonal
	decomposition (of length $2$) of a stable $\infty$-category $\C$ consists of a pair $(\X,\Y)$ of full stable
	subcategories of $\C$, such that the functor of $\infty$-categories
	\[
		\{ x \to y \; |\; \text{$x : \X$, $y : \Y$} \} \lra \C,
	\]
	given by associating to a morphism $x \to y$ its cofiber, is an equivalence. The lax and
	oplax limit cones for $\Sigma(F)$ naturally induce semiorthogonal decompositions with
	components given by the kernels of the projection maps. This yields concretely:
	\begin{enumerate}[label=\arabic*.]
		\item {\em lax limit of }$F$: $(\{(a,0,F(a) \to 0)\}, \{(0, b, 0 \to b)\})$
		\item {\em oplax limit of }$F$: $(\{(0,b, 0 \to b)\},\{(a, b, F(a) \overset{\simeq}{\to}b)\})$
	\end{enumerate}	
	In terms of the terminology developed in \cite{DKSS:spherical}, these two semiorthogonal
	decompositions of $\Sigma(F)$ are related by {\em mutation}, where the first decomposition is
	{\em coCartesian} and the second decomposition {\em Cartesian}. Note that conversely any Cartesian or coCartesian semiorthogonal decomposition with gluing functor $F$ is equivalent to the above described corresponding semiorthogonal decomposition of $\Sigma(F)$. If $F$ further admits a right adjoint $F^R$, then the first semiorthogonal decomposition is also Cartesian, which is reflected by the fact that the lax limit of $F$ is equivalent to the oplax limit of $F^R$, see \cite[Cor.~A.5]{CDW24}.
	
This perspective on length $n=2$ semiorthogonal decompositions with a gluing functor as an (op)lax limit over $\Delta^1$ can be extended to semiorthogonal decompositions of arbitrary length $n$ using (op)lax limits over $2$-categorical versions of $\Delta^{n-1}$. A semiorthogonal decompositions of length $n$ of a stable $\infty$-category $\C$ can equivalently be expressed as $n-1$ iterated semiorthogonal decompositions of length $2$ of an increasing sequence of stable subcategories of $\C$. If these semiorthogonal decompositions of length $2$ are each coCartesian, then $\C$ arises as a lax limit over diagram indexed by the $2$-categorical simplex $\Delta^{n-1}$. In this setting, the single gluing functor from the case $n=2$ is thus replaced by a collection of gluing functors together with gluing natural transformations between composites of the gluing functors. 
\end{rem}

When we wish to make the point that the above universal constructions are equivalent (where the
equivalences are determined by the choices of universal cones in Figure \ref{fig:laxcones}), then we
will use the notation
\[
	\laxsum{\A}{F}{\B}  
\]
to denote any of them, and refer to it as the {\em lax sum} of $\A$ and $\B$ along $F$. Whenever we
would like to refer to one of the above universal lax cones, then we will use the 
notation for the respective universal construction from Figure \ref{fig:laxcones}.


\section{Simplicial totalization and the Dold-Kan correspondence}
\label{sec:dk}

\subsection{The categorical cochain complex of a simplex}
\label{subsec:warmup}

To convey a first impression of the workings of categorical complexes, we present a class of
examples which is easy to describe directly yet illustrates some of the key phenomena.

Let $A$ be an abelian group. We will describe categorifications of the 
simplicial cochain complex $C^{\bullet}(\Delta^n, A)$ of an $n$-simplex. We have 
\[
	C^{k}(\Delta^n, A) = \Map( \{\sigma\colon [k] \hra [n]\},  A)
\]
with differential given by the formula
\[
	(d a)(\sigma) = \sum_{i=0}^n (-1)^i a(\sigma \circ \partial_i).
\]
Let us further focus on the $2$-simplex, where we have 
\[
	\begin{tikzcd}[sep=3em]
		& C^{\bullet}(\Delta^2, A): & & \Cc^{\bullet}(\DDelta^2, \A):    \\
		A  \cong & \{(x_{012})\} & & \{X_{012}\}  & \simeq \A\\
	A^3  \cong & \{(x_{01}, x_{02}, x_{12})\} \ar[swap]{u}{x_{012} = x_{12} - x_{02} + x_{01}}&
		   & \{X_{01} \to  X_{02} \to X_{12}\} \ar[swap]{u}{X_{012} = \tot(X_{01} \to X_{02} \to X_{12})} & \simeq \langle \A,\A,\A \rangle\\
	A^3  \cong &  \{(x_0, x_1, x_2)\} \ar[swap]{u}{x_{ij} = x_j - x_i} &  & \{X_0 \to X_1 \to
	X_2 \} \ar[swap]{u}{X_{ij} = \cone(X_i \to X_j)} & \simeq \langle \A,\A,\A \rangle
	\end{tikzcd}
\]
We note the following features of the categorical complex $\Cc^{\bullet}(\DDelta^2, \A)$:
\begin{enumerate}
	\item To be able to define the first differential, we have adjoined morphisms $X_i \to X_j$ so that we can take
		their ``difference'' according to rule 3).
	\item To define the second differential, the objects $X_{ij}$ need to form a complex so that
		we may totalize it.
	\item The octahedral ``axiom'' (or rather the third isomorphism theorem) guarantees that the
		first differential is well defined (the cones $\cone(X_i \to X_j)$ form a complex)
		and further, that $d^2 \simeq 0$ (since the just-mentioned complex is exact). 
\end{enumerate}

The additional lax data we needed to implement the rules of categorification from \S \ref{subsec:rules-of-categorification}
actually has a natural interpretation: The simplex category $\Delta$ can be defined as the full
subcategory of the $1$-category $\Cat$ of small categories spanned by the standard ordinals. As a result,
given ordinals $[k]$ and $[n]$ we obtain a set $\Delta([k],[n])$ of $k$-simplices in $\Delta^n$
described as the corresponding set of maps in $\Delta$. However, the collection of categories has a
natural $2$-categorical structure taking in to account natural transformations. Defining $\DDelta$
to be the full $2$-subcategory of the $2$-category $\CCat$ instead yields morphism {\em categories},
in fact posets, $\DDelta([k],[n])$.

\begin{defi}\label{defi:catcochain}
Within this context, we can now give the general definition of the complex of the {\em categorical
simplicial cochain complex} $\Cc^{\bullet} \coloneqq  \Cc^{\bullet}(\DDelta^n, \A)$:
\begin{enumerate}
	\item \label{catcochain:1} For $0 \le k \le n$, the $\infty$-category $\Cc^k$ is defined to be the full
		subcategory
		\[
			\Cc^k \subset \Fun(\DDelta([k],[n]), \A)
		\]
		consisting of those diagrams $X$ satisfying:
		\begin{itemize}
			\item for every non-injective map $\tau\colon [k] \to [n]$, we have
				$X(\tau) \simeq 0$. 
		\end{itemize}
	\item \label{catcochain:2} The differential $d\colon \Cc^k \to \Cc^{k+1}$ is given by associating to a diagram $X \in
		\Cc^k$, the diagram
		\[
			dX\colon \DDelta([k+1],[n]) \to \A,\; \sigma \mapsto \tot(X(d_n \sigma) \to
			... \to X(d_0\sigma)).
		\]
\end{enumerate}
\end{defi}

\begin{rem}
	\label{rem:sloppyconvention}
	We comment on the defining formula for the diagram $dX$ in \ref{catcochain:2} which needs to
	be made precise in several ways. 
	\begin{enumerate}[label=\arabic*.]
		\item \label{d:1} First, the formula
		\[
			\tot(X(d_n \sigma) \to ... \to X(d_0\sigma))
		\]
		needs to be explained. The expression 
		\[
			X(d_n \sigma) \to ... \to X(d_0\sigma)
		\]
		refers to an $n+1$-term complex in the stable $\infty$-category $\A$. Formally, this
		is realized as a cubical diagram $\I^{n} \to \A$, with $\I=[1]^\op=\{1\to 0\}$, with zero relations encoded by
		certain vertices of this cube being mapped to zero objects. The (cofiber) totalization of this complex is then obtained by extending with zeros to a punctured $(n+1)$-cube $\I^{n+1}\backslash \{(0,\dots,0)\}$, i.e.~a right Kan extension, and taking the colimit over $\I^{n+1}\backslash \{(0,\dots,0)\}$, i.e.~a left Kan extension to $\I^{n+1}$. 
	\item \label{d:2} Second, we need to explain how the various values $dX(\sigma)$ organize
		into an actual diagram $\DDelta([k+1],[n]) \to  \A$ and, further, how the association $X \mapsto dX$ defines
			a functor in $\Stk$. This is achieved by means of the formalism of Kan extensions as developed in \cite{lurie:htt}
	\end{enumerate}

	For most parts of this paper, we will not delve into the technical details of constructions
	such as \ref{d:1} and \ref{d:2} as this would have a tendency of hiding the main ideas
	behind the (routine) technical details. Rather, we assume that the reader is familiar with
	the techniques alluded to and keep the treatment somewhat informal as in Definition
	\ref{defi:catcochain}.  Typically, it is rather straightforward how to make things formally
	precise so that we hope that not much is lost by this approach.
\end{rem}

\subsection{The categorified Dold--Kan correspondence}
\label{subsec:catdk}

The discussion in $\S$ \ref{subsec:warmup} shows that interesting examples of categorical
complexes arise when studying simplicial (or rather $2$-simplicial) objects. These ideas can be put
in a more systematic context, the {\em categorified Dold--Kan correspondence} as established in
\cite{dyck:dk}, which will be briefly recalled here.

Denote by $\CCat$ the $2$-category of small categories, and further, by $\DDelta \subset \CCat$ the
full sub $2$-category spanned by the standard ordinals $\{[n] | n \ge 0 \}$. Treating $\DDelta$ as
an $(\infty,2)$-category, we define a $2$-stable $\infty$-category to be a functor
\[
	\X\colon \DDelta^{\op} \to \St
\]
of $(\infty,2)$-categories. For example, this can be modelled concretely as a $\sSet$-enriched
functor from $\DDelta$, considered as $\sSet$-enriched by taking the nerve of the categories of
morphisms, to $\sSet$ taking values in stable $\infty$-categories with exact functors.  

We explain how to associate to $\X$ a complex of stable $\infty$-categories by a construction which
is, in a sense, dual to the one discussed in \S \ref{subsec:warmup}. 
To this end, for $n > 0$, consider the cubical diagram
\[
	q\colon \{0,1\}^n \lra \Fun([n-1],[n]), (i_0, ..., i_{n-1}) \mapsto 0 + i_0 \le 1 + i_1 \le ...
	\le (n-1) + i_{n-1}.
\]
Note that $\{0,1\}=[1]$ as posets, we use the former notation in this section to avoid confusions. The $2$-functoriality of $\X$ induces a cubical diagram 
\[
	q_{\X}\colon \{0,1\}^n \lra \Fun(\X_n,\X_{n-1})\,.
\]
We now denote by 
\[
	\overline{\X_n} \coloneqq  \X_n / \{\text{degenerate simplices}\}
\]
the ``Verdier quotient'' by the subcategory of degenerate simplices, i.e., the full subcategory
spanned by the images of all degeneracy functors. The above cube $q_{\X}$ descends to define a
cubical diagram
\[
	\overline{q_{\X}}\colon \{0,1\}^n \lra \Fun(\X_n,\overline{\X}_{n-1})
\]
which represents a homotopy-coherent $(n+1)$-term complex formed by the face maps
\[
	d_n \to d_{n-1} \to ... \to d_0. 
\]
Passing to the totalization of this complex, we obtain a functor
\[
	\tot(\overline{q_{\X}})\colon \X_n \to \overline{\X}_{n-1}
\]
and further, noting that this totalization maps degenerate simplices to zero objects, an induced
functor
\begin{equation} \label{eq:simplicial_differential}
	d\colon \overline{\X_n} \to \overline{\X}_{n-1}
\end{equation}

We leave the proof of the following lemma as an entertaining exercise to the reader:

\begin{lem}
	The differential constructed in \eqref{eq:simplicial_differential} satisfies $d \circ d
	\simeq 0$. 
\end{lem}

\begin{defi} The resulting complex $(\overline{\X_{\bullet}}, d)$ is called the {\em simplicial
	totalization of $\X$}. 
\end{defi}

As shown in \cite{dyck:dk}, the association $\X \mapsto (\overline{\X_{\bullet}}, d)$ defines an
equivalence which can be interpreted as a categorified variant of the classical Dold-Kan
correspondence:

\begin{thm}\label{thm:dk}
	The simplicial totalization functor defines an equivalence of $\infty$-categories
	\[
		\C\colon \St_{\DDelta} \llra \Ch_{\ge 0}(\St)\noloc \Nc
	\]
\end{thm}

Note that in \cite{dyck:dk} a different description of the functor $\C$ is given, but this can be
shown to be equivalent to the above, by investigating the fully faithful adjoints of the localization
functor $\X_n \to \overline{\X_n}$.  

One interesting aspect of the categorified Dold-Kan correspondence is the reverse procedure of
constructing  $2$-simplicial objects from categorical complexes. In the classical context, this
provides a bridge between homological and homotopical data, and we expect the categorified Dold-Kan
correspondence to play a similar role for categorical complexes.

\section{Totalizations of categorical multicomplexes}
\label{sec:tot}

Besides the totalization of simplicial objects, another natural class of examples of categorical
complexes arises by totalizing bicomplexes, or more generally multi-dimensional complexes. In
this section, we discuss these totalization constructions and establish some basic properties. This
will allow us to introduce and investigate several interesting classes of categorical complexes. 

\subsection{Directed pushouts and pullbacks}
\label{subsec:directed}

We introduce two universal $(\infty,2)$-categorical constructions which will be crucial in this section. 
\begin{enumerate}
	\item
The {\em \dpb{}} of a diagram
\begin{equation}
	\label{eq:cospan}
	\begin{tikzcd}
		& \B \ar[d, "S"] \\
	\C \ar[r, swap, "G"] & \D 
	\end{tikzcd}
\end{equation}
in $\Stk$ is an object $\laxpull{\B}{\D}{\C}$ together with a universal cone of the form
\begin{equation}
  \oplaxsquare
  {\laxpull{\B}{\D}{\C}}{\B}{\C}{\D}
  {}{}{S}{G}{}
\end{equation}
It may be expressed in terms of the lax limit from \S \ref{subsec:laxsum} as
\begin{equation}\label{eq:laxpullvialimit}
	\laxpull{\B}{\D}{\C} \simeq (\laxlima{\B}{S}{\D}) \times_{\D} \C
\end{equation}
or explicitly as the fiber product of simplicial sets
\[
	\laxpull{\B}{\D}{\C} = \B \times_{\D} \Fun(\Delta^1, \D) \times_{\D} \C = \{(b,c,\eta)\; | \;
	\text{$b \in \B$, $c \in \C$, $\eta\colon S(b) \to G(c)$} \}.
\]

\item
Dually, the {\em \dpo{}} of a diagram
\begin{equation}
	\label{eq:span}
	\begin{tikzcd}
	\A \ar[d, swap, "R"] \ar[r, "F"] & \B \\
	\C  & 
	\end{tikzcd}
\end{equation}
in $\Stk$ is an object $\laxpush{\B}{\A}{\C}$ together with a universal cone of the form
\begin{equation}
  \label{eq:universalpush}
  \oplaxsquare{\A}{\B}{\C}{\laxpush{\B}{\A}{\C}}{F}{R}{}{}{}
\end{equation}
We may construct it in terms of the oplax colimit via the formula
\[
	\laxpush{\B}{\A}{\C} = \B \amalg_{\A} (\oplaxcolima{\A}{R}{\C})
\]
or as a pushout 
\[
	\B \amalg_{\A} \Fun(\Delta^1,\A) \amalg_{\A} \C
\]
in $\Stk$, where the functors $\A \to \Fun(\Delta^1,\A)$ are given by the formulas
\[
	a \mapsto (0 \to a) \quad \text{and} \quad a \mapsto (a \overset{\id}{\to} a),
\]
respectively.
\end{enumerate}

In the presence of suitable adjoints, we may express the directed pullbacks (resp. pushouts) introduced
in this section in terms of the (op)lax (co)limits of \S \ref{subsec:laxsum}. This will be most
relevant for \S \ref{subsec:beck-chevalley} where we investigate Beck--Chevalley conditions and is
documented in the following proposition.

\begin{prop} 
	\begin{enumerate}
		\item Suppose we are given a diagram as in \eqref{eq:cospan}.
			\begin{enumerate}
			\item There is a canonical equivalence
					\[
						\laxpull{\B}{\D}{\C} \simeq \oplaxlima{\C}{\ra{S} \circ G}{\B}.
					\]
				\item  Suppose further that $G$ has a left adjoint $\la{G}$. Then there is a canonical equivalence
					\[
						\laxpull{\B}{\D}{\C} \simeq \laxlima{\B}{\la{G} \circ S}{\C}.
					\]
				
			\end{enumerate}
		\item Suppose we are given a diagram as in \eqref{eq:span}.
			\begin{enumerate}
			\item Then there is a canonical equivalence
					\[
						\laxpush{\B}{\A}{\C} \simeq \laxcolima{\C}{F \circ \ra{R}}{\B}.
					\]
				\item Suppose further that $F$ has a left adjoint $\la{F}$. Then there is a canonical equivalence
					\[
						\laxpush{\B}{\A}{\C} \simeq \oplaxcolima{\B}{R \circ
						\la{F}}{\C}.
					\]
				
			\end{enumerate}
	\end{enumerate}
\end{prop}

\begin{proof}
  Informally, the equivalences are simply described by transposition
  using the provided adjunctions,
  for example
  \begin{equation}
    (b,c,Sb\to Gc) \leftrightarrow (b,c, b\to \ra{S}Gc)
  \end{equation}
  for the first statement.
  We refer to \cite[Cor.~A.5]{CDW24} for a general \((\infty,2)\)-categorical proof.
\end{proof}

\subsection{Commutative squares}

Let 
\begin{equation}
  \label{eq:csquare}
  \cdsquare{\A}{\B}{\C}{\D}{F}{R}{S}{G}
\end{equation}
be a commutative square in $\Stk$.
We introduce two means of totalizing \eqref{eq:csquare}
to obtain a categorical complex:
\begin{enumerate}
\item We define the {\em product totalization} to be the sequence of functors
  \begin{equation}
    \label{eq:totpsquare}
    \begin{tikzcd}
      \A \ar[r, "d_2"] & \laxpull{\B}{\D}{\C} \ar[r, "d_1" ] & \D
    \end{tikzcd}
  \end{equation}
  where $d_2$ is the canonical functor arising from the square \eqref{eq:csquare},
  interpreted as a directed cone over \eqref{eq:cospan}, i.e.
  \[
    a \mapsto (F(a),R(a),S(F(a)) \simeq G(R(a)))
  \]
  while the functor $d_1$ is given by
  \[
    (b, c, \eta\colon S(b) \to G(c))  \mapsto \fib(\eta).
  \]
  It is immediate that we have $d_1 \circ d_2 \simeq 0$ so that we may interpret
  \eqref{eq:totpsquare} as a categorical complex with $\D$ in degree $0$.
\item We define the {\em coproduct totalization} to be the sequence of functors
  \begin{equation}
    \label{eq:totcsquare}
    \begin{tikzcd}
      \A \ar[r, "d_2"] & \laxpush{\B}{\A}{\C} \ar[r, "d_1" ] & \D
    \end{tikzcd}
  \end{equation}
  where the functor $d_2$ is the cofiber of the natural transformation in the universal
  cone \eqref{eq:universalpush} while $d_1$ is the canonical functor arising from the
  square \eqref{eq:csquare}, interpreted as a directed cone under \eqref{eq:span}.
  Again, it is evident, that $d_1 \circ d_2 \simeq 0$: The postcomposition of the said
  natural transformation from \eqref{eq:universalpush} with $d_1$ is a natural
  equivalence, since the square \eqref{eq:csquare} commutes. Thus, we may interpret
  \eqref{eq:totcsquare} as a categorical complex as well, concentrated in degrees $2$ to $0$.
\end{enumerate}

While the two categorical totalizations of the square \eqref{eq:csquare} yield isomorphic complexes
of abelian groups when passing to $K_0$ (as can be deduced from the existence of semiorthogonal
decompositions $(\B,\C)$ on both directed pushout and directed pullback), they are not in general
equivalent as categorical complexes. Nevertheless, there exists a canonical comparison morphism
which will be investigated next.

We construct a functor
	\begin{equation}
		\label{eq:chi}
		\chi\colon \laxpush{\B}{\A}{\C} \lra \laxpull{\B}{\D}{\C}
	\end{equation}
in terms of the universal properties of both sides:
\begin{enumerate}[label=\arabic*.]
	\item the functor $\B \to \laxpull{\B}{\D}{\C}$ is given by 
		\[
			b \mapsto (b,0,S(b) \to 0)
		\]
	\item the functor $\C \to \laxpull{\B}{\D}{\C}$ is given by 
		\[
			c \mapsto (0,c[1],0 \to G(c)[1])
		\]
	\item the functor $\A \to \Fun(\Delta^1,\laxpull{\B}{\D}{\C})$ is given by 
		\[
		a \mapsto \begin{tikzcd}
			( & F(a) \ar[d], & 0 \ar[d], &  S \circ F(a) \ar[d] \ar[r] & 0 \ar[d] & )\\
			( & 0, & R(a)[1], &  0 \ar[r] & G \circ R(a)[1] & )
		\end{tikzcd}
		\]
		where the square is the biCartesian square exhibiting the equivalence $S \circ F
		\simeq G \circ R$. 
\end{enumerate}

\begin{prop}\label{prop:coprod_prod}
  The functor $\chi$ extends to a morphism of categorical complexes
	\[
	\begin{tikzcd}
		\A \ar[r, "d_2"]\ar[d,"{[1]}"] &
    \laxpush{\B}{\A}{\C} \ar[r, "d_1" ] \ar[d, "\chi"]
    & \D \ar[d, "{\id}"]\\
		\A \ar[r, "{d_2}"] & \laxpull{\B}{\D}{\C} \ar[r, "{d_1}" ] & \D
	\end{tikzcd}
	\]
	between the coproduct and product totalizations of \eqref{eq:csquare}.
\end{prop}
\begin{proof} Direct computation.
\end{proof}

We now provide a criterion for when the functor $\chi$ is an equivalence so that in particular, by Proposition
\ref{prop:coprod_prod}, the coproduct and product totalization of \eqref{eq:csquare} will be
equivalent as categorical complexes. 
The square \eqref{eq:csquare} is called {\em vertically right adjointable}, if the functors $R$ and $S$ have right
adjoints and the resulting mate
\begin{equation}
	\label{eq:bcsquare}
  \laxsquare{\C}{\D}{\A}{\B}{G}{\ra{R}}{\ra{S}}{F}{}
\end{equation}
is a natural equivalence. 

\begin{prop}\label{prop:bcsquare} Suppose that the square \eqref{eq:csquare}
  is vertically right adjointable. Then the functor $\chi$ is an equivalence.
  In particular, 
  the product and coproduct totalizations are canonically equivalent.
\end{prop}
\begin{proof}
	One may explicitly verify that we have a commutative square
	\[
		\begin{tikzcd}
			\laxpush{\B}{\A}{\C} \ar[rr, "\chi\desuspension"]\ar[d, "\simeq"] &  & \laxpull{\B}{\D}{\C} \ar[d, "\simeq"]\\
			\laxcolima{\C}{F \circ \ra{R}}{\B} \ar[r, "\simeq"] & \laxcolima{\C}{\ra{S} \circ G}{\B} \ar[r, "\simeq"] & \oplaxlima{\C}{\ra{S} \circ G}{\B}
		\end{tikzcd}
	\]
	where all equivalences are given by the canonical identifications that arise from the
	various universal cones described in \S\ref{subsec:laxsum}.
\end{proof}

The goal of the subsequent parts of \S \ref{sec:tot} will be to generalize the totalization of
squares to bicomplexes, and finally multicomplexes -- with a special focus on cubes, as these are our
main examples.

\subsection{Categorical bicomplexes}\label{subsec:totofbicomplexes}

A {\em categorical bicomplex} $\A_{\bullet,\bullet}$ consists of the datum of
\begin{itemize}
	\item a family $\{\A_{(i,j)}\}_{(i,j) \in \ZZ^{\op,2}}$ of objects in $\Stk$,
	\item functors $d\colon \A_{(i,j)} \to \A_{(i-1,j)}$ and $\delta\colon \A_{(i,j)} \to \A_{(i,j-1)}$,
	\item equivalences $d\delta \simeq \delta d$.
\end{itemize}
such that
\begin{itemize}
	\item $d^2 \simeq 0$ and $\delta^2 \simeq 0$.
\end{itemize}
Up to contractible choices, we may identify a categorical bicomplex with a functor $\ZZ^{\op,2} \to \St$
(with the properties $d^2 = 0$ and $\delta^2 = 0$). 
Given a categorical bicomplex $\A_{\bullet,\bullet}$, we define the {\em product totalization} to be the categorical complex 
\[
	\C_{\bullet} = \tot^{\times}(\A_{\bullet,\bullet})
\]
given as follows:
\begin{enumerate}
	\item \label{tot:1} for $n \in \ZZ^\op$, the category $\C_n$ is the iterated \dpb{} 
		\[
			\cdots \underset{\A_{1,n-2}}{\overset{\curvearrowright}{\times}}\A_{1,n-1}
				\underset{\A_{0,n-1}}{\overset{\curvearrowright}{\times}} \A_{0,n}
			\underset{\A_{-1,n}}{\overset{\curvearrowright}{\times}}\A_{-1,n+1} 
			\underset{\A_{-2,n+1}}{\overset{\curvearrowright}{\times}}
			\cdots
		\]
		Thus, an object of this category consists of a sequence $\{a_{i,j}\}_{i+j = n}$
		together with, for every $i,j$, a specified morphisms $d a_{i,j} \to \delta a_{i-1,j+1}$ in
		$\A_{i-1,j}$.
	\item \label{tot:2} the differential $d_{\C}\colon \C_n \to \C_{n-1}$ is defined as follows: for an object
		$\{a_{i,j}\}$ of $\C_n$, the component in $\A_{k,l}$ of its image under $d_{\C}$ is given by 
		\[
			\fib( d a_{k+1,l} \to \delta a_{k,l+1})[k]
		\]
		The images of adjacent components of $d_{\C}(\{a_{i,j}\})$ under $d$ and $\delta$ are canonically equivalent
		via the chosen equivalence $d \delta \simeq \delta d$. The components, equipped with
		these equivalences then define an object of $\C_{n-1}$.
\end{enumerate}

Note that, in contrast to its classical counterpart, the definition of the categorical product
totalization is not symmetric in each term (ignoring the differential) with respect to swapping
the coordinates of the bicomplex, due to the appearance of the iterated \dpb{} in
\ref{tot:1}. Our convention is to use the linear order of the coordinates of the bicomplex to
determine the chosen direction.

The coproduct totalization $\C_\bullet'=\tot^{\amalg}(\A_{\bullet,\bullet})$ is defined in analogy to
\eqref{eq:totcsquare}: 
\begin{enumerate}
	\item \label{tot:1} for $n \in \ZZ^\op$, the category $\C_n'$ is the iterated \dpo{} 
		\[
			\cdots \underset{\A_{2,n-1}}{\overset{\curvearrowright}{\amalg}}\A_{1,n-1}
				\underset{\A_{1,n}}{\overset{\curvearrowright}{\amalg}} \A_{0,n}
			\underset{\A_{0,n+1}}{\overset{\curvearrowright}{\amalg}}\A_{-1,n+1} 
			\underset{\A_{-1,n+2}}{\overset{\curvearrowright}{\amalg}}
			\cdots
		\]
	Thus, via the universal property, a functor $f:\C_n'\rightarrow \D$ corresponds to a collection of functors $\alpha_{i,j}:\A_{i,j}\rightarrow \D$ with $i+j=n$ and natural transformations $\eta_{i,j}$ fitting into diagrams of the following form:
	\[
	\begin{tikzcd}
{\A_{i+1,j}} \arrow[r, "d"] \arrow[d, "\delta"'] & {\A_{i,j}} \arrow[d, "\alpha_{i,j}"] \arrow[ld, Rightarrow, "\eta_{i,j}"'] \\
{\A_{i+1,j-1}} \arrow[r, "\alpha_{i+1,j-1}"']         & \D                                         
\end{tikzcd}
	\]
	\item \label{tot:2} We specify the composite of the differential $d_{\C'}\colon \C_{n+1}' \to \C_{n}'$ with any functor $f:\C_n'\rightarrow \D$ as above. Specializing to $f=\on{id}_{\C_n'}$ yields the differential $d_{\C'}$. We describe $f\circ d_{\C'}$ via the universal property of the iterated directed pushout in terms of functors 
	\[ \beta_{i+1,j}=\on{fib}(\eta_{i,j})[i+1]\colon \A_{i+1,j}\longrightarrow \D\]
	and natural transformations 
	\begin{equation}\label{eq:nateqdifferential} \beta_{i+1,j}\circ \delta \simeq \beta_{i,j+1}\circ d\end{equation}
	arising from the identities $d^2\simeq 0$, $\delta^2\simeq 0$ and $\delta\circ d\simeq d\circ \delta$. It follows from the fact that the natural transformation \eqref{eq:nateqdifferential} is a natural equivalence, that $d_{\C'}^2\simeq 0$.
\end{enumerate}

\begin{rem}
With some more effort, one can show that product and coproduct totalizations form functors
\[ 
\on{tot}^{\times},\on{tot}^{\amalg}\colon \on{Ch}_2(\Stk)\rightarrow \on{Ch}(\Stk)\,,
\]
but we omit the details in this work.
\end{rem}

The most important example of a totalization of a bicomplex for us will be the following special case.

\begin{con}
  \label{con:Fib-and-Cof}
  Let $F\colon \A_{\bullet} \to \B_{\bullet}$ be a morphism of categorical complexes, depicted as follows.
  \begin{equation}\label{eq:morphism_bicomplex}
    \begin{tikzcd}
      \dots \arrow[r] & \A_2 \arrow[d, "F_2"] \arrow[r, "d^{\A}"] & \A_{1} \arrow[r, "d^{\A}"] \arrow[d, "F_1"] & \A_0 \arrow[r] \arrow[d, "F_0"] & \dots \\
      \dots \arrow[r] & \B_2 \arrow[r, "d^{\B}"]                  & \B_1 \arrow[r, "d^{\B}"]                    & \B_0 \arrow[r]                  & \dots
    \end{tikzcd}
  \end{equation}
  We may interpret $F$ as a
  bicomplex $\C_{\bullet,\bullet}$ with $\C_{0,\bullet} = \A_{\bullet}$ and $\C_{-1,\bullet} = \B_{\bullet}$. We define the {\em fiber} of $F$ as the product totalization $\Fib(F) \coloneqq  \tot^{\times}(\C_{\bullet,\bullet})$. Explicitly, we have 
  \[
		\Fib(F)_i = \laxpull{\A_{i}}{\B_{i}}{\B_{i+1}}, 
  \]
  The differential $d\colon \Fib(F)_i\rightarrow \Fib(F)_{i-1}$ is given by
  \begin{equation}
    \label{eq:diff_tot} 
    (a, b, \eta\colon F_i(a) \to d(b)) \mapsto (d(a), \fib(\eta), F_{i-1}d(a) \simeq dF_i(a)),
  \end{equation}
  where we note that $d\fib(\eta) \simeq dF_i(a)$. 

  Dually, we define the {\em cofiber} $\Cof(F)$ of $F$ as the coproduct totalization of the
  bicomplex \eqref{eq:morphism_bicomplex}, with a degree shift of $-1$, so that we have
  \[ 
    \Cof(F)_i \simeq \laxpush{\A_{i-1}}{\A_{i}}{\B_{i}}.
  \] 
 \end{con}

 \begin{defi}
   \label{defi:shift-of-complex}
   For a chain complex \((\A_\bullet,d^\A)\) we define its shift
   \(\A[1] \coloneqq \Cof(\A\to 0)\),
   which is explicitly given by
   \begin{equation}
     \A[1]_i \coloneqq \A_{i-1}
     \quad
     \text{and}
     \quad
     d^{\A[1]}_{i+1}\coloneqq d^\A_{i}[1] \colon \A_i\to \A_{i-1}.
   \end{equation}
   Dually we define \(\A[-1]\coloneqq \Fib(0\to \A)\)
   and observe
   \begin{equation}
     \A[1][-1]\simeq \A \simeq \A[-1][1].
   \end{equation}
 \end{defi}

\begin{rem}\label{rem:homotopy_cofiber} The cofiber (resp. fiber) satisfy universal properties which can be formulated in terms
	of the notion of categorical homotopy introduced in \S\ref{sub:mapping_complexes}. For
	example, the cofiber is the universal example of a categorical complex $C_{\bullet}$
	equipped with a morphism $G\colon \B_{\bullet} \to \C_{\bullet}$ together with a categorical zero
	homotopy of the composite $G \circ F$. See also \cite[Section~9]{CDW24} for more details. We expect that this phenomenon will become
	relevant when trying to introduce a notion of {\em derived category} of categorical
	complexes. 
\end{rem}

In light of Remark \ref{rem:homotopy_cofiber} one may anticipate
that fiber and cofiber of a given morphism are equivalent up to a shift.
While this is not the case for a general morphism,
remarkably, there is a natural class of morphisms for which the statement holds and which we
introduce next. 

\begin{defi}
  \label{defi:adjointable-chain-map}
  A chain map $F\colon \A_{\bullet} \to \B_{\bullet}$
  between categorical complexes is called {\em right adjointable} if each square 
	\begin{equation}
		\label{eq:com_square_fun}
    \cdsquare{\A_i}{\A_{i-1}}{\B_i}{\B_{i-1}}
    {d^\A}{F_i}{F_{i-1}}{d^\B}
	\end{equation}
  is vertically right adjointable. 
\end{defi}
\begin{prop}[Beck--Chevalley]
	\label{prop:beckchevalley} 
	Let $F\colon \A_{\bullet} \to \B_{\bullet}$ be a right adjointable morphism. Then there is a
	canonical equivalence of categorical complexes
	\[
		\chi_{\bullet}\colon \Cof(F) \overset{\simeq}{\lra} \Fib(F)[1]
	\]
	where 
	\[
		\chi_i\colon \laxpush{\A_{i-1}}{\A_{i}}{\B_{i}} \lra \laxpull{\A_{i-1}}{\B_{i-1}}{\B_{i}}, 
	\]
	is the functor defined in \eqref{eq:chi} corresponding to the commutative square
	\eqref{eq:com_square_fun}.
\end{prop}
\begin{proof} Each functor $\chi_i$ is an equivalence by Proposition \ref{prop:bcsquare}. The fact
	that $\chi_{\bullet}$ defines a morphism of complexes is verified by direct computation
	using the universal properties of the terms involved. An alternative argument is also provided in \cite[Cor.~8.30]{CDW24}.
\end{proof}

\subsection{Categorical multicomplexes}

Consider a categorical $n$-complex as in \Cref{def:ncplx}. To iteratively define its totalization, we begin by introducing the partial totalization at adjacent coordinates. 

\begin{con}
Let $n\geq 2$ and $\A_\ast:\mathbb{Z}^{\op,n}\rightarrow \St_k$ be a categorical $n$-complex. We consider $\A_\ast$ as an object of $\on{Fun}(\mathbb{Z}^{\op,n-2},\on{Ch}_2(\Stk))$, valued in the bicomplexes describing chosen coordinates $1\leq i,i+1\leq n$. Using the functoriality of the totalization constructions of \Cref{subsec:totofbicomplexes}, we may partially product totalize the chosen coordinates to obtain a categorical $(n-1)$-complex $\tot^{\times}_{i,i+1}(\A_\ast)$, defined as the image under the functor 
\[ \tot^{\times}_{i,i+1}:\on{Fun}(\mathbb{Z}^{\op,n-2},\on{Ch}_2(\Stk))\xlongrightarrow{\on{Fun}(\mathbb{Z}^{\op,n-2},\tot^\times)} \on{Fun}(\mathbb{Z}^{\op,n-2},\on{Ch}(\St_k))\,.\]
The partial coproduct totalization $\tot^{\amalg}_{i,i+1}(\A_\ast)\in \on{Ch}_{n-1}(\Stk)$ is defined similarly.
\end{con}

Given a categorical $n$-complex $\A_\ast$, we obtain its product and coproduct totalizations by iterating the above partial totalizations. For the moment, we fix one particular order for these partial totalizations.

\begin{defi}
Let $\A_\ast$ be a categorical $n$-complex. We define its product totalization $\on{tot}^{\times}(\A_\ast)$ as the iterated partial totalization 
\[ \on{tot}^{\times}_{1,2}\circ \on{tot}^{\times}_{2,3}\circ \dots\circ \on{tot}^{\times}_{n-1,n}(\A_\ast)\,.
\]
Similarly, the coproduct totalization  $\on{tot}^{\amalg}(\A_\ast)$ of $\A_\ast$ is defined as 
\[ \on{tot}^{\amalg}_{1,2}\circ \on{tot}^{\amalg}_{2,3}\circ \dots\circ \on{tot}^{\amalg}_{n-1,n}(\A_\ast)\,.
\]
\end{defi}

With our conventions, both the product and coproduct totalizations of a categorical $n$-cube lie in degrees $n$ to $0$. 

\begin{exa}\label{exa:higherauslander}
	Suppose that $q: \I^n \to \Stk$ is a constant cube with value $\A$. By an inductive
	argument, we obtain the following description of the product totalization of $q$: 
	\begin{itemize}
		\item $\C_{n+1} = \A$,
		\item the category $\C_{n-k}$ consists of diagrams
			\[
				X\colon \Fun([k],[n]) \to \A
			\]
			mapping non-injective maps $[k] \to [n]$ to zero objects,
		\item the differential $d\colon \C_{n-k} \to \C_{n-k-1}$ is given by associating to a
			diagram 
			\[
				X\colon \Fun([k],[n]) \to \A
			\]
			the diagram 
			\[
				X'\colon \Fun([k+1],[n]) \to \A
			\]
			where $X'(\sigma)$ is the total fiber of the complex
			\[
				X(\sigma \circ \partial_n) \to X(\sigma \circ \partial_{n-1}) \to
				... \to X(\sigma \circ \partial_0)
			\]
			in $\A$, where $X'$ is formally constructed from $X$ by means of Kan extensions.
	\end{itemize}
	In particular, we observe that the $n$-cells of this complex are representations in $\A$ of
	higher Dynkin type $\AA$. Note that this complex is an augmented variant of categorical
	cochain complex $\C^{\bullet}(\DDelta, \A)$ from Definition \ref{defi:catcochain}
	(normalized slightly differently there, since the differentials are given by taking the
	total cofiber instead of the fiber). The example illustrates the general phenomenon that,
	even for simple cubical diagrams, the terms of the totalization turn out to be rather
	interesting stable $\infty$-categories. This is somewhat in contrast to the usual
	totalization of cubical diagrams of abelian groups, where the terms of the totalization are
	simply direct sums of the terms of the cube. 
	The coproduct totalization of $q$ admits a dual description, we leave the analogous detailed
	description to the reader. It follows from the Beck--Chevalley property of the cube $q$,
	that the product and coproduct totalizations are in fact equivalent as categorical complexes
	(repeatedly apply \Cref{prop:beckchevalley}). The totalization of $q$ is also an example of a categorical
	Koszul complex, discussed in \S \ref{sec:koszul}. See \Cref{thm:catkoszul} and
	\Cref{cor:duality} for categorical Koszul duality and its implications when applied to $q$.
\end{exa}

\begin{rem}
We expect a fully coherent associativity statement for the totalization of multicomplexes, meaning
that the totalization does not depend on the order of the totalizations of the adjacent coordinates. A systematic analysis of the associativity of the totalizaton would however go beyond the scope of this work. Instead, we sketch how to prove a partial result in \Cref{prop:totasso}
below, needed for concrete applications (e.g. in \S \ref{sec:koszul}).

The totalization of a cube depends on a further choice, namely the given total order of the
coordinates. We leave it as an interesting problem to determine under which conditions the
totalization does not depend on this order, up to equivalence.
\end{rem}

\begin{prop}\label{prop:totasso}
Let $\A_\ast$ be a categorical $n$-cube. Then there exist equivalences of categorical complexes
\[ \on{tot}^{\times}(\A_\ast)\simeq \on{tot}^{\times}_{1,2}\circ \dots \circ \on{tot}^{\times}_{1,2}(\A_\ast)\]
and
\[ \on{tot}^{\amalg}(\A_\ast)\simeq \on{tot}^{\amalg}_{1,2}\circ \dots \circ \on{tot}^{\amalg}_{1,2}(\A_\ast)\,.\]
\end{prop}

Informally, \Cref{prop:totasso} states that iteratively totalizing the last two coordinates is equivalent to iteratively totalizing the first two coordinates.

To prove \Cref{prop:totasso}, we begin with the case of categorical $3$-complexes of the form $\I\times \mathbb{Z}^\op\times \I\rightarrow \Stk$, where $\I=[1]^{\on{op}}\subset \ZZ^\op$.

\begin{lem}\label{lem:totasso}
Let $\on{Fun}_{\on{Ch}}(\I\times \mathbb{Z}^\op\times \I,\Stk)\subset \on{Fun}(\I\times \mathbb{Z}^\op\times \I,\Stk)$ be the full subcategory spanned by categorical $3$-complexes. There exists a natural equivalence between the functors
\[ 
\on{tot}^{\times}_{1,2}\circ \on{tot}^{\times}_{1,2}:\on{Fun}_{\on{Ch}}(\I\times \mathbb{Z}^\op\times \I,\Stk)\rightarrow \on{Ch}(\Stk)
\]
and
\[ 
\on{tot}^{\times}_{1,2}\circ \on{tot}^{\times}_{2,3}:\on{Fun}_{\on{Ch}}(\I\times \mathbb{Z}^\op\times \I,\Stk)\rightarrow \on{Ch}(\Stk)\,.
\]
\end{lem}

\begin{proof}
We can depict a part of $\A_\ast$ as follows, 
\[
\begin{tikzcd}[column sep=tiny]
\A_{1,i+1,1} \arrow[r] \arrow[d]  & \A_{1,i+1,0} \arrow[rrd] &                                      &                            &                          &                \\
\A_{0,i+1,1} \arrow[r] \arrow[rrd]         & \A_{0,i+1,0} \ar[from=u]\arrow[rrd]           & \A_{1,i,1}\ar[from=llu,crossing over]  \arrow[r] & \A_{1,i,0} \arrow[d] \arrow[rrd] &                          &                \\
                                   &                            & \A_{0,i,1} \ar[from=u, crossing over] \arrow[r] \arrow[rrd]           & \A_{0,i,0} \arrow[rrd]           & \A_{1,i-1,1} \ar[from=llu,crossing over] \arrow[r] & \A_{1,i-1,0} \arrow[d] \\
                                   &                            &                                      &                            & \A_{0,i-1,1} \ar[from=u, crossing over] \arrow[r]           & \A_{0,i-1,0}         
\end{tikzcd}
\]
and relabel this part for better readability as follows:
\[
\begin{tikzcd}
\A_2 \arrow[r] \arrow[d]  & \B_1 \arrow[d] \arrow[rrd] &                                      &                            &                          &                \\
\B_4 \arrow[r] \arrow[rrd]         & \C_3 \arrow[rrd]           & \B_2  \ar[from=llu,crossing over] \arrow[r] & \C_1 \arrow[d] \arrow[rrd] &                          &                \\
                                   &                            & \C_4 \ar[from=u, crossing over] \arrow[r] \arrow[rrd]           & \D_3 \arrow[rrd]           & \C_2 \ar[from=llu,crossing over]\arrow[r]& \D_1 \arrow[d] \\
                                   &                            &                                      &                            & \D_4 \ar[from=u, crossing over] \arrow[r]           & \E_3          
\end{tikzcd}
\]
We have by definition 
\[ 
\tot^{\times}_{1,2}\circ \tot^\times_{1,2}(\A_\ast)_{i+1}=\laxpull{(\laxpull{\C_1}{\D_3}{\C_3})}{(\laxpull{\D_1}{\E_3}{\D_3})}{(\laxpull{\C_2}{\D_4}{\C_4})}
\]
and 
\[ \tot^{\times}_{1,2}\circ \tot^{\times}_{2,3}(\A_\ast)_{i+1}=\laxpull{(\laxpull{\C_1}{\D_1}{\C_2})}{(\laxpull{\D_3}{\E_3}{\D_4})}{(\laxpull{\C_3}{\D_3}{\C_4})}\,.\]
Unraveling the definition, we find that both of these $\infty$-categories satisfy the following universal property: a functor from $\X \in \Stk$ corresponds to
\begin{itemize}
\item functors $F_i\colon \X\to \C_j$ in $\Stk$ for all $1\leq j\leq 4$,
\item natural transformations 
\[ \begin{tikzcd}[row sep=small]
                         & { \C_1} \arrow[rd] \arrow[dd, Rightarrow, "\alpha"] &              \\
\X \arrow[ru, "F_1"] \arrow[rd, "F_2"] &                                                 & {\D_1} \\
                         & {\C_2} \arrow[ru]                       &             
\end{tikzcd}
\quad \begin{tikzcd}[row sep=small]
                         & { \C_1} \arrow[rd, ""] \arrow[dd, Rightarrow, "\gamma"] &                \\
\X \arrow[ru, "{F_1[-1]}"] \arrow[rd, "F_3"] &                                                   & {\D_3} \\
                         & {\C_3} \arrow[ru]                           &               
\end{tikzcd}\]
\[ \begin{tikzcd}[row sep=small]
                         & { \C_2} \arrow[rd] \arrow[dd, Rightarrow, "\beta"] &                \\
\X \arrow[ru, "F_2"] \arrow[rd, "F_4"] &                                                 & {\D_4} \\
                         & { \C_4} \arrow[ru]                      &               
\end{tikzcd}
\quad \begin{tikzcd}[row sep=small]
                         & { \C_3} \arrow[rd] \arrow[dd, Rightarrow, "\delta"] &              \\
\X \arrow[ru, "F_3"] \arrow[rd, "F_4"] &                                                   & {\D_3} \\
                         & { \C_4} \arrow[ru]                          &             
\end{tikzcd}
\] 
and a null-homotopy of $\delta\circ \gamma$. This null-homotopy induces a natural transformation $\nu$ from the functor $\X\xrightarrow{F_1}\C_1\to \D_3$ to $\cof(\delta)$. 
\item A natural equivalence between the two composite natural transformations 
\[
\begin{tikzcd}
                                                                   & \C_1 \arrow[rd] \arrow[d, "\nu", Rightarrow] &                                                         &      \\
\X \arrow[ru, "F_1"] \arrow[rr, "\cof(\delta)"'] \arrow[rd, "F_4"'] & {}                                           & \D_3 \arrow[r] \arrow[d, "\simeq", Rightarrow, no head] & \E_3 \\
                                                                   & \C_4 \arrow[r]                               & \D_3 \arrow[ru]                                         &     
\end{tikzcd}\quad\quad
\begin{tikzcd}
                                                        & \C_1 \arrow[rd] \arrow[d, "\alpha", Rightarrow]          &                                                         &      \\
\X \arrow[ru, "F_1"] \arrow[rd, "F_4"'] \arrow[r, "F_2"] & \C_2 \arrow[r] \arrow[rd] \arrow[d, "\beta", Rightarrow] & \D_1 \arrow[r] \arrow[d, "\simeq", Rightarrow, no head] & \E_3 \\
                                                        & \C_4 \arrow[r]                                           & \D_4 \arrow[ru]                                         &     
\end{tikzcd}
\]
between the functors $\X\xrightarrow{F_1}\C_1\to \E_3$ and $\X\xrightarrow{F_4}\C_4\to \E_3$. The lower natural equivalence in the left diagram arises from the fact that the functor $\C_3\to \D_3\to \E_3$ is zero, since $d^2$ is zero.
\end{itemize}
One checks that the differentials of both complexes are identified in terms of degreewise equivalences arising from the matching universal properties. All arising suspensions can be dealt with by including suitable equivalences. Using \eqref{eq:laxpullvialimit}, we can express all directed pullbacks appearing above in terms of $\infty$-categorical limits and the above equivalences in terms of functorial equivalences between these limits, yielding the desired equivalence between $\on{tot}^{\times}_{1,2}\circ \on{tot}^{\times}_{1,2}$ and $\on{tot}^{\times}_{1,2}\circ \on{tot}^{\times}_{2,3}$.
\end{proof}

\begin{proof}[Proof of \Cref{prop:totasso}.]
We only prove the statement about the product totalization, the case of the coproduct totalization is analogous. Repeatedly applying \Cref{lem:totasso}, we find equivalences for all $1<j<n$
\begin{align*} 
\on{tot}^{\times}_{j-1,j}\circ \on{tot}^{\times}_{j,j+1} \dots\circ \on{tot}^{\times}_{j,j+1}(\A_\ast) &\simeq \on{tot}^{\times}_{j-1,j}\circ \on{tot}^{\times}_{j-1,j}\circ \on{tot}^{\times}_{j+1,j+2}\circ \dots\circ \on{tot}^{\times}_{j+1,j+2}(\A_\ast) \\
&\simeq \dots\\
&\simeq  \on{tot}^{\times}_{j-1,j}\circ \dots\circ \on{tot}^{\times}_{j-1,j}\circ \on{tot}^{\times}_{j-1,j}(\A_\ast)\,.
\end{align*}
Combining these equivalences, we have
\begin{align*}
\on{tot}^{\times}_{1,2}\circ\dots\circ \on{tot}^{\times}_{n-1,n}(\A_\ast)& \simeq \on{tot}^{\times}_{1,2}\circ\dots\circ \on{tot}^{\times}_{n-3,n-2}\circ \on{tot}^{\times}_{n-2,n-1}\circ \on{tot}^{\times}_{n-2,n-1}(\A_\ast) \\
& \simeq \on{tot}^{\times}_{1,2}\circ\dots\circ \on{tot}^{\times}_{n-3,n-2}\circ \on{tot}^{\times}_{n-3,n-2}\circ \on{tot}^{\times}_{n-3,n-2}(\A_\ast) \\
& \simeq \dots \\
& \simeq \on{tot}^{\times}_{1,2}\circ\dots \circ \on{tot}^{\times}_{1,2}(\A_\ast)\,,
\end{align*}
as desired.
\end{proof}

\subsection{Beck-Chevalley conditions and spherical functors}
\label{subsec:beck-chevalley}

\begin{defi}\label{def:adjointable}
Consider a commutate square in $\St_k$:
	\begin{equation}\label{eq:BCsqu}
	\begin{tikzcd}
\A \arrow[r, "G"] \arrow[d, "F"] & \B \arrow[d, "F'"] \\
\C \arrow[r, "G'"]               & \D                
\end{tikzcd}
	\end{equation}
	We say that the square is
	\begin{enumerate}[(1)]
	\item  \textit{horizontally right adjointable}, or simply right adjointable, if $G$ and $G'$ admit $k$-linear right adjoints $\ra{G}$ and $\ra{(G')}$ and the natural transformation
	\begin{equation*}
		F\circ \ra{G} \xRightarrow{\on{u}\circ F\circ \ra{G}} \ra{(G')}\circ G' \circ F\circ \ra{G} \xRightarrow{\simeq} \ra{(G')}\circ F'\circ G\circ \ra{G}  \xRightarrow{\ra{(G')}\circ F' \circ \on{cu}}\ra{(G')}\circ F' 
		\end{equation*}	
		is a natural equivalence.
		\item  \textit{horizontally left adjointable}, if $G$ and $G'$ admit $k$-linear left adjoints $\la{G}$ and $\la{(G')}$ and the natural transformation
		\begin{equation}\label{eq:BC2}
		\la{(G')}\circ F' \xRightarrow{\la{(G')}\circ F' \circ \unit} \la{(G')}\circ F'\circ G \circ \la{G} \xRightarrow{\simeq} \la{(G')}\circ G'\circ F \circ \la{G}  \xRightarrow{\counit \circ F\circ \la{G}} F\circ \la{G}
		\end{equation}
		is a natural equivalence.
		\item \textit{vertically right adjointable}, if $F$ and $F'$ admit $k$-linear right adjoints $\ra{F}$ and $\ra{(F')}$ and the natural transformation
		\begin{equation}\label{eq:BC3}
		G\circ \ra{F} \xRightarrow{\on{u}\circ G\circ \ra{F}} \ra{(F')}\circ F' \circ G\circ \ra{F} \xRightarrow{\simeq} \ra{(F')}\circ G'\circ F\circ \ra{F}  \xRightarrow{\ra{(F')}\circ G' \circ \on{cu}}\ra{(F')}\circ G' 
	\end{equation}
		is a natural equivalence.
		\item \textit{vertically left adjointable}, if $F$ and $F'$ admit $k$-linear left adjoints $\la{F}$ and $\la{(F')}$ and the natural transformation
		\begin{equation*}
		\la{(F')}\circ G' \xRightarrow{\la{(F')}\circ G' \circ \unit} \la{(F')} \circ G'\circ F \circ \la{F} \xRightarrow{\simeq} \la{(F')} \circ F'\circ G \circ \la{F}  \xRightarrow{\counit \circ G\circ \la{F}} G\circ \la{F}
		\end{equation*}
		is a natural equivalence.
	\end{enumerate}
\end{defi}

\begin{rem}\label{rem:BCvariant}
If $F$ and $F'$ admit right adjoints and $G$ and $G'$ admit left adjoints, then the natural transformations \eqref{eq:BC2} and \eqref{eq:BC3} are adjoint to another and conditions (2) and (3) from \Cref{def:adjointable} hence equivalent. An analogous statement holds for conditions (1) and (4).
\end{rem} 

Remarkably, all adjointability conditions are equivalent if the square is spherical.

\begin{prop}\label{prop:sphBCagree}
Consider a commutative square in $\Stk$ as in \eqref{eq:BCsqu} and suppose that all functors $F,F',G,G'$ are spherical functors, see \Cref{def:sphadj}. Then all four conditions of \Cref{def:adjointable} are equivalent.
\end{prop}

\begin{proof}
To begin with, we note that spherical functors admit all repeated left and right adjoints, see \cite[Cor.~2.5.16]{DKSS:spherical}. It thus follows that conditions (2) and (3) of \Cref{def:adjointable} are equivalent, as are conditions (1) and (4). It remains to show that conditions (1) and (2) are equivalent. Let $T$ be the cotwist functor of $G\dashv \ra{G}$ and $T'$ the twist functor of $G'\dashv \ra{(G')}$. Then we have $\la{G}\circ T \simeq \ra{G}$ and $T'\circ \la{(G')}\simeq \ra{(G')}$. We find the following commutative diagram,
\[
\begin{tikzcd}[row sep=large, column sep=5.5em]
0 \arrow[r] \arrow[d]                                                                   & T'\circ \la{(G')}\circ F'\circ T \arrow[r] \arrow[d, "{T'\circ \la{(G')} \circ F' \circ \unit'\circ T}"'] \arrow[rd, "\square", phantom]          & 0 \arrow[d]              \\
F\circ \ra{G} \arrow[r, "\on{u}\circ F\circ \ra{G}"] \arrow[d] \arrow[rd, "\square", phantom] & \ra{(G')}\circ G' \circ F\circ \ra{G} \arrow[r, "\ra{(G')}\circ F' \circ \on{cu}"] \arrow[d, "T'\circ \counit' \circ F\circ \la{G} \circ T"] & \ra{(G')}\circ F' \arrow[d] \\
0 \arrow[r]                                                                             & T'\circ F\circ \la{G}\circ T \arrow[r]                                                                                                     & 0                       
\end{tikzcd}
\]
omitting some identifications in the center. Here, $\unit$ is the unit of $G'\dashv \ra{(G')}$ and $\counit$ the counit of $G\dashv \ra{G}$, $\unit'$ the unit of $\la{G}\dashv G$ and $\counit'$ the counit of $\la{(G')}\dashv G'$. The fact that the upper right square and the lower left square are biCartesian, i.e.~both pullback and pushout, follows from the general properties of units and counit of spherical adjunctions, see for instance Remark 2.9 and Lemma 2.10 in \cite{Chr20}. By the pasting laws for biCartesian squares, we find that the horizontal middle composite is a natural equivalence if and only if the vertical middle composite is a natural equivalence. The equivalence of conditions (1) and (2) now follows from the fact that $T$ and $T'$ are invertible.
\end{proof}

\begin{defi}
A categorical $n$-cube is called \textit{Beck-Chevalley} if each rectilinear face is both horizontally and vertically right adjointable. 
\end{defi}

\begin{lem}\label{lem:parBC}
Suppose that $\A_\ast$ is a Beck-Chevalley categorical $n$-cube. Then each partial product totalization 
\[ \on{tot}^{\times}_{i,i+1}\circ \dots\circ \on{tot}^{\times}_{n-1,n}(\A_\ast)\]
and partial coproduct totalization 
\[
\on{tot}^{\amalg}_{i,i+1}\circ \dots\circ \on{tot}^{\amalg}_{n-1,n}(\A_\ast)
\]
with $1\leq i\leq n-1$ is a Beck-Chevalley $(n-i)$-cube.
\end{lem}

\begin{proof}
We only show that $\tot^{\times}_{n-1,n}(\A_\ast)_\ast:\I^{\times n-2}\times [2]^{\op}\rightarrow \St_k$ is Beck-Chevalley. A similar argument applies to $\tot^{\amalg}_{n-1,n}(\A_\ast)$. The general case is proven analogously, by using that the $a$-times repeated partial totalization consists of $a$ stacked cubes, to each of which an analogous argument as below applies.

Let $J\in \I^{n-2}$. For $l=0,2$, we have
\[ \tot^{\times}_{n-1,n}(\A_\ast)_{(J,l)}=\A_{(J,l,l)}\,.\]
If instead $l=1$, let $J^1=(J,1,0),\,J^2=(J,0,1)\,,J^0=(J,0,0)\in \I^n\,.$ Then we have  
\[ \tot^{\times}_{n-1,n}(\A_\ast)_{(J,1)}\simeq \laxpull{\A_{J^1}}{\A_{J^0}}{\A_{J^2}}\,.\]
Consider a rectilinear face $f:\I^2\rightarrow \I^{n-2}\times [2]^{\op}$ with $f(1,1)$ differing from $f(0,0)$ by subtracting $1$ in the entries $i$ and $j$ with $1\leq i,j \leq n-1$ and $i\neq j$. There are a few cases to distinguish, in which similar arguments apply. We highlight two such cases. If $i,j\neq n-1$ and $f(0,0)_{n-1}=0,2\in [2]^{\op}$, then the face of $\tot^{\times}_{n-1,n}(\A_\ast)$ is equivalent to a corresponding face of $\A_\ast$ and hence both vertically and horizontally right adjointable. If $i,j\neq n-1$ and $f(0,0)_{n-1}=1\in [2]^{\op}$, then the face of $\tot^{\times}_{n-1,n}(\A_\ast)$ is of the following form.
\[
\begin{tikzcd}
{\laxpull{\A_{f(0,0)^1}}{\A_{f(0,0)^0}}{\A_{f(0,0)^2}}} \arrow[d] \arrow[r] & {\laxpull{\A_{f(1,0)^1}}{\A_{f(1,0)^0}}{\A_{f(1,0)^2}}} \arrow[d] \\
{\laxpull{\A_{f(0,1)^1}}{\A_{f(0,1)^0}}{\A_{f(0,1)^2}}} \arrow[r]           & {\laxpull{\A_{f(1,1)^1}}{\A_{f(1,1)^0}}{\A_{f(1,1)^2}}}          
\end{tikzcd}
\]
This diagram arises from the functoriality of directed pullbacks applied to a diagram containing two Beck-Chevalley faces of $\A_\ast$. The adjointability properties now  follow from the adjointability properties of these two squares and the observation that the right adjoints of the morphisms in the above diagrams are computed componentwise, see \cite[Appendix~A.2]{CDW24}.
\end{proof}

\subsection{Mapping complexes}%
\label{sub:mapping_complexes}

Given complexes $A_{\bullet}$, $B_{\bullet}$ of abelian groups, there is an associated mapping
complex $\Map(A_{\bullet},B_{\bullet})_{\bullet}$ with
\[
	\Map(A_{\bullet},B_{\bullet})_{n} = \prod_{i \in \ZZ} \Hom(A_i, B_{n+i})
\]
and differential given by the formula
\[
	d(\{f_i\})_k = d_{B} \circ f_k - (-1)^{k-1} f_{k-1} \circ d_{A}.
\]

We now introduce a categorical variant of this construction.
Let $\A_{\bullet}$, $\B_{\bullet} \in \Ch(\Stk)$ be
categorical complexes. The {\em oplax categorical mapping complex} between $\A_{\bullet}$ and
$\B_{\bullet}$ is defined as the product totalization
\begin{equation}\label{eq:catmap}
	\Mapoplax[\bullet]{\A_\bullet}{\B_\bullet}
  \coloneqq  \tot^{\times}\Stk(\A_{-\bullet},\B_{\bullet})
\end{equation}
of the bicomplex obtained by applying $\Stk(-,-)$. Explicity:
\begin{itemize}
\item
  an object of $\Mapoplax[n]{\A}{\B}$ consists of a sequence
  $\{F_i\}_{i \in \ZZ}$ of $k$-linear functors $F_i\colon \A_i \to \B_{n+i}$
  equipped with natural transformations
		\[
			\phi_{i}\colon F_{i-1} \circ d_{\A} \Rightarrow d_{\B} \circ F_{i};
		\]
	\item
    the differential
    \(\Mapoplax[n]{\A_\bullet}{\B_\bullet}\to\Mapoplax[n-1]{\A_\bullet}{\B_\bullet}\)
    associates to this datum the sequence
    \begin{equation}
      \label{eq:formula-diff-oplax}
			\{G_i = \fib(\phi_i)[-i]\}
    \end{equation}
		of functors $G_i\colon \A_i \to \B_{i+n-1}$
    together with the natural equivalences
    \begin{align}
      G_{i-1}d_{\A}&=\fib(\F_{i-2}d_\A\to d_\B F_{i-1})[-i+1]d_\A
      \\
      &\simeq d_\B F_{i-1} d_\A[-i] \simeq d_B\fib(\F_{i-1}d_\A\to d_\B F_{i})[-i]
      \\
      &= d_\B G_{i}
    \end{align}
\end{itemize}
Note that the ``categorical signs'', i.e., the powers of the suspension $[1]$, do not agree with the usual Koszul signs
upon passing to $K_0$. Since we will not use categorical mapping complex systematically in this work, we will not
attempt to address this (purely conventional) discrepancy. 

More pictorially, we may interpret
$\Mapoplax[n]{\A_\bullet}{\B_\bullet}$ as the stable $\infty$-category of
diagrams in $\Stk$ of the form 
\begin{equation}
  \label{eq:oplax-morphisms}
  \begin{tikzcd}
    \dots\ar[r,"d"]
    &
    \A_2\ar[r,"d"]
    \ar[d,"F_2"']
    &
    \A_1\ar[r,"d"]
    \ar[d,"F_1"']
    &
    \A_0\ar[r,"d"]
    \ar[d,"F_0"']
    &
    \A_{-1}\ar[r,"d"]
    \ar[d,"F_{-1}"]
    &
    \dots
    \\
    \dots\ar[r,"d"']
    &
    \B_{n+2}\ar[r,"d"']
    \ar[ur,"\phi_1", shorten >= 5, shorten <= 5 , Leftarrow]
    &
    \B_{n+1}\ar[r,"d"']
    \ar[ur,"\phi_1", shorten >= 5, shorten <= 5 , Leftarrow]
    &
    \B_{n}\ar[r,"d"']
    \ar[ur,"\phi_{-1}", shorten >= 5, shorten <= 5 , Leftarrow]
    &
    \B_{n-1}\ar[r,"d"']
    &
    \dots
  \end{tikzcd}
\end{equation}
which may be called an {\em oplax morphisms} from $\A_{\bullet}$ to $\B_{n+\bullet}$.
The differential then associates to this oplax morphism the (strict) morphism 
\begin{equation}
  \begin{tikzcd}
    \dots\ar[r,"d"]
    &
    \A_2\ar[r,"d"]
    \ar[d,"G_2"']
    &
    \A_1\ar[r,"d"]
    \ar[d,"G_1"']
    &
    \A_0\ar[r,"d"]
    \ar[d,"G_0"']
    &
    \A_{-1}\ar[r,"d"]
    \ar[d,"G_{-1}"']
    &
    \dots
    \\
    \dots\ar[r,"d"']
    &
    \B_{n+1}\ar[r,"d"']
    &
    \B_{n}\ar[r,"d"']
    &
    \B_{n-1}\ar[r,"d"']
    &
    \B_{n-2}\ar[r,"d"']
    &
    \dots
  \end{tikzcd}
\end{equation}
with components $G_i = \fib(\phi_i)[-i]$.

\begin{rem}
  Implicit in the product totalization is a choice of ordering
  of the two directions in the bicomplex \(\Stk(\A_{-\bullet},\B_\bullet)\)
  which corresponds to the choice of direction of the \(2\)-cells
  appearing in the chain maps \eqref{eq:oplax-morphisms}.
  Flipping this choice yields the notion of lax chain morphisms
  and the corresponding lax mapping complex.
\end{rem}

The classical Dold-Kan correspondence offers a means of turning the homological data captured by a connective chain
complex $C_{\bullet}$ into homotopical data described by the Kan complex underlying the associated Dold-Kan nerve of
$C_{\bullet}$. This transformation is of particular interest when applied to mapping complexes where it leads to a model
for the derived category of complexes as a topological category (and further an $\infty$-category by passing to coherent
nerves, cf. \cite{lurie:ha}). 

We may therefore hope to gain insights as to how the derived category of categorical complexes should be defined by
investigating the ``homotopical data'' captured by the categorified Dold--Kan nerve of a categorical complex, and in
particular, of the categorical mapping complex. 
Let $\A_{\bullet}$ and $\B_{\bullet}$ be categorical complexes and let $\M_{\bullet} = \tau_{\ge
0}\Mapoplax{\A_\bullet}{\B_\bullet}\in \Ch(\St)$ be the categorical mapping complex, truncated below and with $\M_{0} =
\ker(d_0)$. Recall (cf. \S \ref{subsec:catdk}), that the categorified Dold--Kan nerve $\X_{\bullet} =
\Nc_{\on{DK}}(\M_{\bullet})$ is a {\em $2$-simplicial stable $\infty$-category}, i.e. a functor
\[
	\X_{\bullet}\colon \DDelta \lra \St
\]
of $(\infty,2)$-categories. We describe its low-dimensional simplices: 

{
\tikzcdset{arrows=Rightarrow}

\begin{itemize}
	\item A vertex corresponds to a morphism $F\colon \A_{\bullet} \to \B_{\bullet}$ of categorical complexes. 
	\item An edge corresponds to the datum of 
		\begin{enumerate}[label=\arabic*.]
			\item a natural transformation 
				\[
					\eta^{(01)}\colon F^{(0)} \Rightarrow F^{(1)},
				\]
				of (strict) morphisms $F^{(0)},F^{(1)}\colon \A_{\bullet} \to \B_{\bullet}$ of categorical
				complexes, 
			\item an oplax morphism
				\[
					H^{(01)}\colon \A_{\bullet} \to \B_{\bullet+1}
				\]
				which we refer to as a {\em categorical homotopy}, 
			\item an exact triangle
				\[
				\begin{tikzcd}
					d H^{(01)} \ar[r] \ar[d] & F^{(0)} \ar[d, "\eta^{(01)}"]\\
					0 \ar[r] & F^{(1)}
				\end{tikzcd}
			\]
		\end{enumerate}
	\item A $2$-simplex corresponds to the datum of 
		\begin{enumerate}[label=\arabic*.]
			\item a diagram 
				\begin{equation}
					\label{eq:fun_diag}
				\begin{tikzcd}
					F^{(0)} \ar[r] &  F^{(1)} \ar[r] &  F^{(2)} 
				\end{tikzcd}
				\end{equation}
				of natural transformations of morphisms $\A_{\bullet} \to \B_{\bullet}$,
			\item a diagram (not necessarily bicartesian)
				\begin{equation}\label{eq:hom_diag}
				\begin{tikzcd}
					H^{(01)} \ar[r]\ar[d]  & H^{(02)}\ar[d] \\
					0 \ar[r]  & H^{(12)}
				\end{tikzcd}
				\end{equation}
				of oplax morphisms $\A_{\bullet} \to \B_{\bullet+1}$,
			\item an oplax morphism
				\[
					H^{(012)}\colon \A_{\bullet} \to \B_{\bullet+2},
				\]
			\item an extension of \eqref{eq:fun_diag} and $d\eqref{eq:hom_diag}$ to a diagram 
				\[
				\begin{tikzcd}
					dH^{(01)} \ar[r]\ar[d]  & dH^{(02)}\ar[d]\ar[r] & F^{(0)}\ar[d]\\
					0 \ar[r]  & dH^{(12)}\ar[d] \ar[r] & F^{(1)}\ar[d]\\
						 & 0 \ar[r] & F^{(2)}\\
				\end{tikzcd}
				\]
				of morphisms $\A_{\bullet} \to \B_{\bullet}$ with all squares biCartesian,
			\item an extension of \eqref{eq:hom_diag} to a biCartesian cube
				\[
				\begin{tikzcd}
					dH^{(012)} \ar[rr]\ar[dr] \ar[dd] & & H^{(01)}\ar[dr]\ar[dd] &  \\
								 & 0 \ar[rr] \ar[dd] & & H^{(02)}\ar[dd] & \\
					0 \ar[rr]\ar[dr] & & 0 \ar[dr] &  \\
								 & 0 \ar[rr] &  & H^{(12)} & \\
				\end{tikzcd}
				\]
				of oplax morphisms $\A_{\bullet} \to \B_{\bullet+1}$.
		\end{enumerate}
\end{itemize}

}

As we observe from this description, the $2$-simplicial object $\M_{\bullet}$ carries meaningful data, such as a
reasonable notion of ``categorical homotopy'' between morphisms of categorical complexes. However, when analyzing
the structural properties of $\M_{\bullet}$ there are many open questions that remain to be investigated in order to
understand the ``higher homotopical content'' of $\M_{\bullet}$ in analogy to its classical counterpart. For
example, the ordinary Dold-Kan nerve is a Kan complex (as every simplicial abelian group) so that it carries intrinsic
homotopical meaning. The categorified Dold--Kan nerve does indeed have certain categorified lifting properties, but
these are much weaker so that, in particular, it is not possible to even ``compose homotopies'' in an obvious way. We hope to analyze this 
context more systematically in future work.

\subsection{Koszul complexes and categorical Koszul duality}
\label{sec:koszul}

Let $R$ be a commutative ring and let $\lambda_1, \ldots, \lambda_n$ be elements of $R$. We may then
form the Koszul complex as the tensor product
\[
	K(\lambda_1, ..., \lambda_n) = \bigotimes_{i=1}^{n} (R \overset{\lambda_i}{\lra} R).
\]
Recall the classical ``self-duality'' of the Koszul complex:

\begin{thm} There is a (canonical) isomorphism of complexes
\[
	K(\lambda_1, ..., \lambda_n)^{\vee}_{n-\bullet} \cong K(\lambda_1, ..., \lambda_n)_{\bullet}
\]
where $(-)^{\vee} = \Hom_R(-, R)$.
\end{thm}

In this section, we explain how the notion of a Koszul complex, as well as the self-duality
statement admit a variant for categorical complexes. 

Let $k$ be a commutative ring spectrum and $\Modk$ the category of $k$-modules. Given an
object $L \in \Modk$ and $\A \in \Stk$, we obtain a functor
\[
	\A \overset{- \otimes_k L}{\lra} \A
\]
which is functorial in $\A$, i.e., it defines a natural transformation on the identity functor of
$\Stk$. We define the categorical two-term complex
\[
	\K(L) \coloneqq  \Modk \overset{- \otimes_k L}{\lra} \Modk,
\]
concentrated in degrees $0$ and $1$, and further, for $k$-module spectra $L_1$, ..., $L_n$, the categorical complex
\[
	\K(L_1, ..., L_n) \coloneqq  \tot^{\amalg}(\K(L_1) \otimes_k \K(L_2) \otimes_k \cdots  \otimes_k
	\K(L_n)).
\]

The main result of this section is the following theorem, which can be regarded as a categorical
variant of the classical ``self-duality'' of the Koszul complex. 

\begin{thm}[Categorical Koszul duality]\label{thm:catkoszul}
	Let $L_1$, ..., $L_n$ be dualizable $k$-module spectra. Then there is a canonical equivalence of
	categorical complexes
	\[
		\K(L_1, ..., L_n)^{\vee}_{n - \bullet} \simeq \K(L_n, ..., L_1)_{\bullet}
	\]
	where $(-)^{\vee} = \Stk(-,\Modk)$ denotes the dual with respect to the symmetric monoidal
	structure on $\Stk$.
\end{thm}
\begin{proof}
	We prove the statement by induction on $n$. For $n=1$, the identification
	\[
		\Stk(\Modk,\Modk) \simeq \Modk
	\]
	via the evaluation functor $F \mapsto F(k)$, extends to an equivalence between the complex
	$\K(L_1)^{\vee}_{1-\bullet}$ and $\K(L_1)$. 

	We now assume given the equivalence
	\[
		\K(L_1,...,L_{n-1})^{\vee}_{n-1-\bullet} \simeq \K(L_{n-1}, ..., L_1)
	\]
	and write $\A_{\bullet} \coloneqq  \K(L_1,...,L_{n-1})$ and $\A'_{\bullet} \coloneqq   \K(L_{n-1}, ..., L_1)$. 
	We have 
	\[
		\K(L_1,L_2,...,L_n) = \tot^{\amalg}(\A_{\bullet} \otimes \K(L_n))
	\]
	We compute 
		\begin{align}
			\tot^{\amalg}(\A_{\bullet} \otimes \K(L_n))^{\vee}_{n-\bullet} 
		&\simeq (\laxpull{\A_{n-i}^{\vee}}{\A_{n-i}^\vee}{\A_{n-i-1}^{\vee}})_{i \in \ZZ}
		\label{kos:1}\\
		&\simeq (\laxpull{\A'_{i-1}}{\A'_{i-1}}{\A'_{i}})_{i \in \ZZ}\label{kos:2}\\
		&\simeq (\laxpush{\A'_{i-1}}{\A'_{i}}{\A'_{i}})_{i \in \ZZ} \label{kos:3}\\
		& \simeq	\tot^{\amalg}(\K(L_n) \otimes \A'_{\bullet}) \label{kos:4}\\
		&\simeq	\K(L_n,...,L_1) \label{kos:5}
		\end{align}
	where
	\begin{itemize}
		\item \eqref{kos:1} holds since the dual of $\tot^{\amalg}$ is equivalent to
			$\tot^{\times}$ of the dual. 
		\item \eqref{kos:2} holds by the induction hypothesis.
		\item \eqref{kos:3} follows from Proposition \ref{prop:beckchevalley} applied to the morphism
			given by tensoring with $L_n$. This has an adjoint given by tensoring with
			$L_n^{\vee}$ (since we assumed the $L_i$ to be dualizable) which is again
			central. 
		\item \eqref{kos:4} holds by definition.
		\item \eqref{kos:5} follows from \Cref{prop:totasso}.
	\end{itemize}
\end{proof}

We give an application of Theorem \ref{thm:catkoszul}: 

\begin{cor}\label{cor:duality} Set $L_i = k$. Then the categorical Koszul complex can be identified
	with the totalization of a constant cube. By Example \ref{exa:higherauslander}, the
	categorical Koszul duality amounts to an equivalence between 
	\begin{itemize}
		\item the category of diagrams 
			\[
				X\colon \Fun([k],[n]) \to \Modk
			\]
			mapping every non-injective map $\tau:[k] \to [n]$ to a zero object in
			$\Mod_k$, and
		\item the category of diagrams 
			\[
				Y\colon \Fun([n-k],[n]) \to \Modk
			\]
			mapping every non-injective map $\tau:[n-k] \to [n]$ to a zero object in
			$\Mod_k$.
	\end{itemize}
\end{cor}

The duality described in Corollary \ref{cor:duality} was first obtained by Beckert \cite{Bec18}
using the theory of derivators. A geometric proof based on Fukaya categories was given in
\cite{DJL}. The novelty in our proof based on categorical Koszul complexes is the inductive nature
which does not feature in the previous proofs (and seems to make the argument both simpler and more
conceptual).

\section{Spherical complexes and perverse schobers}
\label{sec:spherical}

Given an adjunction $F\colon \A\leftrightarrow \B\noloc G$ of stable $\infty$-categories, we define
\begin{itemize}
	\item the {\em twist functor} $T_{\A}:\A\rightarrow \A$ as the cofiber of the unit $u:\id_{\A}\rightarrow FG$ in the stable $\infty$-category $\on{Fun}(\A,\A)$.
	\item the {\em cotwist functor} $T_{\B}:\B\rightarrow \B$ as the fiber of the counit
	$c:FG\rightarrow \id_{\B}$ in the stable $\infty$-category $\on{Fun}(\B,\B)$.
\end{itemize}

\begin{defi}[\cite{AL17, DKSS:spherical}]\label{def:sphadj}
	The adjunction $F \dashv G$ is called {\em spherical} if both $T_{\A}$ and $T_{\B}$ are
	invertible. A functor $F:\A \to \B$ of stable $\infty$-categories is called {\em 
	spherical} if it admits a right adjoint $G$ and the adjunction $F \dashv G$ is spherical.
\end{defi}

Spherical adjunctions were originally conceived by R. Anno, to describe ``family versions'' of the
spherical objects introduced by P. Seidel and R. Thomas in \cite{ST01}. Natural examples
can be found among the various categorical structures that arise within the context of Kontsevich's homological mirror symmetry program.\\

More recently, it has been proposed in \cite{KS14} to interpret spherical functors (and, more
generally, suitable diagram categories built from spherical functors) as categorified analogues of
perverse sheaves (referred to as {\em perverse schobers}). This is motivated by the observation that the abelian category of perverse
sheaves on the complex plane $\CC$, with stratification given by the origin $\{0\}$ and its
complement, is classically known to be equivalent to the category of diagrams 
\begin{equation}\label{eq:ab}
		\begin{tikzcd}
			\Phi \ar[bend left=20]{r}{f} & \Psi \ar[bend left=20]{l}{g},
		\end{tikzcd}
\end{equation}
of vector spaces $\Phi$ and $\Psi$ with $\id - fg$ and $\id - gf$ invertible. 

While a fully satisfying intrinsic definition of perverse schobers still remains open in general and
work in progress in one complex dimension (see \cite{DKS:milnor,DKSS:spherical} for partial results),
in many situations one can guess ad-hoc definitions, based on diagrammatic descriptions of the
respective categories of perverse sheaves such as \eqref{eq:ab}. The resulting notions of perverse
schobers are not intrinsic and depend on auxiliary choices -- however, it is still worthwhile to
construct interesting examples and study applications, since these descriptions
will hopefully become part of a more intrinsic theory in the future.\\

In our current context of categorical complexes, another instance of such an ad-hoc
notion of perverse schober arises from a diagrammatic description of perverse sheaves on $\CC^n$
with respect to the stratification given by the hyperplane inclusions
\[
	\{0\} \subset \CC \subset \CC^2 \subset \CC^3 \subset \cdots \subset \CC^n
\]
and their complements.
The resulting notion of perverse schober in this context is a {\em spherical complex}
(see \S \ref{subsec:spherical}): a categorical complex of length $n$ all of whose differentials are spherical functors. In \S \ref{subsec:sphcub} we also consider the further stratification of $\mathbb{C}^n$ given by the coordinate hyperplanes, their intersections and the complements. In this case, the resulting notion of perverse schober is a Beck-Chevalley categorical $n$-cube, whose edges are described by spherical functors. We call such a categorical cube a \textit{spherical cube}. 

\subsection{Spherical adjunctions and perverse schobers on \texorpdfstring{$\mathbb{C}$}{C}}\label{subsec:sphadj}

As mentioned above, the abelian category of perverse sheaves on $\CC$ with respect to the
stratification given by $\{0\}$ and its complement is equivalent to the category of diagrams
\[
		\begin{tikzcd}
			\Phi \ar[bend left=20]{r}{f} & \Psi \ar[bend left=20]{l}{g},
		\end{tikzcd}
\]
of vector spaces $\Phi$ and $\Psi$ with $\id - fg$ and $\id - gf$ invertible.
Geometrically, the vector spaces $\Phi$ (resp. $\Psi$) correspond to the spaces of {\em vanishing
cycles} (resp. {\em nearby cycles}) associated to a given perverse sheaf. 

\begin{rem}\label{rem:geometric} It is instructive to investigate how this linear algebraic data
	describes the perverse sheaf, when considering it as an object of the derived category of constructible sheaves. To this end, we may also interpret this latter category as the category of
	constructible sheaves on $\CC$ valued in the $\infty$-category of cochain complexes\footnote{When discussing perverse sheaves, we use cochain complexes as this is the standard convention in the literature. In most other places in this paper, we however use chain complexes.} of vector spaces. Such a sheaf $\F$ of cochain complexes may then be described by assembling the linear algebraic data
	into the diagram
	\begin{equation}\label{eq:geometric_perverse}
	\begin{tikzcd}
\Psi \arrow[r, "g"] \arrow[d, "\on{id}"'] & \Phi \arrow[d] \arrow[ld, "f"', dashed] \\
\Psi \arrow[r, "\circlearrowleft_T"']                            & 0                                      
\end{tikzcd}
	\end{equation}
	where
\begin{itemize}
	\item the stalk $\F_0$ of $\F$ at $0 \in \CC$ is given by the cochain complex
		\[
			\begin{tikzcd}
				\Psi \ar{r}{g} & \Phi
			\end{tikzcd}
		\]
		concentrated in degrees $-1$ and $0$, 
	\item the stalk $\F_1$ of $\F$ at $1 \in \CC$ is given by the cochain complex
		\[
			\begin{tikzcd}
				\Psi \ar{r} & 0 
			\end{tikzcd}
		\]
		concentrated in degree $-1$,
	\item the restriction map $\res\colon \F_0 \to \F_1$ corresponds to the morphism represented by
		the commutative square in \eqref{eq:geometric_perverse},
	\item the monodromy of the stalk $\F_1$ is induced by the automorphism $T = \id - fg$ of
		$\Psi$, and finally, 
	\item we may interpret the map $f\colon \Phi \to \Psi$ as a homotopy, as depicted in
		\eqref{eq:geometric_perverse}, between $\res$ and the composite $T \circ \res$. This homotopy expresses the compatibility relations arising from the $\mathbb{C}^{\ast}$-family of restriction maps $\F_0\rightarrow \F_p$, with $p\in \mathbb{C}^{\ast}$, of which we only need to remember $\res\colon \F_0 \to \F_1$.
\end{itemize}
It is further interesting to note that we may apply Verdier duality to describe $\F$ equivalently as
a constructible cosheaf $\F^{\vee}$ valued in cochain complexes. This cosheaf then admits the
following analogous description:
	\begin{equation}\label{eq:geometric_dual}
\begin{tikzcd}
\Phi \arrow[r, "f"]                          & \Psi                                          \\
0 \arrow[r, "\circlearrowleft_T"'] \arrow[u] & \Psi \arrow[u, "\id"'] \arrow[lu, "g", dashed]
\end{tikzcd}
	\end{equation}
\begin{itemize}
	\item the costalk $\F^{\vee}_0$ of $\F^{\vee}$ at $0 \in \CC$ is given by the cochain complex
		\[
			\begin{tikzcd}
				\ar{r}{f} \Phi & \Psi
			\end{tikzcd}
		\]
		concentrated in degrees $0$ and $1$, 
	\item the costalk $\F^{\vee}_1$ of $\F^{\vee}$ at $1 \in \CC$ is given by the cochain complex
		\[
			\begin{tikzcd}
				\ar{r} 0 &  \Psi 
			\end{tikzcd}
		\]
		concentrated in degree $1$,
	\item the corestriction map $\cores\colon \F^{\vee}_1 \to \F^{\vee}_0$ corresponds to the
		commutative square in \eqref{eq:geometric_dual},
	\item the monodromy of the stalk $\F^{\vee}_1$ is induced by the automorphism $T = \id - fg$
	of $\Psi$, and, 
	\item we may interpret the map $g\colon \Psi \to \Phi$ as a homotopy, as depicted in
		\eqref{eq:geometric_dual} between $\cores$ and the composite $\cores \circ T$,
		completing the data needed to define a constructible cosheaf valued in cochain
		complexes. 
\end{itemize}
	Note that the 
	invertibility of $\id - fg$ is in fact equivalent to the invertibility of $\id - gf$,
	explaining why this condition does not appear in an apparent way in the description of $\F$
	or $\F^{\vee}$. 
\end{rem}

\begin{rem}
	The datum of a spherical adjunction 
	\[
		F\colon \A \llra \B \noloc G
	\]
	can be interpreted in a fashion analogous to Remark \ref{rem:geometric}. In other words, we
	may repackage the data into the diagram
	\begin{equation*}
	\begin{tikzcd}
\B \arrow[r, "G"] \arrow[d, "\on{id}"'] & \A \arrow[d] \arrow[ld, "F"', dashed] \\
\B \arrow[r, "\circlearrowleft_T"']                            & 0                                    
\end{tikzcd}
	\end{equation*}
	where the exact triangles
	\[
		\begin{tikzcd}
			T_{\B} \to \id \to FG
		\end{tikzcd}
	\]
	exhibit $F$ as a categorical homotopy, as defined in \Cref{sub:mapping_complexes}. In the categorical
	context, however, it is not the case, that the invertibility of $T_{\B}$ implies the
	invertibility of $T_{\A}$ so that this perspective only partially characterizes spherical
	adjunctions. Nevertheless, it motivates our proposal to interpret spherical functors as
	$2$-term categorical complexes and use this intuition to generalize to {\em spherical complexes} as
	explained below in \S \ref{subsec:spherical}.
\end{rem}

\subsection{Barbacovi's construction and its geometric interpretation}\label{subsec:Barbacovi}

A remarkable result by Ed Segal (\cite{segal:twists}) says that every autoequivalence on a category
arises as a spherical twist functor (cf.~also \cite{Chr20} for an $\infty$-categorical generalization). In particular, the composition of the spherical twists associated to a
pair of spherical adjunctions with the same target must again be a spherical twist. 
The main result of \cite{Bar20} provides an explicit description of a spherical adjunction that
describes this composite twist in terms of the given spherical adjunctions. 

In this section, we explain how Barbacovi's construction fits into our context, provide a geometric
interpretation in terms of perverse schobers, and finally prove a generalization to stable
$\infty$-categories which we will need below.

To begin with, consider the complex plane $\CC$ equipped the stratification given by two points
$\{x,x'\}$, say $x = i$ and $x' = -i$, and their complement. By choosing an arc connecting the two points, we may identify
the category of perverse sheaves with the category of diagrams
\[
	\begin{tikzcd}
		\Phi \ar[bend left=20]{dr}{f} &  \\
			& \Psi \ar[bend left=20]{dl}{g'}\ar[bend left=20]{ul}{g}\\
		\Phi' \ar[bend left=20]{ur}{f'}& 
	\end{tikzcd}
\]
of vector spaces such that both $(f,g)$ and $(f',g')$ satisfy the conditions from \eqref{eq:ab}.
Geometrically, the automorphisms $T = \id - fg$ and $T' = \id -f'g'$ can be interpreted as the
monodromy transformations on the space of nearby cycles $\Psi$ corresponding to loops around the
points $x$ and $x'$, respectively.

As mentioned above, the vector spaces $\Phi$ and $\Phi'$ can be interpreted as the {\em local} vanishing cycles of the given perverse sheaf $\F$. We define the space of {\em global} vanishing cycles
$\Sigma(\F)$ as the fiber 
\[
	\fib(\R\Gamma(\CC,\F) \to \R\Gamma(\{\Re z < 0\}, \F)) \cong \Phi \oplus \Phi'
\]
of the restriction map along the inclusion of the half-plane  $\{\Re z < 0\} \subset \CC$. 
This vector space comes equipped with natural maps
\begin{equation}\label{eq:amalgamate}
\begin{tikzcd}[column sep=huge]
\Sigma(\F) \arrow[r, bend left=15, "{(f,f')}"] & \Psi \arrow[l, bend left=15, "{(gT',g')}"]
\end{tikzcd}
\end{equation}
and we compute
\[
	\id - (f,f') \begin{pmatrix} gT'\\g'\end{pmatrix} = \id - fgT' - f'g' = (\id - fg)
	(\id - f'g') = T T'.
\]
Therefore, the datum \eqref{eq:amalgamate} defines a perverse sheaf on $\CC$ with respect to the
startification given by $\{0\}$ and its complement. We refer to this sheaf as the {\em amalgamate} of $\F$. This amalgamate likely describes the pushforward sheaf of $\F$ under the endomorphism of $\CC$, which contracts a disc containing $x$ and $x'$ to $0$, but we do not verify this here. 

As we will now explain, these geometric considerations can be lifted to analogous constructions for
perverse schobers. To begin, a perverse schober on $\CC$, equipped with the above
stratification and auxiliary choices, can be interpreted as a pair
\[
	\begin{tikzcd}
		\A \ar[bend left=20]{dr}{F} &  \\
			& \B \ar[bend left=20]{dl}{G'}\ar[bend left=20]{ul}{G}\\
		\A' \ar[bend left=20]{ur}{F'}& 
	\end{tikzcd}
\]
of spherical adjunctions as indicated. As a categorical lift of \eqref{eq:amalgamate}, we define the functor
\begin{equation}\label{eq:fbf'} 
	\laxpull{F}{\B}{F'}:\laxpull{\A}{\B}{\A'}\rightarrow \B\,.
\end{equation}

\begin{thm}\label{thm:compsph}
	The functor $\laxpull{F}{\B}{F'}$ is spherical. The
	suspension of the cotwist functor of the adjunction 
	\[
		\laxpull{F}{\B}{F'} \dashv \ra{\left(\laxpull{F}{\B}{F'}\right)}
	\]
	is equivalent to the composite of the cotwist functors of the adjunctions $F'\dashv G'$ and $F\dashv G$. 
\end{thm}

In the following, we informally describe how Theorem \ref{thm:compsph} can be proved using further
ideas related to perverse schobers. The full and somewhat technical proof implementing these ideas
is given in \S \ref{subsec:proofBarbacovi}.

Given a perverse sheaf $\F$ on $(\mathbb{C},0)$, instead of considering the vector spaces of vanishing and nearby cycles, it turns out that we can equally well describe $\F$ in terms
of its nearby cycles $\Psi$ and its vector space of global sections $\Phi_2$ with support on $\mathbb{R}\subset \mathbb{C}$ or even its vector space of global sections $\Phi_n$ with support on any graph embedded in $\mathbb{C}$ with a single $n$-valent vertex at $0$. The resulting description
states that the category of perverse sheaves on $(\mathbb{C},0)$ is equivalent to the category of
diagrams of vector spaces
\begin{equation}\label{eq:psloc}
	f_i\colon \Phi_n\longleftrightarrow \Psi_i\noloc g_i\,,\quad 1\leq i \leq n 
\end{equation}
satisfying that $f_ig_i=\on{id}_{\Psi_i}$, $f_{i+1}g_i$ is invertible and $f_jg_i=0$ for $j\neq i,i+1$,
with $i,j$ considered modulo $n$, see \cite{KS16}. The vector spaces $\Psi_i$, with $1\leq i\leq n$, are all equivalent and may be chosen for concreteness as the stalks of the $n$-th roots of unity. An ad-hoc categorification of this local description is
described in \cite{Chr22}, based on Waldhausen's relative $S_\bullet$-construction. Using these local descriptions, we
can describe perverse sheaves \cite{KS16} or perverse schobers \cite{Chr22} on any surface ${\bf S}$ with $0$-dimensional
strata $P$ and non-empty boundary given the auxiliary choice of a spanning graph ${\bf G}\subset {\bf S}$. This means that the inclusion of ${\bf G}$ into ${\bf S}$ is required to be a homotopy equivalence and each stratum
$p\in P$ is required to be the image of a vertex of ${\bf G}$. A perverse sheaf or perverse schober on ${\bf
S}$ is encoded as a constructible sheaf and cosheaf on ${\bf G}$, which restricts at each vertex of ${\bf G}$ to a
diagram as in \eqref{eq:psloc} or its categorification. The global sections of the
constructible sheaf on ${\bf G}$ describe the first cohomology of the derived global sections with support on
${\bf G}$ of the perverse sheaf or perverse schober on ${\bf S}$. 

To prove \Cref{thm:compsph}, we consider the perverse schober $\F$ on $(\mathbb{C},\{x,x'\})$ defined by the
pair of spherical adjunctions $F\dashv G$ and $F'\dashv G'$ and describe it as a constructible sheaf on the ribbon graph ${\bf G}$ with two vertices at $x,x'$, depicted as follows (in blue):
\begin{center}
\begin{tikzpicture}
\draw[color=blue,very thick] (-3,0)--(3,0);
\draw[color=black, very thick, dotted] (0,0) ellipse (3cm and 2cm);
\node() at (-0.6,0.35){$x$};
\fill[color=blue] (-0.6,0) circle (0.1cm);
\node() at (0.6,0.35){$x'$};
\fill[color=blue] (0.6,0) circle (0.1cm);
\fill[color=gray] (2.3,0) circle (0.1cm);
\fill[color=gray] (-2.3,0) circle (0.1cm);
\node() at (1.4,0.3){${\bf G}$};
\node() at (2.3,0.35){$q$};
\node() at (-2.3,0.35){$p$};
\node() at (2.5,1.8){$\mathbb{C}$};
\end{tikzpicture}
\end{center}
The points $p$ and $q$ lie left or right of $x$ and $x'$ as indicated, but are otherwise arbitrary. We denote the global sections of $\F$ by
$\R\Gamma^1({\bf G},\F)\in \Stk$. Restriction to the two points $p,q$ on ${\bf G}$ defines a functor 
\[(\on{ev}_{p},\on{ev}_q)\colon \R\Gamma^1({\bf G},\F)\rightarrow \B^{\times 2}\,,\]
The right adjoint is denoted 
$(\ra{\on{ev}_{p}},\ra{\on{ev}_q})$ and the left adjoint is denoted $(\la{\on{ev}_{p}},\la{\on{ev}_q})$. These functors turn out to be much easier to describe than the functor \eqref{eq:fbf'}. The functors $\ra{\on{ev}_p},\ra{\on{ev}_q},\la{\on{ev}_p},\la{\on{ev}_q}$ are all fully faithful and, as one can compute, there are equivalences $\ra{\on{ev}_q}\simeq \la{\on{ev}_p}\circ T$ and $\ra{\on{ev}_p}\simeq \la{\on{ev}_q}\circ T'$, with $T$ the cotwist functor of $F\dashv G$ and $T'$ the cotwist functor of $F'\dashv G'$. Using these facts, we will deduce that that $(\on{ev}_{p},\on{ev}_{q})$ is a spherical functor. 

To show the sphericalness of the functor $\laxpull{F}{\B}{F'}$, we use the fact that this is equivalent to the
assertion that there exist a stable $\infty$-category with a $4$-periodic semiorthogonal
decomposition $(\laxpull{\A}{\B}{\A'},\B)$ with gluing functor $\laxpull{F}{\B}{F'}$, see
\cite{HLS16,DKSS:spherical}. Remarkably, we can find this $4$-periodic semiorthogonal decomposition as a semiorthogonal decomposition of the $\infty$-category $\R\Gamma^1({\bf G},\F)$. Here, $\laxpull{\A}{\B}{\A'}\subset
\R\Gamma^1({\bf G},\F)$ arises as the full subcategory of global sections with support on
${\bf G}$ which vanish at $p$, i.e.~the kernel of
$\on{ev}_p:\R\Gamma^1({\bf G},\F)\rightarrow \B$. The $\infty$-category $\B\subset
\R\Gamma^1({\bf G},\F)$ arises as the image of $\ra{\on{ev}_p}$. 

Theorem \ref{thm:compsph} is  a special case of a more general phenomenon exhibited by perverse
schobers on a stratified surface ${\bf S}$ with non-empty boundary.  We again choose a spanning
graph ${\bf G}$ of ${\bf S}$. By a ${\bf G}$-parametrized perverse schober, we mean a constructible sheaf of stable $\infty$-categories on ${\bf G}$ encoding a perverse schober on ${\bf S}$, as explained above and defined in \cite{Chr22}. Given an edge $e$ of ${\bf G}$, we denote by $\F(e)$ the stalk of $\F$ at any point on that edge.

\begin{thm}\label{thm:schoberbdry}
Let $\F$ be a ${\bf G}$-parametrized perverse schober on ${\bf S}$. Let $E^\partial$ be the set of
external edges of ${\bf G}$ and \[ \prod_{e\in E^\partial}
\on{ev}_e:\R\Gamma^1({\bf G},\F)\rightarrow \prod_{e\in E^\partial}\F(e)\] the restriction
functor of global sections with support on ${\bf G}$ to the stalks at these boundary edges. The functor
$\prod_{e\in E^\partial} \on{ev}_e$ is spherical.   
\end{thm}

In the study of partially wrapped Fukaya categories, the functor $\prod_{e\in E^\partial} \on{ev}_e$ is also
called the cap functor, it is adjoint to the Orlov or cup fuctor, see for instance \cite{Syl19}. We give some hints as to how the proof of \Cref{thm:compsph} generalizes to a proof of \Cref{thm:schoberbdry} in \S\ref{subsec:proofBarbacovi}.

\subsection{Spherical complexes and perverse schobers on \texorpdfstring{$\mathbb{C}^n$}{Cn}}\label{subsec:spherical}

Let $n \ge 1$. We begin with a linear algebraic description of classical perverse sheaves on
$\CC^{n}$ with respect to the stratification given by the hyperplane inclusions
\begin{equation}\label{eq:perverse-cn}
	\{0\} \subset \CC \subset \CC^2 \subset \CC^3 \subset \cdots \subset \CC^n\,,
\end{equation}
each setting the last coordinate to zero.

\begin{thm}\label{thm:perverse-cn}
	The category of perverse sheaves on $\CC^n$ with respect to the stratification
	\eqref{eq:perverse-cn} is equivalent to the category of diagrams
	\begin{equation}\label{eq:perverse_cnA}
		\begin{tikzcd}
			A_0 \ar[bend left=20]{r}{d} & A_1 \ar[bend left=20]{l}{\delta}\ar[bend
			left=20]{r}{d} & A_2  \ar[bend left=20]{l}{\delta} & ... & 
			A_{n-1} \ar[bend left=20]{r}{d} & \ar[bend left=20]{l}{\delta} A_n
		\end{tikzcd}
	\end{equation}
	of vector spaces subject to the conditions
	\begin{enumerate}[label=\arabic*.]
		\item $d^2 = 0$,
		\item $\delta^2 = 0$,
		\item for every $0 \le k \le n$, the endomorphisms $\id - d \delta$ and $\id -
			\delta d$ of $A_k$ are invertible. 
	\end{enumerate}
\end{thm}

\begin{proof} 
The stratification of $\CC^n$ given by the coordinate hyperplanes, their intersections and the complements is a refinement of the stratification \eqref{eq:perverse-cn}. By \Cref{thm:perverse-cube} below, a perverse sheaf on $\CC^n$ with the former stratification amounts to a certain cube of linear maps. One readily finds that such a perverse sheaf also defines a perverse sheaf with respect to the stratification \eqref{eq:perverse-cn} if and only if all entries of this cube vanish, except for a sequence of entries as in \eqref{eq:perverse_cnA}. 
\end{proof}

Alternatively, \Cref{thm:perverse-cn} can also be directly deduced from iterated application of
Beilinson's gluing formula for categories of perverse sheaves (\cite{beilinson:howtoglue}). 

\begin{rem}\label{rem:geometric-cn} We explain how to interpret the linear algebraic data from
	Theorem \ref{thm:perverse-cn} geometrically, in terms of the corresponding perverse sheaf as
	a constructible sheaf $\F$ valued in the $\infty$-category of cochain complexes of
	vector spaces, in analogy to Remark \ref{rem:geometric}. As
	in Remark \ref{rem:geometric}, we assemble the linear algebraic data into a diagram
	\begin{equation}\label{eq:geometric_perverse-cn}
\begin{tikzcd}
                             & \A_n \arrow[r, "\delta"] \arrow[d] & \A_{n-1} \arrow[d] \arrow[ld, "d"', dashed] \arrow[r] & \dots \arrow[ld, dashed] & \A_2 \arrow[d] \arrow[r, "\delta"] & \A_1 \arrow[d] \arrow[r, "\delta"] \arrow[ld, "d"', dashed] & \A_0 \arrow[d] \arrow[ld, "d"', dashed] \\
{\circlearrowleft_{T_{1}}}   & \A_n \arrow[r, "\delta"] \arrow[d] & \A_{n-1} \arrow[d] \arrow[ld, "d"', dashed] \arrow[r] & \dots \arrow[ld, dashed] & \A_2 \arrow[d] \arrow[r, "\delta"] & \A_1 \arrow[d] \arrow[r] \arrow[ld, "d"', dashed]           & 0 \arrow[d]                             \\
{\circlearrowleft_{T_{2}}}   & \A_n \arrow[r, "\delta"] \arrow[d] & \A_{n-1} \arrow[d] \arrow[r] \arrow[ld, dashed]       & \dots                    & \A_2 \arrow[r]                     & 0 \arrow[r]                                                 & 0                                       \\
\vdots                       & \vdots                             & \vdots                                                & \ddots                   & \vdots \arrow[d]                   & \vdots \arrow[d]                                            & \vdots \arrow[d]                        \\
{\circlearrowleft_{T_{n-1}}} & \A_n \arrow[r, "\delta"] \arrow[d] & \A_{n-1} \arrow[d] \arrow[ld, "d"', dashed]           & \dots \arrow[r]          & 0 \arrow[d] \arrow[r]              & 0 \arrow[d] \arrow[r]                                       & 0 \arrow[d]                             \\
{\circlearrowleft_{T_n}}     & \A_n \arrow[r]                     & 0                                                     & \dots \arrow[r]          & 0 \arrow[r]                        & 0 \arrow[r]                                                 & 0                                      
\end{tikzcd}
	\end{equation}
where
\begin{itemize}
	\item the rows of \eqref{eq:geometric_perverse-cn} correspond to the stalks of $\F_i$ at the points $x^0,x^1, x^2, ..., x^n \in \CC^n$ with
		\[
			(x^i)_j = \begin{cases} 
					1 & \text{for $j \leq i$}\\
					0 & \text{for $j > i$}
				  \end{cases}
		\]
		where the part of the complex depicted by the corresponding row of
		\eqref{eq:geometric_perverse-cn} is concentrated in degrees $-n,\dots,-i$. Thus, $\A_i$ lies in degree $-i$. 
	\item the restriction maps $\res_i:\F_{i-1} \to \F_{i}$, $1\leq i\leq n$, correspond to the commutative rectangles
		between the corresponding rows in \eqref{eq:geometric_perverse-cn} where all
		vertical maps are either $0$ or $\id$,
	\item we note without proof that the monodromies of each stalk about the ``previous'' hyperplane in
		\eqref{eq:perverse-cn} are induced by the chain automorphism $T_i = \id_\bullet - d_\bullet \delta_\bullet -
		\delta_\bullet d_\bullet$ of $\F_i$, 
	\item we may interpret the maps $d\colon A_\ast \to A_{\ast+1}$ as defining homotopies, as depicted in
		\eqref{eq:geometric_perverse-cn}, between $\res_i$ and the composite $T_i \circ \res_i$.
\end{itemize}
\end{rem}

\begin{rem}\label{rem:sphericalhomotopy}
	It is interesting to note here that the invertibility of the map $T = \id - d \delta -
	\delta d$ is equivalent to the invertibility of the maps $\id - d \delta$ and $\id - \delta
	d$ by virtue of the formula
	\begin{equation}\label{eq:idddelta}
		\id - d \delta - \delta d = (\id - d \delta) (\id - \delta d) = (\id - \delta d) (\id
		- d\delta)
	\end{equation}
	We leave the Verdier dual interpretation of the linear algebraic data in analogy to Remark
	\ref{rem:geometric} to the reader, only noting that the roles of $d$ and $\delta$ get
	swapped. This gives an explanation why both $d$ and $\delta$ need to square to $0$. \Cref{thm:perverse-cn} shows that this data gives a full description of $\F$. 
\end{rem}

Inspired by Theorem \ref{thm:perverse-cn}, we introduce the following concept of spherical categorical complexes. 

\begin{defi}
	A categorical complex $\A_{\bullet}\in \on{Ch}(\St_k)$ is called {\em spherical} if the
	differential $d:\A_i\rightarrow\A_{i-1}$ is a spherical functor for all $i\in \mathbb{Z}$.
\end{defi}

Spherical categorical complexes concentrated in degrees $n,\dots,0$ can thus be regarded as perverse
schobers on $\CC^n$ with respect to the stratification \eqref{eq:perverse-cn}. We conclude this section with some comments on the twist functors of spherical categorical complexes.

We fix a categorical complex $A_{\bullet}\in \on{Ch}(\St_k)$ and $i\in \mathbb{Z}$. The adjunctions
$d_i\dashv \delta_i$ and $d_{i+1}\dashv \delta_{i+1}$ give rise to the unit
$u:\on{id}_{\A_i}\rightarrow \delta_id_i$ and the counit
$c:d_{i+1}\delta_{i+1}\rightarrow \on{id}_{\A_i}$. We can compose these two natural
transformations to obtain the following two commutative diagrams in $\Stk(\A_i,\A_i)$:
\begin{equation}\label{eq:tottw}
\begin{tikzcd}
d_{i+1}\delta_{i+1} \arrow[r, "\on{cu}"] \arrow[d] & \on{id}_{\A_i} \arrow[d, "\unit"] \\
d_{i+1}\delta_{i+1}\delta_id_i \arrow[r]           & \delta_id_i                       
\end{tikzcd}\quad \quad\quad
\begin{tikzcd}
d_{i+1}\delta_{i+1} \arrow[r, "\on{cu}"] \arrow[d] & \on{id}_{\A_i} \arrow[d, "\on{u}"] \\
\delta_id_id_{i+1}\delta_{i+1} \arrow[r]           & \delta_id_i                       
\end{tikzcd}
\end{equation}
Note that $d_{i+1}\delta_{i+1}\delta_id_i\simeq 0$ and $\delta_id_id_{i+1}\delta_{i+1}\simeq 0$ since $\delta^2\simeq 0$ and $d^2\simeq 0$. 

\begin{lem}
Let $T_i$ be the twist of $d_i\dashv \delta_i$ and $T_i'$ the cotwist of $d_{i+1}\dashv \delta_{i+1}$. The totalization of the left square in \eqref{eq:tottw} is equivalent to $T_i'T_i$ and the totalization of the right square in \eqref{eq:tottw} is equivalent to $T_iT_i'$. 
\end{lem}

\begin{proof}
Immediate. 
\end{proof}

The totalization of the left diagram in \eqref{eq:tottw} categorifies the expression $\on{id}-d\delta-\delta d$. Since the totalization is $T_i'T_i$, we find a categorification of the expression $\on{id}-T_i'T_i=d\delta+\delta d$, expressing that $\delta$ describes a homotopy between the chain maps $T_\bullet'T_\bullet$ and $\on{id}_{A_\bullet}$. Similarly, the right diagram in \eqref{eq:tottw} expresses that $\delta$ is a homotopy between $\on{id}_{\A_\bullet}$ and $T_\bullet T_{\bullet}'$. The following lemma shows that $T_\bullet T_{\bullet}'\simeq T_\bullet' T_{\bullet}$, categorifying the identity \eqref{eq:idddelta}.

\begin{lem}
Suppose that $\A_\bullet$ is spherical. The two diagrams in \eqref{eq:tottw} are equivalent. In
particular, the resulting equivalence of their totalizations shows that the autoequivalences $T_i$ and
$T_i'$ commute.
\end{lem}

\begin{proof}
Using that composition with an exact functor defines an exact functor between stable $\infty$-categories of functors, we find that the two commutative squares in $\Stk(\A_i,\A_i)$ can be extended to commutative diagrams where the horizontal sequences are fiber and cofiber sequences:
\[
\begin{tikzcd}
{T_i'[-1]} \arrow[r] \arrow[d, "{T_i'[-1]\on{u}}"] & d_{i+1}\delta_{i+1} \arrow[r, "\on{cu}"] \arrow[d] & \on{id}_{\A_i} \arrow[d, "\on{u}"] \\
{T_i'\delta_id_i[-1]} \arrow[r]                    & d_{i+1}\delta_{i+1}\delta_id_i \arrow[r]           & \delta_id_i                       
\end{tikzcd}
\quad\quad\quad
\begin{tikzcd}
{T_i'[-1]} \arrow[r] \arrow[d, "{\on{u}T_i'[-1]}"] & d_{i+1}\delta_{i+1} \arrow[r, "\on{cu}"] \arrow[d] & \on{id}_{\A_i} \arrow[d, "\on{u}"] \\
{\delta_id_iT_i'[-1]} \arrow[r]                    & \delta_id_id_{i+1}\delta_{i+1} \arrow[r]           & \delta_id_i                       
\end{tikzcd}
\]
We begin by showing that the morphism $T_i'[-1]\unit$ and $\unit T_i'[-1]$ in the upper diagrams are equivalent. Consider the adjunction $T_i'\delta_i\dashv d_i(T_i')^{-1}$. Its unit is given by $T_i'\unit (T_i')^{-1}$. Using that $T_i'\delta_i\simeq \fib(d_{i+1}\delta_{i+1}\delta_i\rightarrow \delta_i)$ and $\delta^2\simeq 0$, we get $T_i'\delta_i\simeq \delta_i[-1]$ and hence also $d_i(T_i')^{-1}\simeq d_i[1]$. It follows that the unit $T_i'\unit (T_i')^{-1}$ of $T_i'\delta_i\dashv d_i(T_i')^{-1}$ is equivalent to the unit of $\delta_i[-1]\dashv d_i[1]$, which is equivalent to $\unit$. We thus obtain $T_i'[-1]\unit \simeq \unit T_i'[-1]$.

Since $d_{i+1}\delta_{i+1}\delta_id_i\simeq \delta_id_id_{i+1}\delta_{i+1}\simeq 0$, it is now clear that the upper left squares in the above two diagrams are equivalent. Hence the entire diagrams (which are recovered as cofibers) are also equivalent.
\end{proof}

\subsection{Spherical cubes and perverse schobers on \texorpdfstring{$\mathbb{C}^n$}{Cn}}\label{subsec:sphcub}

We begin by recalling another classical description of the category of perverse sheaves on the stratified space $\CC^n$ with strata given by the coordinate hyperplanes and their (iterated) intersections and complements. Let $[[1]]$ be the $1$-category with two objects $0,1$ and morphisms freely generated by two morphisms $1\rightarrow 0$ and $0\rightarrow 1$. We identify the set of objects of $[[1]]^n$ with the set $\Pc([n])^\op$, by identifying $J\subset [n]$ with its characteristic function in $[[1]]^n$. By a cubical double diagram, we mean a functor $[[1]]^n\to \on{Vect}_k$, which amounts to the datum of
\begin{itemize}
\item a vector space $V_J$ assigned to each $J\in \Pc([n])^\op$,
\item a pair of linear maps $f_i:V_{J\cup \{i\}}\leftrightarrow V_{J}:g_i$ assigned to each $J\in \Pc([n])^\op$ and $i\in [n]\backslash J$,
\item satisfying $f_if_j=f_jf_i$, $g_ig_j=g_jg_i$ and $g_if_j=f_jg_i$ for all $i,j\in [n]$, $i\neq j$, whenever these composites are defined.
\end{itemize}
With this terminology, we may then formulate the following classical description of perverse sheaves on
$\CC^n$:

\begin{thm}[{\cite{GMP85}}]\label{thm:perverse-cube}
	The category of perverse sheaves on $\CC^n$ with respect to the coordinate hyperplane 
	stratification is equivalent to the category of cubical double diagrams 
	\begin{equation}\label{eq:double_cube}
			[[1]]^n \to \on{Vect}_k
	\end{equation}
satisfying that, for all pairs of maps
			$f_i\colon V_{I\cup \{i\}} \leftrightarrow V_I\noloc g_i$, the endomorphisms 
			\[ g_if_i-\on{id}_{V_{I\cup \{i\}}}\quad\text{and}\quad f_ig_i-\on{id}_{V_{I}}\] are invertible. 
\end{thm}

As in the $1$-dimensional case, we ask the pairs of linear maps $(f_i,g_i)$ to become pairs of
adjoint functors upon categorification. Remarkably, the commutativity conditions $g_if_j=f_jg_i$
then correspond to the Beck--Chevalley conditions already introduced and studied in \S
\ref{subsec:beck-chevalley} (where the motivation for introducing them came from their effect on
totalization). Inspired by Theorem \ref{thm:perverse-cube}, we thus introduce the following:

\begin{defi} A Beck-Chevalley categorical $n$-cube $\A_\ast:\I^n \to \Stk$ is called {\em
	spherical} if every rectilinear edge is a spherical functor.
\end{defi}

A natural categorification of a perverse sheaf on $\mathbb{C}^n$ with the above stratification is thus given by a spherical categorical $n$-cube. We conclude this section, by showing that totalization takes spherical categorical cubes to spherical categorical complexes.

\begin{prop}\label{prop:sphtot}~
\begin{enumerate}[(1)]
\item Consider a Beck-Chevalley chain map $F:\A_\bullet\rightarrow \B_{\bullet}$ between spherical categorical complexes, satisfying that $F_i:\A_i\rightarrow \B_i$ is a spherical functor for all $i\in \mathbb{Z}$. Then $\Fib(F)$ is a spherical complex.
\item Consider a spherical categorical $n$-cube $\A_\ast$. Then the product totalization $\on{tot}^{\times}(\A_\ast)$ and the coproduct totalization $\on{tot}^{\amalg}(\A_\ast)$ are equivalent and spherical categorical complexes.  
\end{enumerate}
\end{prop}

\begin{proof}
We begin by proving part (1). We depict $F$ as follows:
\[
\begin{tikzcd}
\dots \arrow[r] & \A_2 \arrow[r, "d^{\A}_2"] \arrow[d, "F_2"] & \A_1 \arrow[d, "F_1"] \arrow[r, "d^\A_1"] & \A_0 \arrow[d, "F_0"] \arrow[r] & \dots \\
\dots \arrow[r] & \B_2 \arrow[r, "d^\B_2"]                    & \B_1 \arrow[r, "d^{\B}_1"]                & \B_0 \arrow[r]                  & \dots
\end{tikzcd}
\]
Using \Cref{thm:compsph}, we obtain a spherical functor 
\[
\laxpull{F_1}{\B_1}{d_2^{\B}}:\laxpull{\A_1}{\B_1}{\B_2}\longrightarrow \B_1\,.
\]
The functor 
\[
\tilde{d}_1^{\A}\coloneqq \laxpull{d_1^{\A}}{\A_1}{0}:\laxpull{\A_1}{\B_1}{\B_2}\xlongrightarrow{\pi_{\A_1}} \A_1\xlongrightarrow{d_1^{\A}}\A_0\,,
\]
where $\pi_{\A_1}$ denotes the left adjoint of the inclusion of $\A_1$, is furthermore also
spherical. To see this, we first note that it is immediate by the fully faithfulness of
$\ra{(\pi_{\A_1})}$ that the cotwist functor of the adjunction $\tilde{d}_1^{\A}\dashv
\ra{(\tilde{d}_1^{\A})}$ is equivalent to the cotwist functor of the adjunction $d_1^{\A}\dashv
\ra{(d_1^{\A})}$. To show that the twist is an autoequivalence, we apply \cite[Lemma
4.13]{dyck:dk}. Clearly, the twist acts on $\A_1\subset \displaystyle \laxpull{\A_1}{\B_1}{\B_2}$ as the twist of $d_1^{\A}\dashv \ra{(d_1^{\A})}$ and is thus invertible. On $\B_2\subset \displaystyle \laxpull{\A_1}{\B_1}{\B_2}$, the twist clearly acts as the suspension functor, which is also invertible. It thus remains to show that the twist preserves Cartesian edges, which are of the form $a\xrightarrow{\ast} \la{(d_2^{\B})}\circ F_1(a)$. By the sphericalness of $d^{\B}_2$
we have $\la{(d^{\B}_2)}\simeq T\circ \ra{(d^{\B}_2)}$ with $T$ the twist functor of $\la{d^{\B}_2}\dashv d^{\B}_2$. The Beck-Chevalley condition and $d^2=0$ thus imply
\begin{equation}\label{eq:d20forCart} \la{(d_2^{\B})}\circ F_1\circ \ra{(d_1^{\A})}\circ d_1^{\A}\simeq 0\,.\end{equation}
We get a cofiber sequence in $\displaystyle \laxpull{\A_1}{\B_1}{\B_2}$
\[
\begin{tikzcd}
a \arrow[d, "\ast"] \arrow[r]                         & \ra{(d_1^{\A})}\circ d_1^{\A}(a) \arrow[r] \arrow[d, "\ast"] & T(a) \arrow[d, "\ast"]                           \\
\la{(d_2^{\B})}\circ F_1(a) \arrow[r] & 0 \arrow[r]                                                   & {\la{(d_2^{\B})}\circ F_1(a)[1]}
\end{tikzcd}
\]
where \eqref{eq:d20forCart} implies that the middle vertical arrow is Cartesian and the last vertical morphism is Cartesian as the cofiber of Cartesian morphisms. This shows that the twist preserves Cartesian edge and is hence an equivalence. 

The functor 
\[
\mathcal{A}_0\xlongrightarrow{\ra{(\tilde{d}_1^{\A})}} \laxpull{\A_1}{\B_1}{\B_2}\xlongrightarrow{\laxpull{F_1}{\B_1}{d_2^{\B}}} \B_1
\]
is clearly equivalent to $F_1\circ \ra{(d_1^{\A})}$. It follows that $\laxpush{\A_0}{\laxpull{\A_1}{\B_1}{\B_2}}{\B_1}\simeq \laxpush{\A_0}{\A_1}{\B_1}$ and the Beck-Chevalley property further implies that $\laxpush{\A_0}{\A_1}{\B_1}\simeq \laxpull{\A_0}{\B_0}{\B_1}$. We find that the differential $d_1$ of $\Fib(F)$ is equivalent to the the composite of these equivalences with the functor
\[
\laxpush{\tilde{d}_1^{\A}}{\displaystyle \laxpull{\A_1}{\B_1}{\B_2}}{\left(\laxpull{F_1}{\B_1}{d_2^{\B}}\right)}:\laxpull{\A_1}{\B_1}{\B_2}\longrightarrow \laxpush{\A_0}{\displaystyle \laxpull{\A_1}{\B_1}{\B_2}}{\B_1}\,.
\]
This functor is again spherical by \Cref{thm:compsph}. The same argument applied to each degree shows that all differentials of $\Fib(F)$ are spherical functors, meaning that $\Fib(F)$ is a spherical categorical complex. This concludes the proof of part (1).

For part (2), we begin by noting that the equivalence of the product and coproduct totalizations follows from repeated application of \Cref{prop:beckchevalley}. We prove the sphericalness of the product totalization via an induction on $n$. The case $n=2$ follows from part (1). Fix $n\geq 3$. We can consider the categorical $n$-cube $\A_\ast$ as a morphism between two categorical $(n-1)$-cubes $\B_\ast\rightarrow \B_\ast'$, whose totalizations are spherical categorical complexes by the induction assumption. The arising morphism $\beta:\on{tot}^{\times}(\B_\ast)\rightarrow \on{tot}^{\times}(B_\ast')$ is Beck-Chevalley by \Cref{lem:parBC}. Repeated application of \Cref{lem:laxsph} shows that $\beta$ is a spherical functor in each degree. The sphericalness of $\on{tot}^{\times}(\A_\ast)$ thus follows again from part (1) applied to the morphism of categorical complexes $\beta$, concluding the proof.
\end{proof}

\begin{lem}\label{lem:laxsph}
Consider a commutative Beck-Chevalley diagram in $\St_k$ of the following form.
\[ 
\begin{tikzcd}
\A \arrow[r, "F"] \arrow[d, "\alpha"] & \C \arrow[d, "\alpha'"] \\
\B \arrow[r, "F'"]               & \D                
\end{tikzcd}
\]
If $F$ and $F'$ are spherical functors, then the induced functor 
\[ F\overset{\rightarrow}{\times}{F'}:\laxlima{\A}{\alpha}{\B}\longrightarrow \laxlima{\C}{\alpha'}{\D}\]
is also spherical.
\end{lem}

\begin{proof}
  Let $G$ and $G'$ be the right adjoints of $F$ and $F'$, respectively. Let $T_\mathcal{A}$ and $T_{\mathcal{B}}$ be the twist functors of $F\dashv G$ and $F'\dashv G'$, respectively.
  Adjoints on lax limits are determined componentwise,
  see \cite[Appendix~A.2]{CDW24}.
  We thus have $F\overset{\rightarrow}{\times}{F'}\dashv G\overset{\rightarrow}{\times}{G'}$ and the twist functor of this adjunction can be identified with the induced functor $T_{\mathcal{A}}\overset{\rightarrow}{\times}T_{\mathcal{B}}:\laxlima{\A}{\alpha}{\B}\rightarrow \laxlima{\A}{\alpha}{\B}$. This functor is invertible with inverse given by $T_{\mathcal{A}}^{-1}\overset{\rightarrow}{\times}T_{\mathcal{B}}^{-1}$. An analogous argument shows that the cotwist functor $F\overset{\rightarrow}{\times}{F'}\dashv G\overset{\rightarrow}{\times}{G'}$ is invertible, showing the desired sphericalness of the adjunction. 
\end{proof}

\subsection{The proof of \Cref{thm:compsph}}\label{subsec:proofBarbacovi}

Let $F\colon \A\leftrightarrow \B\noloc G$ and $F'\colon \A'\leftrightarrow\B\noloc G'$ be two spherical adjunction. We denote $\C=\laxpull{\A}{\B}{\A'}$. Consider further the lax limits $\laxlima{\A}{F}{\B}$ and $\laxlima{\A'}{F'}{\B}$. There are two canonical functors $\on{ev}_{\B},\on{rcof}\colon \laxlima{\A}{F}{\B}\rightarrow \B$, acting on objects via $\on{ev}_{\B}(a\rightarrow b)=b$ and $\on{rcof}(a\rightarrow b)=\cof(F(a)\rightarrow b)$. There are similar functors $\on{ev}_{\B}',\on{rcof}'\colon \laxlima{\A'}{F'}{\B}\rightarrow \B$. Denote by $\D$ the limit of the following diagram in $\St_k$:
\begin{equation}\label{eq:climdiag}
\begin{tikzcd}
            & \laxlima{\A}{F}{\B} \arrow[rd, "\on{ev}_{\B}"] \arrow[ld, "\on{rcof}"] &             & \laxlima{\A'}{F'}{\B} \arrow[ld, "\on{ev}_{\B}'"'] \arrow[rd, "\on{rcof}'"'] &             \\
\mathcal{B} &                                                         & \mathcal{B} &                                                                & \mathcal{B}
\end{tikzcd}
\end{equation}
We further denote by $B_1,B_2\colon \D\rightarrow \B$ the functors contained in limit cone, going to the leftmost and rightmost copy of $\B$ and $B=(B_1,B_2)\colon \D\rightarrow \B^{\times 2}$. The $\infty$-category $\D$ is a concrete model for the $\infty$-category $\R\Gamma^1({\bf G},\F)$ mentioned in \S\ref{subsec:Barbacovi} and the functor $B$ describes the functor $(\on{ev}_p,\on{ev}_q)$.

\begin{lem}\label{lem:B1FF'}
The fiber of $\D\xrightarrow{B_2}\B$ is equivalent to $\C$. The functor $\C\hookrightarrow \D\xrightarrow{B_1}\B$ is equivalent to $\laxpull{F}{\B}{F'}$.
\end{lem}

\begin{proof}
The first part follows from the observation that the fiber of $\on{rcof}':\laxlima{\A'}{F'}{\B}\rightarrow \B$ is equivalent to $\A'$. The second part can be checked for instance by describing all involved functors using the universal properties of the lax limits. 
\end{proof}

\begin{con}
We construct two functor $A,C:\B^{\times 2}\rightarrow \D$, which we show in \Cref{lem:abcadj} to be left and right adjoint to $B=(B_1,B_2):\D\rightarrow \B^{\times 2}$.

Consider the Grothendieck construction $p$ of the diagram \eqref{eq:climdiag} and denote by $\L$ the $\infty$-category of sections of $p$. The $\infty$-category $\D$ can be identified with full subcategory of $\L$ spanned by coCartesian sections.
We denote 
\begin{itemize}
\item by $\E_1$ the full subcategory of $\L$ spanned by $p$-relative left Kan extensions of their restriction to $\laxlima{\A}{F}{\B}$.
\item by $\E_2$ the full subcategory of $\L$ spanned by $p$-relative left Kan extensions of their restriction to the central copy of $\B$.
\item by $\E_3$ the full subcategory of $\L$ spanned by $p$-relative left Kan extensions of their restriction to $\laxlima{\A'}{F'}{\B}$.
\end{itemize}
The restriction functor $\E_1\rightarrow \laxlima{\A}{F}{\B}$ is a trivial fibration by \cite[4.3.2.15]{lurie:htt}.
Choosing a section, we can compose to a functor
$H_1:\B\xrightarrow{\ra{\on{rcof}}}\laxlima{\A}{F}{\B}\rightarrow \E_1\subset \L$. Similarly,
restriction defines a trivial fibration $\E_2\rightarrow \B$ and we obtain a functor
$H_2:\B\rightarrow \E_2\subset \mathcal{L}$ by choosing a section. The functor $H_2$ is left adjoint
to the restriction functor $\on{res}_{\B}:\L\rightarrow \B$ at the central copy of $\B$, see
\cite[4.3.2.17]{lurie:htt}, so that we obtain a counit natural transformation $H_2\circ \on{res}_{\B}\rightarrow \on{id}_{\L}$. Precomposition with $H_1$ yields a natural transformation
\[ \eta_1:H_2\circ \on{ev}_{\B}\circ \ra{\on{rcof}}\simeq H_2\circ \on{res}_{\B}\circ H_1\rightarrow H_1\,.\]
Analogous to the above, we define a functor $H_3:\B\xrightarrow{
\ra{(\on{ev}_{\B}')}}\laxlima{\A'}{F'}{\B}\rightarrow \E_3\subset \L$ and a natural transformation 
\[ \eta_2:H_2\simeq H_2\circ \on{ev}_{\B}'\circ \ra{(\on{ev}_{\B}')}\rightarrow H_3\,.\]
The equivalence above arises from the fact that $\on{ev}_{\B}'$ is a reflective localization.  Finally, we define the functor $C_1:\B\rightarrow \L$ as the colimit of the diagram 
\begin{equation}\label{eq:hcolim}
\begin{tikzcd}
    & H_2\circ \on{ev}_{\B}\circ \ra{\on{rcof}} \arrow[ld, "\eta_1"'] \arrow[rd, "\eta_2\circ \on{ev}_{\B}\circ \ra{\on{rcof}}"] &          \\
H_1 &                                                        & H_3\circ \on{ev}_{\B}\circ \ra{\on{rcof}}
\end{tikzcd}
\end{equation}
in $\Stk(\mathcal{B},\mathcal{L})$. It is straightforward to verify that $C_1$ factors through $\D\subset \L$, and we consider $C_1$ as a functor $\B\rightarrow \D$ in the following.

Exchanging the roles of $F$ and $F'$, or equivalently reflecting diagram \eqref{eq:climdiag} along the vertical axis, and reperforming the above construction, we obtain a functor $C_2:\B\rightarrow \D$. We denote 
\[ C=(C_1,C_2):\B^{\times 2}\rightarrow \D\,.\]   

Performing the above construction of $C$ again, but replacing all right adjoints by left adjoints, we obtain the functor $A=(A_1,A_2):\B^{\times 2}\rightarrow \D$. \end{con}

\begin{lem}\label{lem:abcadj}
The functor $A$ is left adjoint to $B$ and the functor $C$ is right adjoint to $B$.
\end{lem}

\begin{proof}
To determine the left adjoint of $\B\xrightarrow{C_1} \D\subset \L$, we use that passing to
adjoints defines an exact functor $\ra{(\mhyphen)}:\on{Fun}^R(\B^{\times 2},\L)\rightarrow \on{Fun}^L(\L,\B^{\times 2})^\op$ between the functor categories of right or left adjoint functors, and the description of $C_1$ as the colimit of \eqref{eq:hcolim}. The left adjoint of $H_1$ is by \cite[4.3.2.17]{lurie:htt} given by the composite $\L\rightarrow \laxlima{\A}{F}{\B}\xrightarrow{\on{rcof}} \B$ of the restriction functor to $\laxlima{\A}{F}{\B}$ and $\on{rcof}$. Composing this functor with $\D\subset \L$ yields the functor $B_1$. A similar argument shows that the left adjoints of $H_2\circ \on{ev}_{\B}\circ \ra{\on{rcof}}$ and $H_3\circ \on{ev}_{\B}\circ \ra{\on{rcof}}$ both restrict on $\D$, up to equivalence, to the functor
\[ E:\D\xlongrightarrow{\on{res}_\B}\B\xlongrightarrow{\on{rcof}\circ \la{(\on{ev}_{\B}')}}\B\,.\]
Using that restricting along $\D\subset \L$ is also exact, we obtain that the left adjoint of $C_1$ is given by $B_1\simeq B_1\times_{E}E$. 

Similar arguments show $B_2\dashv C_2$, $A_1\dashv B_1$ and $A_2\dashv B_2$. From this it follows that $A\dashv B$ and $B\dashv C$.
\end{proof}

\begin{rem}
The functors $A_1,A_2,C_1,C_2$ correspond to $\la{\on{ev}_p},\la{\on{ev}_q},\ra{\on{ev}_p},\ra{\on{ev}_q}$ from \S\ref{subsec:Barbacovi}. These functors are all fully faithful. This can be checked for instance by an explicit computation of the derived Homs, using their pushout description in \eqref{eq:hcolim} and the facts that $H_1,H_2,H_3$ are fully faithful.
\end{rem}

\begin{prop}\label{prop:abcsph}
The adjunctions $A\dashv B$ and $B\dashv C$ are spherical.
\end{prop}

\begin{proof}
We begin by checking that the twist functor of $A\dashv B$ and the cotwist functor of $B\dashv C$ are invertible. We have the following equivalences
\begin{align*}
\on{id}_{\B}\simeq \on{ev}_{\mathcal{B}}\circ \ra{\on{rcof}} & \quad\quad\quad \on{id}_{\B}\simeq  \on{ev}_{\mathcal{B}}'\circ \ra{(\on{rcof}')}\\
T_{F\dashv G}\simeq \on{rcof}\circ \ra{\on{ev}_{\mathcal{B}}} & \quad\quad\quad T_{F'\dashv G'}\simeq \on{rcof}'\circ \ra{(\on{ev}_{\mathcal{B}}')} 
\end{align*}
where $T_{F\dashv G}$ and $T_{F'\dashv G'}$ denote the cotwist functors of the adjunctions $F\dashv G$ and $F'\dashv G'$, see \cite[Lemma 3.7]{Chr22} for a detailed verification of this. Unraveling the construction of $C$, we thus find 
\[ BC\simeq \on{id}_{\B^{\times 2}}\oplus \begin{pmatrix} 0 & T_{F'\dashv G'}\\ T_{F\dashv G} & 0\end{pmatrix}\]
as an endofunctor of $\B^{\times 2}$. The second summand describes the cotwist functor $T'$ of $B\dashv C$ and is clearly an autoequivalence. A similar description holds for the twist functor of $A\dashv B$. 

To deduce the sphericalness of $A\dashv B$, we apply \cite[Prop.~4.5]{Chr20}, by which it suffices to check that the twist functor $T$ of $A\dashv B$ commutes pointwise with the unit of $A\dashv B$ and that the essential image of $A$ agrees with the essential image of $C$. The former is immediate, since the unit is pointwise the inclusion of a direct summand. Inspecting the constructions, we find equivalences $A_1\circ T_{F\dashv G}(b)\simeq C_2(b)$ and $A_2\circ T_{F'\dashv G'}\simeq C_1(b)$ for all $b\in \B$, which shows the latter. The sphericalness of $B\dashv C$ follows from the sphericalness of $A\dashv B$, see \cite[Cor.~2.6]{Chr20}, concluding the proof.
\end{proof}

\begin{proof}[Proof of \Cref{thm:compsph}]
We show that $\D$ has a $4$-periodic semiorthogonal decomposition with gluing functor
$\laxpull{F}{\B}{F'}$, in the sense of \cite{DKSS:spherical}. We check that there is a semiorthogonal decomposition $(\on{Im}(C_2),\fib(B_2))$ of $\D$, where $\on{Im}(C_2)\subset \D$ denotes the stable subcategory given by the essential image of $C_2$. Given $d\in \fib(B_2)$ and $C_2(b)\in \on{Im}(C_2)$, we have $\on{Map}_{\D}(d,C_2(b))\simeq \on{Map}_{\B}(B_2(d),b)\simeq \ast$. Further, given $d\in \D$, the unit of $B_2\dashv C_2$ defines a morphism $u_d\colon d\rightarrow C_2B_2(d)$ with $B_2(u_d)$ an equivalence, since $C_2$ is fully faithful. This shows $\fib(u_d)\in \fib(B_2)$. 

Similar arguments show the existence of semiorthogonal decompositions $(\fib(B_2),\on{Im}(A_2))$
$(\on{Im}(C_1),\fib(B_1))$, $(\fib(B_1),\on{Im}(A_1))$. Using that
$\on{Im}(C_1)=\on{Im}(A_2)$ and $\on{Im}(C_2)=\on{Im}(A_1)$, we obtain the desired $4$-periodic
semiorthogonal decomposition of $\D$. The gluing functors of the involved semiorthogonal
decompositions are hence spherical, see \cite{DKSS:spherical}.

The left gluing functor of the semiorthogonal decomposition $(\on{Im}(C_2),\fib(B_2))$ is defined as the composite $\fib(B_2)\hookrightarrow \D\xrightarrow{\pi} \on{Im}(C_2)$, where $\pi$ is right adjoint to the inclusion $\on{Im}(C_2)\subset \D$. Under the equivalence $C_2\colon \B\simeq \on{Im}(C_2)$, $\pi$ identifies with $T_{F\dashv G}\circ B_1$. Under the further equivalence $\fib(B_2)\simeq \C$, the gluing functor identifies by \Cref{lem:B1FF'} with $T_{F\dashv G}\circ \laxpull{F}{\B}{F'}$. This shows the sphericalness of $\laxpull{F}{\B}{F'}$. 

It remains to describe the cotwist functor of $\laxpull{F}{\B}{F'}\dashv \ra{\left(\laxpull{F}{\B}{F'}\right)}$. Let $\iota:\C\subset \D$ denote the composite of $\C\simeq \on{fib}(B_2)$ with the inclusion. The equivalence $B_2C_1\simeq T_{F'\dashv G'}$ induces by passing to the adjoint $C_2$ of $B_2$ a natural transformation $\eta:C_1\rightarrow C_2\circ T_{F'\dashv G'}$, satisfying that $\fib(\eta):\B\rightarrow \D$ factors through $\iota:\C\subset \D$. We denote the factorization by $\fib(\eta)':\B\rightarrow \C$. Using that $B_2\circ \iota\simeq 0$, we find 
\begin{align*}
\fib(\eta)'& \simeq \ra{\iota}\circ\iota\circ \fib(\eta)' \simeq \ra{\iota}\circ \fib(\eta)\simeq \ra{\left(\la{(\fib(\eta)}\circ \iota \right)}\\
& \simeq \ra{\left(\cof(T_{F'\dashv G'}^{-1}\circ B_2\rightarrow B_1)\circ \iota \right)} \simeq \ra{\left( B_1\circ \iota\right)}\simeq \ra{\left(\laxpull{F}{\B}{F'}\right)}\,.
\end{align*}
We thus have an equivalence
\[ \left(\laxpull{F}{\B}{F'}\right)\circ \ra{\left(\laxpull{F}{\B}{F'}\right)}\simeq B_1\circ \fib(\eta)\simeq \fib(B_1C_1\rightarrow B_1C_2T_{F'\dashv G'})\simeq \fib(\on{id}_{\B}\rightarrow T_{F\dashv G}\circ T_{F'\dashv G'})\,. \]
One finds the counit of the adjunction $\laxpull{F}{\B}{F'}\dashv \ra{\left(\laxpull{F}{\B}{F'}\right)}$ to arise as the apparent natural transformation $\fib(\on{id}_{\B}\rightarrow T_{F\dashv G}\circ T_{F'\dashv G'})\rightarrow \on{id}_{\B}$, which implies that the cotwist functor of the adjunction is equivalent to $T_{F\dashv G}\circ T_{F'\dashv G'}[-1]$, as desired. This concludes the proof.
\end{proof}

\begin{proof}[Proof sketch of \Cref{thm:schoberbdry}]
Let $i$ be the number of boundary components of ${\bf S}$ and $E^\partial=\coprod_{1\leq j\leq i}E^\partial_j$ the decomposition of the set $E^\partial$ of external edges of ${\bf G}$ into the sets of edges ending on a given boundary component of ${\bf S}$. For each $1\leq j\leq i$, the functor $\prod_{x\in E^\partial_j}\on{ev}_e:\R\Gamma({\bf G},\F)\rightarrow \prod_{x\in E^\partial_j}\F(x)$ is spherical. Though we omit the details, we note that this can be shown by a notationally more involved variant of the proof of \Cref{prop:abcsph} if $|E^\partial_j|\geq 2$ or the proof of \Cref{thm:compsph} if $|E^\partial_j|=1$. 

Let now $1\leq j',j\leq i$ with $j'\neq j$. It is easy to see that $\prod_{x\in E^\partial_j}\on{ev}_e\circ \prod_{x\in E^\partial_{j'}}\ra{\on{ev}_e}\simeq 0$. Hence, the twist functors of the adjunctions $\prod_{x\in E^\partial_j}\on{ev}_e \dashv \prod_{x\in E^\partial_j}\ra{\on{ev}_e}$ and $\prod_{x\in E^\partial_{j'}}\on{ev}_e \dashv \prod_{x\in E^\partial_{j'}}\ra{\on{ev}_e}$ commute. The twist functor of the adjunction $\prod_{x\in E^\partial}\on{ev}_e \dashv \prod_{x\in E^\partial}\ra{\on{ev}_e}$
is equivalent to the composite of the commuting twist functors of the adjunctions $\prod_{x\in E^\partial_l}\on{ev}_e \dashv \prod_{x\in E^\partial_l}\ra{\on{ev}_e}$, $1\leq l\leq i$, and hence an equivalence. The cotwist functor of $\prod_{x\in E^\partial}\on{ev}_e \dashv \prod_{x\in E^\partial}\ra{\on{ev}_e}$ acts on $\F(e)$, with $e\in E^\partial_l$ an external edge, as the cotwist functor of the adjunction $\prod_{x\in E^\partial_l}\on{ev}_e \dashv \prod_{x\in E^\partial_l}\ra{\on{ev}_e}$ and is hence also invertible. This shows the sphericalness of the functor $\prod_{x\in E^\partial}\on{ev}_e$.
\end{proof}

\section{Calabi-Yau complexes}
\label{sec:CY}

The goal of this section is to introduce a notion of Calabi--Yau structure on a categorical complex.
In \S
\ref{subsec:totalHochschild}, we recall the definitions of Hochschild homology and negative cyclic
homology of $k$-linear $\infty$-categories and introduce the total Hochschild homology of a categorical
complex. In \S \ref{subsec:CYstructures}, we define left Calabi--Yau structures on categorical
complexes. In \S \ref{subsec:CYlaxlim}, we discuss some ways in which Calabi--Yau structures arise
on lax limits. In the final \S \ref{subsec:CYtot}, we introduce Calabi--Yau structures on
categorical cubes and show that they induce Calabi--Yau structures on their totalizations. This will
be our main technique to construct examples of categorical Calabi--Yau complexes in \S
\ref{sec:examples}.

\subsection{Total Hochschild and negative cyclic homology}\label{subsec:totalHochschild}

Let $\C$ be a compactly generated $k$-linear $\infty$-category. As explained in \S
\ref{subsec:basic_stable}, the dual of $\C$ with respect to the
monoidal structure on $\Stk$ can be described as
\[
	\C^{\vee} = \Ind(\C_0^{\op}),
\]
where $\C_0 \subset \C$ denotes the subcategory of compact objects. We thus have 
\[
	\St_k(\D,\C) \simeq \Stk(\D\otimes_k\C^\vee,\Modk)
\]
for any $\D \in \Stk$. In particular, we have $\Stk(\C,\C) \simeq
\Stk(\C\otimes_k\C^{\vee},\Modk)$ and can define the diagonal bimodule
$\Delta_{\C}\colon\C\otimes_k\C^{\vee}\rightarrow \Modk$ as the image of the identity functor
$\on{id}_{\C}$ under this equivalence. The $k$-linear Hochschild homology of $\C$ is defined as the trace
\[
	\HH(\C) \coloneqq (\Delta_{\C} \circ \phi(\Delta_{\C}))(k)\in \Modk,
\] 
where $\phi$ denotes the equivalence $\Stk(\C\otimes_k \C^{\vee},\Modk) \simeq
\Stk(\Modk,\C\otimes_k \C^{\vee})$. When $k$ is the sphere spectrum, we
obtain topological Hochschild homology and when $k$ is a field, we obtain the usual Hochschild
complex. Hochschild homology comes equipped with an action of the circle group $S^1$,
see \cite{HSS17}, and we denote its fixed points by $\HH^{S^1}(\C)$. We call $\HH^{S^1}(\C)$ the
negative cyclic homology of $\C$. Hochschild homology and negative cyclic homology form functors
\[
	\HH,\HH^{S^1}\colon\St_k^{\cpt}\rightarrow \Modk
\]
on the subcategory $\St_k^{\cpt} \subset \Stk$ of compactly generated $\infty$-categories and compact objects preserving functors, see
\cite{HSS17}. There is a canonical natural transformation $\HH^{S^1} \to \HH$.

The $k$-linear stable $\infty$-category $\C$ is called smooth if it is compactly generated and
$\Delta_{\C}$ admits a bimodule left dual $\Delta_{\C}^{!}\in \Stk(\Modk,\C\otimes_k \C^{\vee})$. In this case, we sometimes consider $\Delta_{\C}^{!}$ as an object in $\Stk(\C\otimes_k \C^{\vee},\Modk)$ using the
equivalence $\phi$. In this case, the $k$-linear Hochschild homology of $\C$ is equivalent to the $k$-linear morphism object $\on{Mor}_{\Stk(\C\otimes_k \C^{\vee},\Modk)}(\Delta_{\C}^{!},\Delta_{\C})$, see for instance \cite[Lem.~2.27]{Chr23}. Given a $k$-linear,
compact objects preserving functor $F\colon\C\rightarrow \D$, we denote
$\HH(\D,\C)=\cof\HH(F)$ and $\HH^{S^1}(\D,\C)=\cof\HH^{S^1}(F)$. Supposing that $\C$ and $\D$ are smooth, the arising morphism 
\[\on{Mor}_{\Stk(\C\otimes_k \C^{\vee},\Modk)}(\Delta_{\C}^{!},\Delta_{\C})\simeq \HH(\C)\xrightarrow{\HH(F)} \HH(\D)\simeq \on{Mor}_{\Stk(\D\otimes_k \D^{\vee},\Modk)}(\Delta_{\D}^{!},\Delta_{\D})\]
admits a concrete description, see \cite[Prop.~4.4]{BD21} for $k$ a field and \cite[Prop.~2.30]{Chr23} for $k$ arbitrary. It maps a degree $m$-morphism $\alpha\colon \Delta^{!}_{\C}\rightarrow \Delta_{\C}[m]$ to the morphism $\Delta^{!}_{\D}\rightarrow \Delta_{\D}[m]$ obtained as the composite 
\[
\begin{tikzcd}
{\Delta_{\D}^!} \arrow[r, "{\unit}"] & {F_!(\Delta_{\C}^!)} \arrow[d, "F_!(\alpha)"] &           \\
                                    & {F_!(\Delta_{\C})[m]} \arrow[r, "{\counit[m]}"]             & {\Delta_\D[m]}
\end{tikzcd}
\]
where $F_!$ is the functor
\begin{equation}\label{eq:F_*def} \Stk(\C\otimes_k \C^{\vee},\Modk)\simeq \Stk(\C,\C)  \xlongrightarrow{F\circ \mhyphen \circ G} \Stk(\D, \D)\simeq 
\Stk(\D\otimes_k \D^{\vee},\Modk)
\end{equation}
with $G$ the $k$-linear right adjoint of $F$. Above, the natural transformation $\counit$ is the counit of $F\dashv G$ and $\unit$ is the defined as the composite
\[ \Delta_{\D}^!\rightarrow \Delta_{\D}^!\circ \Delta_{\C}\circ \Delta_{\C}^! \rightarrow \Delta_{\D}^!\circ \Delta_{\D}\circ (F\otimes_k \on{Ind}f^{\on{op}})\circ \Delta_{\C}^!\rightarrow (F\otimes_k \on{Ind}f^{\on{op}})\circ \Delta_{\C}^!\simeq F_!(\Delta_{\C}^!)\,,\]
with $\on{Ind}f^{\on{op}}\colon\C^{\vee}\rightarrow \D^{\vee}$ obtained by restricting $F$ to compact objects, passing to opposite categories and $\on{Ind}$-completing.

\begin{defi}
A categorical $n$-complex $\A_*\in \on{Ch}_n(\Stk)$ is called smooth if it consists of smooth
$k$-linear stable $\infty$-categories and all differentials preserve compact objects.
\end{defi}

\begin{defi}
Let $\A_\bullet\in \on{Ch}(\Stk)$ be a bounded, smooth categorical complex. The total negative
cyclic homology $\HH^{S^1,\tot}(\A_\bullet)$ is defined as the total cofiber of the negative cyclic
homology cube, obtained by applying $\HHS$ to the corresponding categorical cube of $\A_\bullet$,
see \S \ref{subsec:complexesofstables} (see \cite{DJW19}[A.2] for terminology and basic results on
total fibers and cofibers). The total Hochschild homology $\HH^{\tot}(\A_\bullet)$ is defined
similarly. 
We have normalized the suspensions in the totalization such that, if $\A_\bullet=\A[0]$ is concentrated in
degree $0$, then $\HH^{S^1,\tot}(\A_\bullet)\simeq \HHS(\A)$ and $\HH^{\tot}(\A_\bullet)\simeq
\HH(\A)$. 
\end{defi}

\begin{defi}
Let $\A_\bullet\in \on{Ch}(\St_k)$.
\begin{itemize}
\item The lower truncation $\tau_{\geq i}(\A_\bullet)$ at $i\in \mathbb{Z}$ of $\A_\bullet$ is defined as the complex 
\[
\begin{tikzcd}
\A_\bullet                & \dots \arrow[r] & \A_{i+1} \arrow[r] \arrow[d] & \A_{i} \arrow[r] \arrow[d] & \A_{i-1} \arrow[d] \arrow[r] & \A_{i-2} \arrow[d] \arrow[r] & \dots \\
\tau_{\geq i}(\A_\bullet) & \dots \arrow[r] & \A_{i+1} \arrow[r]           & \A_i \arrow[r]             & 0 \arrow[r]                           & 0 \arrow[r]                           & \dots
\end{tikzcd}
\]
which is identical to $\A_\bullet$ in degrees $i$ or larger and vanishes in degrees less than $i$. Lower truncation forms a functor $\tau_{\geq i}\colon\on{Ch}(\St_k)\rightarrow \on{Ch}(\St_k)$ and is equipped with a canonical natural transformation $\on{id}_{\on{Ch}(\St_k)}\rightarrow \tau_{\geq i}$.
\item The upper truncation $\tau^{\leq i}(\A_\bullet)$ at $i\in \mathbb{Z}$ of $\A_\bullet$ is defines as the complex 
\[
\begin{tikzcd}
\A_\bullet                & \dots \arrow[r] & \A_{i+2} \arrow[r] \arrow[d] & \A_{i+1} \arrow[r] \arrow[d] \arrow[rd, "\ulcorner", phantom] & \A_{i} \arrow[r] \arrow[d] & \A_{i-1} \arrow[d] \arrow[r] & \dots \\
\tau^{\leq i}(\A_\bullet) & \dots \arrow[r] & 0 \arrow[r]                  & 0 \arrow[r]                                                   & \A_i/\A_{i+1} \arrow[r]    & \A_{i-1} \arrow[r]           & \dots
\end{tikzcd}
\]
which vanishes in degrees larger than $i$, is in degree $i$ given by the cofiber of $\A_{i+1}\xrightarrow{d}\A_i$ and in the other degrees identical to $\A_\bullet$. Upper truncation forms a functor $\tau^{\leq i}\colon\on{Ch}(\St_k)\rightarrow \on{Ch}(\Stk)$ and is equipped with a canonical natural transformation $\on{id}_{\on{Ch}(\Stk)}\rightarrow \tau^{\leq i}$. Note that $\tau^{\leq i}$ preserves smooth categorical complexes, as smoothness of categories is preserved under quotients along compact objects preserving functors. 
\item The two-sided truncation $\tau^{\leq i}_{\geq j}(\A_\bullet)$ is defined as $\tau^{\leq i}\circ \tau_{\geq j}(\A_\bullet)$. Note that the upper and lower truncation commute.
\end{itemize}
\end{defi}

\begin{lem}\label{lem:hhsq}
Let $\A_\bullet\in \on{Ch}(\St_k)$ be a bounded, smooth categorical complex. There is a canonical biCartesian square in $\Modk$:
\[
\begin{tikzcd}
\HH^{S^1,\tot}(\tau_{\geq i}(\A_\bullet)) \arrow[d] \arrow[r] \arrow[rd, "\square", phantom] & \HH^{S^1,\tot}(\tau_{\geq i}^{\leq i+1}(\A_\bullet)) \arrow[d] \arrow[r, "\simeq ", no head] & {\HH^{S^1}(\A_i,\A_{i+1}/\A_{i+2})[i]} \\
\HH^{S^1,\tot}(\tau_{\geq i+1}(\A_\bullet)) \arrow[r]                                        & \HH^{S^1,\tot}(\tau_{\geq i+1}^{\leq i+1}(\A_\bullet)) \arrow[r, "\simeq", no head]         & {\HH^{S^1}(\A_{i+1}/\A_{i+2})[i+1]}            
\end{tikzcd}
\]
\end{lem}

\begin{proof}
This follows from the pasting law for biCartesian squares and the commutative diagram 
\[
\begin{tikzcd}
\HH^{S^1,\tot}(\tau_{\geq i}(\A_\bullet)) \arrow[d] \arrow[r] & {\HH^{S^1}(\A_i,\A_{i+1}/\A_{i+2})[i]} \arrow[d] \arrow[r] \arrow[rd, "\square", phantom] & 0 \arrow[d]                   \\
\HH^{S^1,\tot}(\tau_{\geq i+1}(\A_\bullet)) \arrow[r]         & {\HH^{S^1}(\A_{i+1}/\A_{i+2})[i+1]} \arrow[r]                                                   & {\HH^{S^1}(\A_i)[i+1]}
\end{tikzcd}
\]
in which the right square and the outer square are biCartesian.
\end{proof}

Repeatedly applying \Cref{lem:hhsq}, we obtain that $\HH^{S^1,\tot}(\A_{\bullet})$ is equivalent to the limit of the diagram
\[
\begin{tikzcd}[column sep=small]
                                                                                          &                                                                                     & {\HH^{S^1}(\A_{i-1},\A_{i}/\A_{i+1})[i-1]} \arrow[d] \arrow[r] & \dots \\
                                                                                          & {\HH^{S^1}(\A_i,\A_{i+1}/\A_{i+2})[i]} \arrow[d] \arrow[r] & {\HH^{S^1}(\A_{i}/\A_{i+1})[i]}                                         &       \\
{\HH^{S^1}(\A_{i+1},\A_{i+2}/\A_{i+3})[i+1]} \arrow[r] \arrow[d] & {\HH^{S^1}(\A_{i+1}/\A_{i+2})[i+1]}                                 &                                                                                         &       \\
\dots                                                                                     &                                                                                     &                                                                                         &      
\end{tikzcd}
\]
The limit of the above diagram is equivalent to an equalizer as follows:

\begin{prop}\label{prop:hhtot}
Let $\A_\bullet\in \on{Ch}(\St_k)$ be a bounded, smooth categorical complex.  Then $\HH^{S^1,\tot}(\A_\bullet)$ is equivalent to the equalizer of the diagram in $\Modk$
\[
\begin{tikzcd}
{\bigoplus_{i\in \mathbb{Z}}\HH^{S^1}(\A_i,\A_{i+1}/\A_{i+2})[i]} \arrow[r, "\bigoplus \delta_i", shift left] \arrow[r, "\bigoplus \pi_i"', shift right] & {\bigoplus_{i\in \mathbb{Z}}\HH^{S^1}(\A_{i}/\A_{i+1})[i]\,,}
\end{tikzcd}
\]
where \[\delta_i\colon\HH^{S^1}(\A_i,\A_{i+1}/\A_{i+2})[i]\rightarrow \HH^{S^1}(\A_{i+1}/\A_{i+2})[i+1]\] is the fiber map and 
\[ \pi_i\colon\HH^{S^1}(\A_i,\A_{i+1}/\A_{i+2})[i]\rightarrow \HH^{S^1}(\A_{i}/\A_{i+1})[i]\] 
is the map induced by the quotient functor $\A_{i}\rightarrow \A_{i}/\A_{i+1}$.\\
The above statement also holds when replacing $\HH^{S^1}$ with $\HH$ and $\HH^{S^1,\tot}$ with $\HH^{\tot}$.
\end{prop}

\begin{rem}
\Cref{prop:hhtot} expresses, that a Hochschild or total negative cyclic homology class consists of a compatible
family of relative Hochschild or relative negative cyclic homology classes of the functors $\A_{i+1}/\A_{i+2}\rightarrow
\A_i$, shifted into the appropriate degrees.
\end{rem}

\subsection{Left Calabi--Yau structures}\label{subsec:CYstructures}

Let $F\colon\C\rightarrow \D$ be a compact objects preserving, $k$-linear functor between smooth $k$-linear $\infty$-categories. As explained in \S\ref{subsec:totalHochschild}, a $k$-linear Hochschild class $\sigma\colon k[n]\to \HH(\D,\C)$ defines a diagram
\[
\begin{tikzcd}
{\Delta_{\D}^!} \arrow[r] & {F_!(\Delta_{\C}^!)} \arrow[d]               &                       \\
                                         & {F_!(\Delta_{\C})[n-1]}\arrow[r] & {\Delta_{\C}[n-1]}
\end{tikzcd}
\]
together with a choice of null-homotopy. It hence gives rise to a diagram with horizontal fiber and cofiber sequences as follows:
\[
\begin{tikzcd}
{\Delta_{\D}^!} \arrow[r] \arrow[d] & {F_!(\Delta_{\C}^!)} \arrow[d] \arrow[r] & \cof \arrow[d] \\
\fib \arrow[r]                       & {F_!(\Delta_{\C})[n-1]} \arrow[r]                    & {\Delta_{\C}[n-1]}       
\end{tikzcd}
\]
We call the Hochschild class $\sigma$ {\em non-degenerate} if all vertical maps in the above diagram are equivalences. 

\begin{defi}
\begin{enumerate}
\item A weak left $n$-CY structure on the functor $F$ consists of a non-degenerate Hochschild class $\sigma\colon k[n]\to \HH(\D,\C)$. 
\item A left $n$-CY structure on the functor $F$ consists of a negative cyclic class $\eta\colon k[n]\to
	\HH^{S^1}(\D,\C)$, whose image under $\HH^{S^1}(\D,\C)\to \HH(\D,\C)$ defines a
	non-degenerate Hochschild class. 
\end{enumerate}
\end{defi}

We refer to \Cref{ex:A2CY}, \Cref{thm:gorcy} and \Cref{thm:mfd} for examples of left Calabi--Yau
structures. 

\begin{defi}\label{def:CYcplx}
Let $\A_\bullet\in \Ch(\Stk)$ be a bounded, smooth categorical complex. 
\begin{enumerate}
\item A weak left $n$-Calabi--Yau structure on $\A_\bullet$ consists of a class $\sigma\colon k[n]\rightarrow \HH^{\tot}(\A_\bullet)$, whose composite with the morphism from \Cref{prop:hhtot} 
\[ \HH^{\tot}(\A_\bullet) \longrightarrow \bigoplus_{i\in \mathbb{Z}}\HH(\A_{i},\A_{i+1}/\A_{i+2})[i]\] 
defines a collection of weak left $(n-i)$-Calabi--Yau structures on the functors\newline
$\A_{i+1}/\A_{i+2}\rightarrow \A_i$.
\item A left $n$-Calabi--Yau structure on $\A_\bullet$ consists of a class $\eta\colon k[n]\rightarrow \HH^{S^1,\tot}(\A_\bullet)$, whose composite with the morphism from \Cref{prop:hhtot} 
\[ \HH^{S^1,\tot}(\A_\bullet) \longrightarrow \bigoplus_{i\in \mathbb{Z}}\HHS(\A_{i},\A_{i+1}/\A_{i+2})[i]\] 
defines a collection of left $(n-i)$-Calabi--Yau structures on the functors $\A_{i+1}/\A_{i+2}\rightarrow \A_i$.
\end{enumerate}
\end{defi}

\begin{rem}
There is an apparent analogue of \Cref{def:CYcplx} for bounded complexes of proper $k$-linear $\infty$-categories, obtained by replacing left Calabi--Yau structures with right Calabi--Yau structures.
\end{rem}

\subsection{Calabi--Yau structures and lax limits}\label{subsec:CYlaxlim}

Consider a colimit diagram in $\Stk$
\[
\begin{tikzcd}
                      &                                                                   & \A'' \arrow[d] \\
                      & \A \arrow[d] \arrow[r] \arrow[rd, "\ulcorner", phantom] & \B' \arrow[d]  \\
\A' \arrow[r] & \B \arrow[r]                                             & \C           
\end{tikzcd}
\]
where all appearing $\infty$-categories are smooth and all functors preserve compact objects. Suppose further that the functors $\A\times \A'\rightarrow \B$ and $\A\times \A''\rightarrow \B'$ carry left
$n$-Calabi--Yau structures, which are compatible at $\A$. 

\begin{thm}\label{thm:glueCY}
The functor $\A'\times \A''\rightarrow \C$ inherits a left $n$-Calabi--Yau structure. 
\end{thm}
\begin{proof}
See \cite[Thm.~6.2]{BD19} if $k$ is a field and \cite[Thm.~3.14]{Chr23} for $k$ an arbitrary $\mathbb{E}_\infty$-ring spectrum, see . 
\end{proof}

\Cref{thm:glueCY} admits a refinement, in the form of a relative Calabi--Yau structure on the
\dpo{}. Before stating that result in \Cref{lem:laxcy}, we need to consider the following class of
examples of relative left Calabi--Yau structures.

\begin{exa}\label{ex:A2CY}
Let $\A$ be a smooth $k$-linear $\infty$-category with a left $(n-1)$-Calabi--Yau structure $\eta\colon k[n-1]\to \HH^{S^1}(\A)$. Consider the functor 
\[ 
	G_{A_2}=(G_1,G_2,G_3)\coloneqq (\on{ev}_{0},\cof,\on{ev}_1)\colon \Fun(\Delta^1,\A)\longrightarrow \A^{\times 3}\,,
\]
with components given by the evaluations at $i \in \{0,1\}$ and the cofiber, respectively. The
notation $A_2$ comes from regarding $\Fun(\Delta^1,\A)$ as the $\infty$-category of $\A$-valued representations of the $A_2$-quiver. The left adjoint
\[ 
	F_{A_2}=(F_1,F_2,F_3)\colon \A^{\times 3}\rightarrow \Fun(\Delta^1,\A)
\] 
of $G_{A_2}$ admits a left $n$-Calabi--Yau structure which restricts to the given left
$(n-1)$-Calabi--Yau structure $\eta^{\oplus 3}:k^{\oplus 3}\to \HH^{S^1}(\A)^{\oplus 3}\simeq
\HH^{S^{1}}(\A^{\times 3})$ on $\A^{\times 3}$. 

This can be seen as follows. For $\A=\Modk$ with the apparent left
$0$-Calabi--Yau structure, this left $1$-Calabi--Yau structure on $F_{A_2}$ is described in \cite[Thm.~5.14]{BD19} and \cite[Lemma.~5.6]{Chr23}. For a general $k$-linear $\infty$-category $\A$, we have $\Fun(\Delta^1,\Modk) \otimes
\A \simeq \Fun(\Delta^1,\A)$ and $\Modk^{\times 3}\otimes \A\simeq \A^{\times 3}$ using the
symmetric monoidal structure on $\Stk$. The desired left Calabi--Yau structure now arises as the
image of the above pair of Calabi--Yau structures under the canonical map 
\[
	\HH_1^{S^1}(\Fun(\Delta^1,\Modk),\Modk^{\times 3})\times
	\HH_{n-1}^{S^1}(\A)\rightarrow \HH_n^{S^1}(\Fun(\Delta^1,\A),\A^{\times 3})\,.
\]
\end{exa}

\begin{lem}\label{lem:laxcy}
The functor 
\begin{equation}\label{eq:laxcy}
\A\times \A'\times \A''\rightarrow \laxpush{\B}{\A}{\B'}
\end{equation} 
inherits a left $n$-Calabi--Yau structure which induces the left $n$-Calabi--Yau structure on
${\displaystyle \A'\times \A''}\rightarrow \C$ from Theorem \ref{thm:glueCY} when forming the
pushout along the functor $\A\rightarrow 0$.
\end{lem}

\begin{proof}
We have the following colimit diagram in $\St_k$:
\[
\begin{tikzcd}[column sep=large]
                                          &                                                                                      & \A \arrow[d, "F_2"]                      \\
                                          & {\displaystyle \A\times \A} \arrow[r, "{(F_1,F_3)}"] \arrow[d] \arrow[rd, "\ulcorner", phantom, near end] & {\Fun(\Delta^1,\A)} \arrow[d] \\
{\displaystyle \A'\times \A''} \arrow[r] & {\displaystyle \B\times \B'} \arrow[r]                                             & {\displaystyle \laxpush{\B}{\A}{\B'}}             
\end{tikzcd}
\]
Applying \Cref{thm:glueCY}, we find that the left $n$-Calabi--Yau structures on $F_{A_2}$ and on $\A\times \A\times A'\times A''\to \B\times \B'$ glue to a left $n$-Calabi--Yau structure on the functor \eqref{eq:laxcy}.
\end{proof}

We end this section with the following further observations concerning left Calabi--Yau structures on lax limits.

\begin{lem}\label{lem:laxcy2}
Let $F\colon\A\rightarrow \B$ be a $k$-linear, compact objects preserving functor between smooth $k$-linear stable $\infty$-categories.
\begin{enumerate}[(1)]
\item Let $\pi\colon\A'\rightarrow \A$ be a $k$-linear functor which preserves compact object and admits a fully faithful right adjoint. A class 
\[ \eta\colon k[n]\lra \HHS(\B,\A')\]
determines a left $n$-Calabi--Yau structure on $F\circ \pi$ if and only if its composite with 
\[
\HHS(\B,\A')\lra \HHS(\B,\A)
\]
determines a left $n$-Calabi--Yau structure on $F$.
\item Assume that there are $k$-linear functors $F_i\colon\A\rightarrow \B_i$ with $1\leq i\leq n$, such that 
\[ \B\simeq \laxpush{(\laxpush{(\laxpush{\B_1}{\A}{\B_2})}{\A}{\dots})}{\A}{\B_n}\]
is an iterated directed pushout and $F$ is equivalent to the iteratively induced functor into the directed pushout. Compatible left $n$-Calabi--Yau structures on the functors $F_i$ for all $1\leq i\leq n$ canonically determine a left $n$-Calabi--Yau structure on $F$.
\item Conversely, in the setting of (2), a left $n$-Calabi--Yau structure on $F$ determines compatible left $n$-Calabi--Yau structures on the functors $F_i$, $1\leq i\leq n$; these assignments are inverse to each other.
\end{enumerate}
\end{lem}

\begin{proof}
Part (1) follows from the observation that the functor $\pi_*$, see \eqref{eq:F_*def}, maps the diagonal bimodule $\Delta_{\A'}$ to $\Delta_{\A}$, by the fully faithfulness of $\ra{\pi}$

Part (2) follows from repeated application of \Cref{lem:laxcy} (with $\A'=\A''=0$).

For part (3), we note that the projection map $\pi_i\colon\B\twoheadrightarrow \B_i$ admits a fully faithful right adjoint for each $1\leq i\leq n$. Again $(\pi_i)_*$ maps $\Delta_{\B}$ to $\Delta_{\B_i}$ and the arising map 
\[ \HHS(\B,\A)\lra \HHS(\B_i,\A)\]
thus left $n$-Calabi--Yau structures to left $n$-Calabi--Yau structures.
\end{proof}

\subsection{Cubical Calabi--Yau structures}\label{subsec:CYtot}\label{subsec:CYcube}

Recall that $\I=[1]^\op=\{1\to 0\}$. Let 
\[
	\A_*\colon \I^N \lra \Stk
\]
be a smooth, cubical diagram of $k$-linear $\infty$-categories. Applying $\HHS$ yields a cubical diagram 
\[
	\HHS(\A_*)\colon \I^N \lra \Modk.
\]
We denote by $\HH^{S^1,\tot}(\A_*)\in \Modk$ the total cofiber of $\HHS(\A_*)$. To clarify the
chosen grading, note that if $\A_\ast$ assigns $0$ to all vertices of the cube except $(0,\dots,0)$,
then we have $\HH^{S^1,\tot}(\A_*)\simeq \HH^{S^1}(\A_{(0,\dots,0)})$. For $L \in \Modk$, we refer
to a map $\eta\colon L \to \HH^{S^1,\tot}(\A_*)$ as an $L$-class. Using that the total cofiber is
the $N$-fold suspension of the total fiber, we may identify an $L$-class $\eta$ with a
natural transformation in $\Fun(\I^{N},\Mod_k)$
\[
	\eta\colon L[-N]_0 \lra \HHS(\A_*)
\]
where, for $E \in \Modk$, we denote by $E_0$ the cubical diagram
\[
	E_0\colon \I^N \lra \Modk, \quad x \mapsto \begin{cases} 
						E & \text{for $x = (1,1,...,1)$},\\ 0 &
						\text{else.} 
					      \end{cases}
\]
From this perspective, we obtain from $\eta$ the following data:
\begin{enumerate}[label=\arabic*.]
	\item for every $1 \le i \le N$, a morphisms in $\Fun(\I^{N-1},\Mod_k)$
		\[
			\partial_i \eta\colon L[-N]_0 \lra \HHS(\A_{\partial_i *})
		\]
		where $\A_{\partial_i *}$ denotes the restriction of the cube $\A_*$ to the face
		$\I^{N-1} \to \I^N$ obtained by setting the $i$th coordinate to $1$, and hence a
		class
		\begin{equation}\label{eq:faceclass}
			L[-1] \lra \HHStot(\A_{\partial_i *}),
		\end{equation}
	\item a diagram
		\[
			\begin{tikzcd} 
				L[-1]  \ar{rr}\ar{d} & & 0\ar{d}\\
				\colim\nolimits_{<(0,\dots,0)} \HHS(\A_{\ast}) \ar{r} & \HHS(\colim\nolimits_{<(0,\dots,0)} \A_{*})\ar{r} &
				\HHS(\A_{(0,\ldots,0)}),
			\end{tikzcd}
		\]
		and hence a relative class
		\begin{equation}\label{eq:matchingclass}
			L \lra \HHStot(\A_{\ast}) \simeq \HHS(\A_{(0,\dots,0)},\colim\nolimits_{<(0,\dots,0)}\A_\ast),
		\end{equation}
		where above $\colim\nolimits_{<(0,\dots,0)}$ denotes the colimit over the punctured cube $\I^N_{<(0,\dots,0)}=\I^N \setminus \{(0,\ldots,0)\}$. 
\end{enumerate}

We recursively define an $n$-Calabi--Yau structure on cubical diagrams as follows.

\begin{defi}
	Let $\A_*\colon \I^N \lra \Stk$ be a smooth cubical diagram in $\Stk$. A class
	$\eta\colon k[n] \to \HH^{S^1,\tot}(\A_*)$ is called an $n$-Calabi--Yau structure on $\A_*$, if
	\begin{enumerate}
		\item for every $1 \le i \le N$, the class 
			\[
				k[n-1] \to \HH^{S^1,\tot}(\A_{\partial_i *})
			\]
			from \eqref{eq:faceclass} defines a left $(n-1)$-Calabi--Yau structure on the $(N-1)$-cube $\A_{\partial_i *}$,
		\item the class $k[n]\to \HH^{S^1}(\A_{(0,\dots,0)},\colim\nolimits_{<(0,\dots,0)}\A_\ast)$
			from \eqref{eq:matchingclass} defines a left $n$-Calabi--Yau structure	on the functor 
			\begin{equation}\label{eq:CYcolimC}
				\colim\nolimits_{<(0,\dots,0)} \A_{*} \lra \A_{(0,\ldots,0)}.
			\end{equation}
	\end{enumerate}
\end{defi}

\begin{exa}
	Let $\A_\bullet \in \Ch(\Stk)$ be categorical complex concentrated in degrees $N$, $N-1$, ...,
	$0$ so that $\A_\bullet$ gives rise to a categorical $N$-cube $\A_\ast$ (see
	\eqref{eq:coherent_cube} in \S \ref{subsec:complexesofstables}). Then $\HHStot(\A_\bullet) \simeq
	\HH^{S^1,\tot}(A_\ast)$ and a class $\eta\colon k[n]\to \HH^{S^1,\tot}(\A_\bullet)$ describes a
	left $n$-Calabi-Yau structure on $\A_\bullet$ if and only if it describes a left
	$n$-Calabi-Yau structure on $\A_\ast$.
\end{exa}

\begin{thm}\label{thm:CYcubetot} 
Let $\A_*$ be a smooth cubical diagram in $\Stk$ equipped with a left $n$-Calabi-Yau structure. Then
its coproduct totalization inherits a canonical left $n$-Calabi--Yau structure. 
\end{thm}

Before proving \Cref{thm:CYcubetot}, we collect some preparatory results. We begin with an analysis
of the differentials in the totalization of a categorical cube. 
For notational convenience, we will use the reparameterization of the poset
$\I^{N+1}=([1]^\op)^{N+1}$ by the poset $\Pc([N])^\op$ of subsets of the set $[N]=\{0,1,...,N\}$ (by
means of associating to a subset its characteristic function). For instance, $(1,\dots,1)$ is
identified with $[N]\in \Pc([N])^\op$, $(1,0,\dots,0)$ with $\{1\}\in  \Pc([N])^\op$ and $(0,\dots,0)$ with $\emptyset\in \Pc([N])^\op$.

Let $\A_*:\Pc([N])^\op\rightarrow \St_k$ be a smooth categorical $N$-cube with $N\geq 2$. When
taking product and coproduct totalizations, we equip the subset of $\Pc([N])^\op$ of subset of
cardinality $i$ with the canonical total order, where $I<J$ if $\on{min}(I\backslash
J)<\on{min}(J\backslash I)$, for $I,J\subset [N]$. 

\begin{lem}\label{lem:totcube=ipb}
There exist equivalences in $\Stk$
\[ 
	\tot^{\times}(\A_*)_{i}\simeq \bigtimes^{\curvearrowright}_{\tot^{\times}(\A_*)_{i-1}}\A_J
\]
and 
\[ 
	\tot^{\amalg}(\A_*)_{i}\simeq \coprod^{\curvearrowright}_{\tot^{\amalg}(\A_*)_{i+1}}\A_J
\]
where the iterated \dpb{} and the iterated \dpo{} run over all $J\subset [N]$ with $|J|=i$, with the order as specified above.
\end{lem}

\begin{proof}
We only verify the formula for the product totalization, the coproduct totalization can be treated
analogously.
The proof is by induction on $N$ using the iterative definition of the product totalization. The
case $N=2$ is clear. We proceed with the induction step. Given a categorical $(N+1)$-cube $\A_{*}$,
we can consider it as describing a chain map $f_*:\B_*\rightarrow \C_*$ between $N$-cubes. The
induction step now follows from combining the following two observations.
\begin{itemize}
\item We have that $\tot^{\times}(\B_*)_i\simeq
	\bigtimes^{\curvearrowright}_{\tot^{\times}(\A_*)_{i-1}}\B_J$ and a similar statement for
	$\C_*$. This follows from the fact that the functor $\tot^{\times}(\A_*)_{i-1}\rightarrow
	\tot^{\times}(\B_*)_{i-1}$ is a reflective localization.
\item Spelling out the definition of the iterated \dpb{}, one finds that the order of the bracketing
	is irrelevant (i.e.~satisfies associativity), implying that 
	\[ 
		\left(
	\bigtimes^{\curvearrowright}_{\tot^{\times}(\A_*)_{i-1}}\B_J \right)
\times^{\curvearrowright}_{\tot^{\times}(\A_*)_{i-1}} \left(
\bigtimes^{\curvearrowright}_{\tot^{\times}(\A_*)_{i-1}}\C_J  \right) \simeq
\bigtimes^{\curvearrowright}_{\tot^{\times}(\A_*)_{i-1}}\A_J\,.
\]
\end{itemize}
\end{proof}

Let $J\subset [N]$ with $|J|=i$. Our next goal is to find a description of the functor
\[ \tot^{\amalg}(\A_*)_{i+2}\xrightarrow{d} \tot^{\amalg}(\A_*)_{i+1}\simeq \coprod^{\curvearrowright}_{\tot^{\amalg}(\A_*)_{i+1}}\A_{J'}\twoheadrightarrow \A_{J}\,.\]

\begin{defi}\label{con:kercokerJ}
Fix $J\subset [N]$ with $|J|=i$. 
\begin{enumerate}[(1)]
\item We define the stable $\infty$-category $P^J$ as the full subcategory  
\[ \P^J\coloneqq \bigtimes^{\curvearrowright}_{\tot^{\times}(\A_*)_{i-2},\,j\in J}\A_{J\backslash\{j\}}\subset \bigtimes^{\curvearrowright}_{\tot^{\times}(\A_*)_{i-2}}\A_{J'}\,,\] meaning the stable subcategory generated by the components $\A_{J\backslash \{j\}}$ of the iterated lax product. We further define $\on{ker}(d)^J \subset P^J$ as the fiber in $\Stk$ of 
\[ P^J\subset \tot^{\times}(\A_*)_{i-1} \xlongrightarrow{d} \tot^{\times}(\A_*)_{i-2}\,.\]
\item The dual version $P_J$ is defined as the full subcategory
\[ \P_J\coloneqq \coprod^{\curvearrowright}_{\tot^{\amalg}(\A_*)_{i+2},\,j\in [N]\backslash J}\A_{J\cup \{j\}}\subset \coprod^{\curvearrowright}_{\tot^{\amalg}(\A_*)_{i+2}}\A_{J'}\simeq \tot^{\amalg}(\A_*)_{i+1} \,.\]
The left adjoint of this inclusion is denoted by $\tot^{\amalg}(\A_*)_{i+1}\twoheadrightarrow \P_J$. We further define $\on{coker}(d)_J$ as the cofiber in $\Stk$ of
\[ \tot^{\amalg}(\A_*)_{i+2}\xlongrightarrow{d} \tot^{\amalg}(\A_*)_{i+1} \twoheadrightarrow P_J\,.\]
\end{enumerate}
\end{defi}

\begin{lem}\label{lem:coker(d)colim}
Let $J\subset [N]$.
\begin{enumerate}
\item[(1)] The $k$-linear stable $\infty$-category $\on{ker}(d)^J$ is equivalent to the limit in $\Stk$ of the restriction of $\A_*$ to the full subcategory $\Pc([N])^\op_{>J}$ of $\Pc([N])^\op$ spanned by objects $J'\neq J$ satisfying $J'\subset J$.
\item[(2)] The $k$-linear $\infty$-category $\on{coker}(d)_J$ is equivalent to the colimit in $\Stk$ of the restriction of $\A_*$ to the full subcategory $\Pc([N])^\op_{<J}$ of $\Pc([N])^\op$ spanned by objects $J'\neq J$ satisfying $J\subset J'$.
\end{enumerate}
\end{lem}

\begin{proof}
Part (2) is dual to part (1), i.e.~follows from passing to right adjoints. Part (1) follows from an induction over $|J|$ and by decomposing the limit over $\Pc([N])^\op_{>J}$ via \cite[4.2.3.10]{lurie:htt}. 
\end{proof}

Via the universal properties of these limits and colimits, we obtain functors $d^J\colon \A_J\rightarrow \on{ker}(d)^J$ and $d_J\colon \on{coker}(d)_J\rightarrow \A_J$.

\begin{lem}\label{lem:factordircoprod}
Let $\A_*$ be a categorical $n$-cube and let $J\subset [N]$ with $|J|=i$.
\begin{enumerate}[(1)]
\item There is a commutative diagram in $\Stk$ 
\[
\begin{tikzcd}
\A_J \arrow[d, hook] \arrow[r, "d^J"]                 & \on{ker}(d)^J \arrow[rd, "0"] \arrow[d, hook]      &                                     \\
\tot^{\times}(\A_*)_{i} \arrow[r, "d"] & \tot^{\times}(\A_*)_{i-1} \arrow[r, "d"] & \tot^{\times}(\A_*)_{i-2}
\end{tikzcd}\]
with fully faithful vertical functors. Furthermore, if $\A_*$ takes values in limit preserving functors, then all functors in the above diagram also preserve limits.
\item There is a commutative diagram in $\Stk$
\[
\begin{tikzcd}
\tot^{\amalg}(\A_*)_{i+2} \arrow[r, "d"] \arrow[rd, "0"'] & \tot^{\amalg}(\A_*)_{{i+1}} \arrow[r, "d"] \arrow[d, two heads] & \tot^{\amalg}(\A_*)_{i} \arrow[d, two heads] \\
                                                                    & \on{coker}(d)_J \arrow[r, "d_J"]                                               & \A_J                                                
\end{tikzcd}\]
such that the vertical functors admit fully faithful right adjoints. Furthermore, if $\A_*$ takes values in compact objects preserving functors, then all functors in the above diagram also preserve compact objects.
\end{enumerate}
\end{lem}

\begin{proof}
We only prove part (1), part (2) is dual. Using the notation of \Cref{con:kercokerJ}, we have a commutative diagram in $\Stk$:
\[
\begin{tikzcd}
                                                            & \on{ker}(d)^J \arrow[d] \arrow[r] \arrow[rdd, "\lrcorner", phantom, near start] & 0 \arrow[dd]                        \\
\A_J \arrow[d, hook] \arrow[r] \arrow[r] \arrow[ru, dotted] & P^J \arrow[rd] \arrow[d, hook]                                      &                                     \\
\tot^{\times}(\A_*)_{i} \arrow[r, "d"]              & \tot^{\times}(\A_*)_{i-1} \arrow[r, "d"]                  & \tot^{\times}(\A_*)_{i-2}
\end{tikzcd}
\]
The top right square is pullback, the dotted arrow arises via the universal property of this pullback. Unraveling the definition of $d$, one can further show that the dotted arrow is described by the functor $d^J$. If $\A_*$ takes values in limit preserving functors, the pullback is also equivalent to the pushout of the $k$-linear left adjoint diagram in $\Stk$. Passing again to right adjoints, we find that all functors in the diagram preserve limits.
\end{proof}

\begin{exa}
Consider a categorical $3$-cube $\A_*$. Then 
\[ \on{coker}(d)_{\{j\}}\simeq \A_{\{j,j'\}}\amalg_{\A_{[3]}} \A_{\{j,j''\}}\]
with $j,j',j''\in [3]$ pairwise distinct. The $\infty$-category $\on{coker}(d)_{\emptyset}$ is equivalent to the colimit in $\Stk$ of the diagram
\[\begin{tikzcd}
{\A_{[3]}} \arrow[rr] \arrow[dd] \arrow[rd] &                                      & {\A_{\{1,2\}}} \arrow[dd] \arrow[rd] &            \\
                                            & {\A_{\{1,3\}}} \arrow[dd] \arrow[rr] &                                      & \A_{\{1\}} \\
{\A_{\{2,3\}}} \arrow[rr] \arrow[rd]        &                                      & \A_{\{2\}}                           &            \\
                                            & \A_{\{3\}}                           &                                      &           
\end{tikzcd}
\]
which is also equivalent to the iterated (usual, not directed) pushout
\[
\left(\A_{\{1\}}\amalg_{\A_{\{1,2\}}} \A_{\{2\}}\right)\amalg_{\A_{\{1,3\}}\amalg_{\A_{[3]}} \A_{\{2,3\}}}\A_{\{3\}}\,.
\]
\end{exa}

\begin{prop}\label{prop:CYtot}
Let $\A_\bullet\in \Ch_N(\Stk)$ be a smooth cube. A class 
\[ \eta\colon k[n]\to \HHStot(\tot^{\amalg}(\A_\bullet))\,,\] 
which induces left $(n-|J|)$-Calabi--Yau structures on the functors $d_J:\on{coker}(d)_J\rightarrow \A_J$, for all $J\subset [N]$, describes a left $n$-Calabi--Yau structure on $\tot^{\amalg}(A_\bullet)$.
\end{prop}

\begin{proof}
Combine parts (1) and (2) of \Cref{lem:laxcy2} with \Cref{lem:totcube=ipb} and \Cref{lem:factordircoprod}.
\end{proof}

\begin{proof}[Proof of \Cref{thm:CYcubetot}]
There is a canonical equivalence $\HHStot(\A_\ast)\simeq \HHStot(\tot^{\amalg}(\A_\ast))$, which can be obtained using the recursive definition of the coproduct totalization. Consider a left $n$-Calabi-Yau structure $\eta\colon k[n]\to \HHStot(\A_\ast)$ of $\A_\ast$. Under the above equivalence, the arising top non-degeneracy condition of $\HHStot(\A_{\partial_i\dots \partial_j \ast})$ implies that the functor $d_{J}:\on{coker}(d)_{J}\rightarrow \A_{J}$ with $J=\{i,\dots,j\}$ inherits a left $(n-|J|)$-Calabi--Yau structure from $\eta$. The Theorem thus follows from \Cref{prop:CYtot}.
\end{proof}

\begin{rem}
In view of part (3) of \Cref{lem:laxcy2}, the proof of \Cref{thm:CYcubetot} also shows the 'converse' of \Cref{thm:CYcubetot}, namely that any left $n$-Calabi--Yau structure of $\on{tot}^\amalg(\A_\ast)$ canonically determines a left $n$-Calabi-Yau structure on $\A_\ast$.
\end{rem}

\begin{rem}
Then tensor product of a left $n$-Calabi--Yau functor with a left $m$-Calabi--Yau $\infty$-category yields an $(m+n)$-Calabi--Yau functor. We expect that this generalizes in the sense that the tensor product of two or more Calabi--Yau functors yields a Calabi--Yau square or cube. 
\end{rem}

\section{Examples}
\label{sec:examples}

\subsection{Normal crossings divisors and categorical intersection complexes}
\label{subsec:completeintersection}

As a standing assumption in this section, we require that all schemes are separated, reduced and of
finite type over a field $k$ of characteristic zero. Morphisms between schemes are assumed to be
separated. The category of such schemes is denoted $\Schk$.

Let $k$ be a field of characteristic $0$, $Z$ a smooth scheme over $k$ and $i\colon D \subset Z$ a
normal crossing divisor, by which we mean a union $D=\bigcup_{1\leq i\leq N} D_i$ of smooth divisors
$D_i$, intersecting transversely. We may organize the intersections
\[
	D_I \coloneqq \bigcap_{i \in I} D_i,
\]
for the various subsets $I \subset \{1,...,N\}$, along with their inclusions $D_I \subseteq D_J$ for
$I \supseteq J$, into a cubical diagram $D_*$ of smooth schemes. Here, we interpret the empty
intersection $D_{\emptyset}$ as the ambient scheme $Z$.
For example, for $N=2,3$, the resulting cube can be depicted as follows:
\[
\begin{tikzcd}
D_1\cap D_2 \arrow[r] \arrow[d] & D_1 \arrow[d] \\
D_2 \arrow[r]                   & Z 
\end{tikzcd} \hspace{3em}
\begin{tikzcd}[column sep=small]
D_1\cap D_2\cap D_3 \arrow[rr] \arrow[rd] \arrow[dd] &                                   & D_1\cap D_3 \arrow[rd] \arrow[dd] &                     \\
                                                     & D_2\cap D_3 \arrow[dd] \arrow[rr] &                                   & D_3 \arrow[dd]      \\
D_1\cap D_2 \arrow[rd] \arrow[rr]                    &                                   & D_1 \arrow[rd]                    &                     \\
                                                     & D_2 \arrow[rr]                    &                                   & Z
\end{tikzcd}
\]
Consider the functor $\IndCoh \colon N(\Schk)\rightarrow \St_k$, which assigns to a
scheme $Y$ its $k$-linear $\infty$-category $\IndCoh(Y)$ of $\on{Ind}$-coherent sheaves and to a
morphism $f\colon Y\rightarrow Y'$ of schemes the $k$-linear functor $f_*\colon
\IndCoh(Y)\rightarrow \IndCoh(Y')$, see \cite{GR17} for details. If $f$ is proper, then the functor $f_*$ admits a right
adjoint $f^!$, which preserves colimits, thus defining a morphism in $\Stk$. Applying $\IndCoh$ to
the cube $D_*$, we obtain a categorical $N$-cube $\IndCoh(D_\ast)$, for example in
the case $N=2$ given by:
\begin{equation}\label{eq:indcohsq}
	\begin{tikzcd}
	\IndCoh(D_1\cap D_2) \arrow[d] \arrow[r] & \IndCoh(D_1) \arrow[d] \\
	\IndCoh(D_2) \arrow[r]                   & \IndCoh(Z)            
	\end{tikzcd}
\end{equation}
This categorical cube is Beck-Chevalley by base change, see \cite[Cor.~3.1.4]{GR17}. 

\begin{defi}
We call the coproduct totalization $\tot^{\amalg}(\IndCoh(D_*))$ of $\IndCoh(D_\ast)$ the {\em
categorical intersection complex} of $\{D_i\}_{1\leq i\leq N}$.
\end{defi}

For the square \eqref{eq:indcohsq}, the categorical intersection complex has $3$ nontrivial terms given by 
\[
\IndCoh(D_1\cap D_2)\longrightarrow \laxpush{\IndCoh(D_1)}{\IndCoh(D_1\cap D_2)}{\IndCoh(D_2)}\longrightarrow \IndCoh(Z)\,.
\]
The $\infty$-category 
\[
	\laxpush{\IndCoh(D_1)}{\IndCoh(D_1\cap D_2)}{\IndCoh(D_2)}
\]
describes the ``lax gluing'' of the two schemes $D_1$ and $D_2$ along $D_1 \cup D_2$. We proceed with
describing situations in which the categorical intersection complex is spherical or allows a
Calabi--Yau structure. 

\begin{rem}
	\label{rem:lunts} The categories arising as the various terms in the
categorical intersection complex are closely related to the poset schemes
introduced in \cite{lunts:poset} (which correspond to the fully lax colimit of
the punctured cubical diagram $\IndCoh(D_*))$.  The data captured by the
categorical intersection complex seems to be somewhat more refined, and it
should be interesting to explore its role in the theory of noncommutative
resolutions of singularities. 
\end{rem}

\begin{prop}\label{prop:intcplxisspherical}
Suppose that $Z$ is a smooth, projective variety. Then $\IndCoh(D_\ast)$ is a spherical categorical
cube. In this case, the categorical intersection complex is hence a spherical categorical complex. 
\end{prop}
\begin{proof}
The sphericalness of the cube follows from the fact that the inclusion of any smooth divisor into a
smooth projective variety induces a spherical functor when passing to $\Ind$-coherent sheaves,
see for instance \cite{Add16}. The statement about the totalization follows from \Cref{prop:sphtot}.
\end{proof}

Recall the following Theorem from \cite{BD19}:

\begin{thm}[{\cite[Thm.~5.13]{BD19}}]\label{thm:gorcy}
Suppose that $Z$ is a Gorenstein scheme of dimension $n$ and $i \colon D \subset Z$ an anticanonical
divisor. Then $\IndCoh(D)$ admits a canonical left $(n-1)$-Calabi--Yau structure and the functor
$i_*\colon \IndCoh(D)\to \IndCoh(Z)$ admits a canonical compatible left $n$-Calabi--Yau structure.
\end{thm}

In terms of categorical complexes, this result equips the categorical $2$-term complex
\[
	\IndCoh(D) \overset{i_*}{\lra} \IndCoh(Z)
\]
associated to the anticanonical divisor $D \subset Z$ with an $n$-dimensional Calabi--Yau structure. We will
now show that, if $D = \bigcup_{1 \le i \le N} D_i$ is an anticanonical normal crossings divisor,
then the statement can be refined to provide an $n$-dimensional Calabi--Yau structure on the
categorical intersection complex of $\{D_i\}_{1\leq i\leq N}$. For the proof, we require the
sphericalness of the cube, and to this end we assume that $Z$ is smooth projective. 

\begin{thm}\label{thm:anticancubeCY}
	Let $Z$ be a smooth projective scheme of dimension $n$ and let $D = \bigcup_{1 \le i \le N}
	D_i$ be an anticanonical normal crossings divisor. Then the categorical $N$-cube 
	\[ 
		\IndCoh(D_*)
	\] 
	is spherical and admits a canonical left $n$-Calabi--Yau structure. 
\end{thm}

\begin{cor}\label{cor:anticantotCY}
	Let $Z$ be a smooth projective scheme of dimension $n$ and let $D = \bigcup_{1 \le i \le N}
	D_i$ be an anticanonical normal crossings divisor. Then the categorical intersection complex 
	\[ 
		\tot^{\amalg}(\IndCoh(D_*))
	\] 
	is spherical and admits a canonical left $n$-Calabi--Yau structure.
\end{cor}

\begin{proof}
Combine \Cref{thm:anticancubeCY}, \Cref{prop:sphtot} and  \Cref{thm:CYcubetot}.
\end{proof}

\begin{rem}
Theorem \ref{thm:anticancubeCY} and Corollary \ref{cor:anticantotCY} provide possible answers to the
question raised by Katzarkov-Kontsevich-Pantev, as to what kind of Calabi--Yau and spherical
categorical structures arise from an anticanonical normal crossings divisor, see Remark 4.36 in
\cite{KKP08}.
\end{rem}

Before proving \Cref{thm:anticancubeCY}, we recall some results from \cite[Section~5.2]{BD19}. Given a scheme $Z$, we denote by $\on{RHom}_Z(\mhyphen,\mhyphen)$ the $k$-linear derived Hom in
$\IndCoh(Z)$. We denote by $\omega_Z^\bullet$ the dualizing complex of $Z$, defined as $\pi^!(k)$
with $\pi\colon Z\to \ast$. Using the diagonal map $\Delta\colon Z\rightarrow Z\times Z$, one finds
a canonical morphism 
\[ \on{RHom}_{Z}(\O_Z,\omega^\bullet_Z)\rightarrow \on{RHom}_{Z\times Z}(\Delta_*\O_Z,\Delta_*\omega^\bullet_Z)\simeq \on{HH}(\IndCoh(Z))\,,\]
natural in $Z$.

A scheme $Z$ is called Cohen-Macaulay of dimension $d$ if $\omega_X= \omega_X^\bullet[-d]$ is a coherent sheaf. In this case, there exists an isomorphism $H_d\on{HH}(\IndCoh(Z))\simeq H_d\on{HH}^{S^1}(\IndCoh(Z))$. In the setting of \Cref{thm:gorcy}, we further have $i^!\omega_Z\simeq \omega_D[-1]$ and a cofiber sequence
\begin{equation}\label{eq:DZcofseq} \O_Z\rightarrow i_*\O_D\rightarrow \omega_Z[1]\,.\end{equation}
A straightforward computation shows that the morphism $i_*\O_D\rightarrow \omega_Z[1]$ is adjoint to an equivalence $\O_D\simeq \omega_D$, see \cite{BD19}. Furthermore, the morphism $i_*\O_D\to \omega_Z[1]$ factors through the equivalence $i_*\O_D\simeq i_*\omega_D$ via the counit $\counit$ of $i_*\dashv i^!$: 
\[i_*\O_D\simeq i_*\omega_D \simeq i_*i^!\omega_Z[1]\xrightarrow{\counit} \omega_Z[1]\,.\] 
It follows that the cofiber sequence \eqref{eq:DZcofseq} gives rise to a relative negative cyclic homology class describing a Calabi--Yau structure on $i_*$.\\

Recall our notation $D_J\coloneqq \bigcap_{j\in J} D_j$, for $J \subset [N]$. We further denote
$\tilde{D}_{J}\coloneqq \bigcup_{j\in [N]\backslash J} D_{J\cup \{j\}} \subset D_J$. The inclusions
$D_J\subset Z$ is denoted $f_J$. 

Each component $D_i \subseteq Z$ is cut out by a section $s_i\colon  \O_Z \to \O(D_i)$. We dualize $s_i$ to obtain a map $s_i^{\vee}:
\O(-D_i) \to \O_Z$ which is part of a cofiber sequence
\[
	\begin{tikzcd}
		\O(-D_i) \ar{r}\ar{d} & \O_Z \ar{d}\\
		0 \ar{r} & (f_i)_*\O_{D_i}\,.
	\end{tikzcd}
\] 
The tensor product of the morphisms
$s_i^{\vee}$, $1 \le i \le n$, yields a cubical diagram
\[
	p\colon  \Pc([N])^{\op} \to \IndCoh(Z),\; J \mapsto \O( - \sum_{i \in J} D_i).
\]
Further, by passing to iterated cofibers, or equivalently directly tensoring the cofiber morphisms, we obtain a ``reflected'' cubical diagram
\begin{equation}\label{eq:qcube}
	q\colon \Pc([N]) \to \IndCoh(Z),\; J \mapsto (f_J)_*\O_{D_J}\,.
\end{equation}
By construction, the total fiber of the cube $q$ is $\omega_Z$ so that its total cofiber is $\omega_Z[N]$. 
Applying the functor $\on{RHom}_{Z}(-,\omega_Z^\bullet)$ to the cube $q$, we obtain the cube
\[
	q_{\omega}\colon  \Pc([N])^{\op} \to \Modk,\; J \mapsto \on{RHom}_{Z}((f_J)_*\O_{D_J},\omega_Z^\bullet)\,.
\]
Finally, composing with the pushforward along the diagonal map $\Delta\colon Z\to Z\times Z$, we obtain the cube 
\[
	Q_{\omega}\colon  \Pc([n])^{\op} \to \Modk,\; J \mapsto \on{RHom}_{Z \times
	Z}(\Delta_*(f_J)_*\O_{D_J},\Delta_*\omega_Z^\bullet)\,.
\]

\begin{lem}\label{lem:equivcubes}
There is a natural equivalence of cubes
	\begin{equation}\label{eq:lemident}
		Q_{\omega} \simeq \HH(\IndCoh(D_*))\,.
	\end{equation}
\end{lem}

\begin{proof}
	We have equivalences 
	\[
		\on{RHom}_{Z \times Z}(\Delta_*(f_J)_* \O_{D_J},\Delta_*\omega_Z^\bullet) \simeq \on{RHom}_{D_J\times D_J}(\Delta_*\O_{D_J},\Delta_*\omega_{D_J}^\bullet)
	\]
	arising from the various adjunctions $(-)_* \dashv (-)^!$ and base change associated to the pullback square
	\[
		\begin{tikzcd}
			D_J  \ar{r}{\Delta}\ar{d}{f_J} & D_J \times D_J \ar{d}{f_J \times f_J}\\
			Z \ar{r}{\Delta} & Z\times Z 
		\end{tikzcd}
	\]
	along with the equivalence $f_J^!\omega_{Z}^\bullet \simeq \omega_{D_J}^\bullet$.
	These can be combined with the natural equivalences
	\[
		\on{RHom}_{D_J \times D_J}(\Delta_*\O_{D_J},\Delta_*\omega_{D_J}^\bullet) \simeq
		\HH(\IndCoh(D_J))
	\]
	to the desired equivalences of cubes. The remainder of the proof justifies that these equivalences coherently assemble into an equivalence of cubes.
	
	 In the cases $N=1,2$, the involved coherence issues can be deal with directly. For $N>2$, we show jointly the equivalence \eqref{eq:lemident} and the statement that $\HH$ preserves the colimit over $I^{N}_{<(0,\dots,0)}$ by induction. The fact that $\HH$ preserves the colimit in the case $N=2$ follows from the argument presented in the proof of \Cref{lem:totHHintcplx} below. We apply the decomposition of colimits \cite[4.2.3.10]{lurie:htt} and the induction assumption, to find that $\colim\nolimits_{<(0,\dots,0)}\HH(\IndCoh(D_\ast))$ is equivalent to the pushout of the diagram
	\[
	\begin{tikzcd}
{\HH(\IndCoh(\bigcup_{1\leq i\leq N-1}D_{i}\cap D_{1}))} \arrow[d] \arrow[r] & {\HH(\IndCoh(\bigcup_{1\leq i\leq N-1}D_{i}))} \\
\HH(\IndCoh(D_{\{1\}})\,,                                   &                         
\end{tikzcd}
\]
and thus given by $\HH(\Ind\Coh(\tilde{D}_{[N]}))$. Applying the induction assumption once more, we find an equivalence between the cubes $Q_\omega$ and $\HH(\IndCoh(D_\ast))$ punctured at the final vertices. Via the universal property of the colimit, we can extend this equivalence to the entire cubes by virtue of the commutative diagram 
\[
\begin{tikzcd}
{\on{RHom}_{Z \times Z}(\Delta_*i_* \O_{\tilde{D}_{[N]}},\Delta_*\omega_Z^\bullet)} \arrow[r, "\simeq"] \arrow[d] & {\HH(\Ind\Coh(\tilde{D}_{[N]})} \arrow[d] \\
{\on{RHom}_{Z \times Z}(\Delta_*\O_{Z},\Delta_*\omega_Z^\bullet)} \arrow[r, "\simeq"]                           & {\HH(\IndCoh(Z))\,.}          
\end{tikzcd}
\]
Here $i:\tilde{D}_{[N]}\subset Z$ denotes the inclusion. We have constructed the equivalence of cubes \eqref{eq:lemident}, concluding the induction step.
\end{proof}

\begin{lem}\label{lem:totHHintcplx}
	\begin{enumerate}[label=(\arabic*)]
	\item The total Hochschild homology $\HH^{\tot}(\IndCoh(D_*))$ is equivalent to the cofiber of 
\[ \on{RHom}_{Z\times Z}(\Delta_*i_*\O_{D},\Delta_*\omega^\bullet_{Z}) \xlongrightarrow{\circ \Delta_* \unit} \on{RHom}_{Z\times Z}(\Delta_*\O_Z,\Delta_*\omega^\bullet_{Z})\,,\]
with $\unit$ the canonical morphism $\O_Z\to i_*\O_D$.
\item The cofiber sequence in $\IndCoh(Z)$
\[ \O_Z\to i_*\O_D\to \omega_Z\]
determines a total negative cyclic homology class $\eta\colon k[n]\to \HH^{S^1,\tot}(\IndCoh(D_*))$.
\end{enumerate}
\end{lem}

\begin{proof}
Given any categorical $N$-cube $\A_\ast$, one has an equivalence
\[ \HH^{\tot}(\Ind\Coh(A_\ast))\simeq \cof(\colim\nolimits_{J<(0,\dots,0)} \HH(\Ind\Coh(A_J))\to \HH(\Ind\Coh(A)))\,.\]
Using \Cref{lem:equivcubes}, part (1) thus follows from the equivalence
\begin{align*}
\colim\nolimits_{<(0,\dots,0)}\HH(\Ind\Coh(D_\ast))& \simeq \colim\nolimits_{<(0,\dots,0)}Q_\omega \\ & \simeq \colim\nolimits_{J<(0,\dots,0)}\on{RHom}_{Z \times
	Z}(\Delta_*(f_J)_*\O_{D_J},\Delta_*\omega_Z^\bullet) \\
& \simeq		\on{RHom}_{Z \times
	Z}(\Delta_*\lim\nolimits_{(\I^N_{<(0,\dots,0)})^{\op}} (f_J)_*\O_{D_J},\Delta_*\omega_Z^\bullet)\\
& \simeq \on{RHom}_{Z \times
	Z}(\Delta_* i_*\O_D,\Delta_*\omega_Z^\bullet)\\
& \simeq \HH(\Ind\Coh(D))\,.
\end{align*}

We proceed with part (2). Pushing the given cofiber sequence forward along $\Delta_*$, we obtain a cofiber sequence
\[
			\begin{tikzcd}
\Delta_* \O_Z \arrow[r] \arrow[d] & \Delta_*  i_*\O_D \arrow[d] \\
0 \arrow[r]              & \Delta_* \omega_Z         
\end{tikzcd}
			\]
which determines by part (1) a class $\sigma\colon k[n]\to \HH^{\tot}(\IndCoh(D_*))$. This class lifts uniquely to $\HHStot(\IndCoh(D_*))$, as follows from the above mentioned equivalences 
\[ H_{n-|J|}\HH^{S^1}(\IndCoh(D_J))\simeq H_{n-|J|}\HH(\IndCoh(D_J))\,.
\]  
\end{proof}

\begin{rem}\label{rem:totHHintcplx}
There is an apparent analogue of \Cref{lem:totHHintcplx} with the essentially same proof: 

Let $J\subset [N]$ and denote by $\Pc([N])^\op_{\geq J}$ the subposet of  $\Pc([N])^\op$ consisting of those $J'\subset [N]$ satisfying that $J\subset J'$. Denote by $\IndCoh(D_{\ast\geq J})$ the restriction of $\IndCoh(D_{\ast})$ to $\Pc([N])^\op_{\geq J}$. The total Hochschild homlogy 
\[ \HHtot(\Ind\Coh(D_{\geq J})\]
is equivalent to the cofiber of 
\[ \on{RHom}_{D_J\times D_J}(\Delta_*(i_J)_*\O_{\tilde{D}_J},\Delta_*\omega^\bullet_{D_J}) \xlongrightarrow{\circ \Delta_* \unit} \on{RHom}_{D_J\times D_J}(\Delta_*\O_{D_J},\Delta_*\omega^\bullet_{D_J})\,,\]
with $i_J:\tilde{D}_J\subset D_J$ the inclusion. Any commutative diagram 
\[
\begin{tikzcd}
{\O_{D_J}} \arrow[d] \arrow[r] & {(i_J)_*(\O_{\tilde{D}_J})} \arrow[d] \\
0 \arrow[r]                  & \omega_{D_J}                       
\end{tikzcd}
\]
thus determines a total Hochschild class $\eta\colon k[n-|J|]\to \HHStot(\IndCoh(D_{\ast\geq J}))$. 
\end{rem}

\begin{proof}[Proof of \Cref{thm:anticancubeCY}] 
We show that the negative cyclic class $\eta$ from part (2) of \Cref{lem:totHHintcplx} describes the desired left $n$-Calabi-Yau structure of $\IndCoh(D_*)$. 

The recursive nondegeneracy conditions can be checked before pushing forward along the diagonal morphism $\Delta$. They amount to the following statements:
		\begin{enumerate}
			\item As in the argument for the proof of \Cref{thm:gorcy} recalled above, the first nondegeneracy condition amounts to the fact that the commutative diagram 
			\[
			\begin{tikzcd}
\O_Z \arrow[r] \arrow[d] & i_*\O_D \arrow[d] \\
0 \arrow[r]              & \omega_Z         
\end{tikzcd}
			\]
			giving rise to $\eta$ is a cofiber sequence. 
		\item Let $J=\{j_1,\dots,j_l\} \subset [N]$. Consider $(N+1)$-cube $\tilde{q}$, with $q^{\on{aug}}|_{\Pc([N+1])^\op_{>\{N+1\}}}=q$ (using that $\Pc([N+1])^\op_{>\{N+1\}}\simeq \Pc([N])$), $q^{\on{aug}}(\emptyset)=\omega_Z^\bullet$ and $q^{\on{aug}}(x)=0$ for all other $x\in \Pc([N+1])^\op$, exhibiting $\omega_Z^\bullet$ as the total cofiber of $q$. We obtain an $(N-|J|+1)$-cube $\tilde{q}^{\on{aug}}_J$ by restricting $q^{\on{aug}}$ to $\Pc([N+1])^\op_{>J\cup \{N+1\}}$. 
		Passing to partial colimits, we obtain the commutative diagram 
		\begin{equation}\label{eq:diagOJomegaS}
		\begin{tikzcd}
(f_{J})_*\O_{D_J} \arrow[d] \arrow[r] & (f_{J})_*(i_J)_*\O_{\tilde{D}_J} \arrow[d] \\
0 \arrow[r]                        & \omega_Z[|J|]                    
\end{tikzcd}
		\end{equation} 
		with $i_J\colon \tilde{D}_J\subset D_J$ the inclusion. The diagram \eqref{eq:diagOJomegaS} is adjoint under $(f_J)_*\dashv f_J^!$ to the following diagram:
		\begin{equation}\label{eq:diagOJomegaSadj}
		\begin{tikzcd}
O_{D_J} \arrow[d] \arrow[r] & (i_{J})_*\O_{\tilde{D}_J} \arrow[d] \\
0 \arrow[r]               & \omega_{D_J}\,,   
\end{tikzcd}
		\end{equation}
This diagram encodes as in \Cref{rem:totHHintcplx} the induced Hochschild homology class \[ \partial_{j_1}\dotsm\partial_{j_l}\eta\colon k[n-|J|]\to \HH^{S^1,\tot}(\IndCoh(D_{\ast\geq J}))\,.\]
The inductive argument below shows that the square \eqref{eq:diagOJomegaSadj} is biCartesian, yielding the non-degeneracy condition of this Hochschild class.
		
Let $j\in J$ and assume that the non-degeneracy of the class $\partial_{j_1}\dotsm\partial_{j_{l-1}}\eta$ has been shown. The restriction of $q$ to $\Pc([N])_{\geq J\backslash \{j_l\}}\coloneqq (\Pc([N])^\op_{\geq J\backslash \{j_l\}})^\op$ is given by the pushforward $(f_{J\backslash \{j_l\}})_*$ of a cube denoted $q_{J\backslash \{j_l\}}$. The cube $q_{J\backslash \{j_l\}}$ can be extended to a colimit cube $q^{\on{aug}}_{J\backslash \{j_l\}}$ , exhibiting $\omega_{D_{J\backslash \{j_l\}}}^\bullet$ as the total cofiber of $q_{J\backslash \{j_l\}}$; this colimit cube is adjoint under $(f_{J\backslash \{j_l\}})_*\dashv f_{J\backslash \{j_l\}}^!$ to $\tilde{q}^{\on{aug}}_{J\backslash \{j_l\}}$. 

The cube $q_{J\backslash \{j_l\}}$ arises from applying the unit of the adjunction $i_{J,j_l}^*\dashv (i_{J,j_l})_*$, with $i_{J,j_l}:D_J\subset D_{J\backslash \{j_l\}}$, to its face $\Pc([N])_{\geq J\cup \{j\}}\backslash \Pc([N])_{\geq J}$. We thus consider the $(N-J+1)$-cube $q_{J\backslash \{j_l\}}$ as a morphism between this face and the opposite face, and can pass to the cofiber morphism to obtain another $(N-J+1)$-cube $\on{cof}_{q_{J\backslash \{j_l\}}}$. The cube $\on{cof}_{q_{J\backslash \{j_l\}}}$ in turn arises from applying the counit of the adjunction $(i_{J,j})_*\dashv i_{J,j}^!$ to its face given by the cofiber $(N-|J|)$-cube, as follows from the sphericalness of $(i_{J,j})_*\dashv i_{J,j}^!$, see \cite[Lemma~2.5.15]{DKSS:spherical}. Partially totalizing $\on{cof}_{q_{J\backslash \{j_l\}}}$, we produce from $\on{cof}_{q_{J\backslash \{j_l\}}}$ the restriction of $q^{\on{aug}}_{J\backslash \{j_l\}}$ to one of its faces, and thus upon passing to a partial colimit the commutative square 
\[
		\begin{tikzcd}
(i_{J,j_l})_*\O_{D_{J}} \arrow[d] \arrow[r, "\alpha"] & (i_{J,j_l})_*(i_J)_*\O_{\tilde{D}_{J}} \arrow[d] \\
0 \arrow[r]                        & \omega_{D_{J\backslash \{j_l\}}}[1]                    
\end{tikzcd}
\]
which expresses a counit map $\on{cof}(\alpha)\to  \omega_{D_{J\backslash \{j_l\}}}[1]$. The adjoint square under $(i_{J,j_l})_*\dashv i_{J,j_l}^!$, given by \eqref{eq:diagOJomegaSadj}, is thus biCartesian, concluding the argument.
 \end{enumerate}
\end{proof}

\begin{rem}
The proof of \Cref{thm:anticancubeCY} shows that the divisor $\tilde{D}_J \subset D_J$ is
anticanonical for any $J\subset [N]$. This can also be checked directly using the adjunction
formula. 
\end{rem}

\begin{rem}\label{rem:orientationstrata} Note, that in the context of \Cref{thm:anticancubeCY}, the
	datum that defines the Calabi--Yau structure on the categorical intersection complex 
	\[ 
		\tot^{\amalg}(\IndCoh(D_*))
	\] 
	of the normal crossings divisor $D = \bigcup_{1 \le i \le N} D_i$ is identical to the datum
	that defines the relative Calabi--Yau structure on the functor
	\[
		\IndCoh(D) \overset{i_*}{\lra} \IndCoh(Z).
	\]
	Namely, in both cases this is a class in the relative Hochschild homology of $i_*$. This may be
	regarded as a noncommutative analog of the fact that an orientation on a manifold with
	corners induces a compatible system of orientations on all boundary strata.
\end{rem}

We conclude the section by noting that the above construction of the categorical intersection
complex of a normal crossings divisor can be generalized to the context of cubical resolutions of
schemes (see e.g. \cite{GNPP88}). 

\begin{defi}
	Let $S_\ast$ be a {\em $\I^n$-scheme}, i.e., a functor $\I^n \to \Schk$. A {\em
	$2$-resolution of $S_\ast$} consists of a commutative diagram of $\I^n$-schemes
	\[
	\begin{tikzcd}
	S_{1,1,\ast} \arrow[d] \arrow[r, "a_\ast"] & S_{1,0,\ast} \arrow[d, "f_\ast"] \\
	S_{0,1,\ast} \arrow[r, "b_\ast"]           & S_{0,0,\ast}               
	\end{tikzcd}
	\]
	with $S_{0,0,\ast}=S_\ast$ and satisfying for all $i\in [1]^n$ that $a_i$ and $b_i$ are closed
	immersions, $f_i$ is proper, $S_{1,0,i}$ is smooth and $f_i$ restricts to an isomorphism of schemes
	between $S_{1,0,i}\backslash f_i^{-1}(S_{0,1,i})$ and $S_{0,0,i}\backslash S_{0,1,i}$.
\end{defi}

\begin{exa}
Consider a scheme $S$ with smooth blowup $\on{BL}_Z(S)$ at a closed subvariety $Z$:
\[
\begin{tikzcd}
E \arrow[r] \arrow[d] \arrow[rd, "\lrcorner", phantom, near start] & \on{BL}_Z(S) \arrow[d] \\
Z \arrow[r]                                            & S                     
\end{tikzcd}
\]
Then the above diagram defines a $2$-resolution of $S$.  
\end{exa}

Let $S$ be any scheme. Then there exists a $2$-resolution $S^2$ of $S$ obtained as the pullback square of any resolution of $S$, see \cite[Thm 2.6]{GNPP88}. We may thus choose a $2$-resolution of $S^2$. We consider the first column of this $2$-resolution as a $\I$-scheme, which in turn admits a further $2$-resolution. We choose one such, denoted $S^3$, which we consider as a $\I^3$-scheme. We proceed in this way for all $n< N$, choosing $S^{n+1}$ as a $2$-resolution of the $\I^{n-1}$-scheme obtained from restricting $S^n$ along $\I^{n-1}\times \{1\} \hookrightarrow \I^n$. From this, we extract an apparent $\I^N$-scheme $S_\ast$ satisfying 
\begin{align*} 
S_{0,\dots,0}& = S\,,\\
S_{1,0,\dots,0}&= S^{2}_{1,0}\,,\\
& ~\vdots\\
S_{\ast,1,0,\dots,0}& = S^{j+2}_{\ast,1,0}\,\,\text{ for }\ast\in \I^{j}\,,\\
& ~\vdots\\
S_{\ast,1}& = S^N_{\ast,1}\,.
\end{align*}

\begin{defi}
Let $S$ be a scheme an $S_\ast$ an $\I^N$-scheme obtained as above. We call $S_\ast$ an $S$-augmented cubical hyperresolution if $S_i$ is smooth for all $i\in \I^N\backslash \{(0,\dots,0)\}$.
\end{defi}

\begin{exa}\label{ex:intres}
For $S=X_1\cup X_2\cup X_3$ the union of three smooth subschemes with smooth intersections, a cubical hyperresolution arises from restricting the following diagram consisting of two $2$-resolutions, to the 'outer' $3$-cube.
\[
\begin{tikzcd}[column sep=small]
X_1\cap X_2\cap X_3 \arrow[rr] \arrow[rd] \arrow[dd] &                                   & X_1\cap X_3 \arrow[rd] \arrow[dd]           &                        &                &                     \\
                                                     & X_2\cap X_3 \arrow[dd] \arrow[rr] &                                             & X_3 \arrow[dd]         &                &                     \\
X_1\cap X_2 \arrow[rd] \arrow[rr]                    &                                   & X_1\cap (X_2\cup X_3) \arrow[rd] \arrow[rr] &                        & X_1 \arrow[rd] &                     \\
                                                     & X_2 \arrow[rr]                    &                                             & X_2\cup X_3 \arrow[rr] &                & X_1\cup X_2\cup X_3
\end{tikzcd}
\]
More generally, given a scheme $S$, written as the union of smooth subschemes with smooth intersections, there is an apparent $\I^N$-scheme $S_\ast$ with $S_J=\bigcap_{j\in J}X_j$ for all $J\in \Pc([N])^\op\simeq \I^N$. The $\I^N$-scheme $S_\ast$ is a cubical hyperresolution of $S$. The $\I^N$-scheme associated to a normal crossing divisor inside a scheme $Z$ at the beginning of the section can be obtained from this hyperresolution by changing the value at the terminal vertex of the $N$-cube from $S=\bigcup_{j\in [N]}X_j$ to $Z$.  
\end{exa}

Given an $\I^N$-scheme $S_\ast$, we obtain a functor $\IndCoh(S_\ast)\colon
\I^N\rightarrow \St_k$ by composing with $\IndCoh$.  The diagram
$\IndCoh(S_\ast)$ describes a categorical cube, which is in general neither
spherical, nor does its totalization admit a Calabi--Yau structure. In the case
where $S_\ast$ is an $S$-augmented cubical hyperresolution, it seems an
intriguing question to explore the implications of considering the truncation
$\tau_{\geq 1}\tot^{\amalg}\IndCoh(S_\ast)$ as a resolution of the stable
$\infty$-category $\IndCoh(S)$ in terms of smooth $k$-linear stable
$\infty$-categories. 

\subsection{Picard-Lefschetz theory and Fukaya-Seidel complexes}\label{subsec:FukayaSeidelcplx}

We begin by explaining how classical Picard--Lefschetz theory can be used to construct certain cell
complexes modelling the cohomology of affine varieties. 

Let $X \subset \CC^{N}$ be an $n$-dimensional smooth affine subvariety. Let $X_1 \subset X$ be a
generic hyperplane section of $X$ in $\CC^{N}$. The pair $(X,X_1)$ induces an exact triangle 
\begin{equation}\label{eq:lefschetzles}
	\begin{tikzcd}
		S_{\bullet}(X) \ar{r} & S_{\bullet}(X,X_1) \ar{r}{\partial} & S_{\bullet-1}(X_1) \ar{r}{+1} &
		\phantom{x}
	\end{tikzcd}
\end{equation}
in singular homology. In this context, the classical Lefschetz hyperplane theorem, due to, in the
given affine setup, to Andreotti-Frankel, implies that the complex $S_{\bullet}(X,X_1)$ 
\begin{enumerate}[label=\arabic*.]
	\item has homology concentrated in degree $n$, and
	\item $H_n(X,X_1)$ is generated by Lefschetz thimbles.
\end{enumerate}

In particular, the exact triangle \eqref{eq:lefschetzles} provides a means of computing the homology
$H_{\bullet}(X)$ in terms of the group $H_n(X,X_1)$ and the homology $H_{\bullet}(X_1)$. We may
iterate this construction, choosing a sequence of generic hyperplane sections
\[
	X \supset X_1 \supset X_2 \supset \cdots \supset X_n
\]
so that the the corresponding exact triangles of consecutive pairs combine to give a description of
the singular homology $H_{\bullet}(X)$ as the homology of the complex
\begin{equation}\label{eq:cellcomplex}
		H_n(X,X_1) \to  H_{n-1}(X_1,X_2) \to  ...  \to  H_0(X_n) 
\end{equation}

Under suitable technical assumptions and with suitable choices of symplectic structures (cf.
\cite{Sei08}, P. Seidel showed that a categorical variant of the above discussion can be implemented
to provide an effective means for computing Fukaya categories of affine varieties in terms of
Lefschetz fibrations. Namely, suppose that the hyperplane section $X_1 = \pi_1^{-1}(\{1\})$ is given
as the fiber of a Lefschetz fibration $\pi_1\colon X \to \CC$ (such as the restriction of a generic
linear function on $\CC^N$). Then we have a ``left exact'' sequence of categories
\[
	\Fuk(X) \hra \FS(\pi_1) \overset{\partial}{\lra} \Fuk(X_1)
\]
i.e., the category $\Fuk(X)$ is the kernel of the functor $\partial$. As shown in \cite{Sei08}, the
category $\FS(\pi_1)$ admits an exceptional collection of objects given by Lagrangian Lefschetz
thimbles. 

Again, we may iterate this consideration to produce a categorical complex 
\begin{equation}\label{eq:fscomplex}
	\FS(\pi_{\bullet}) \coloneqq  \FS(\pi_1) \to \FS(\pi_2) \to ... \to \Fuk(X_n)
\end{equation}
where $\FS(\pi_k)$ is the Fukaya--Seidel category of a Lefschetz fibration $\pi_k\colon  X_{k-1} \to \CC$. 

\begin{rem}
	\label{rem:fsdecategorify}
	Note that, since each of the Fukaya--Seidel categories $\FS(\pi_k)$ has an exceptional
	collection given by Lagrangian Lefschetz thimbles of $\pi_k$, we have, by the Lefschetz
	hyperplane theorem, an isomorphism of abelian groups $K_0(\FS(\pi_k)) \cong
	H_n(X_{k-1},X_{k})$. Further, the complex $K_0(\FS(\pi_{\bullet}))$ reproduces the complex
	\eqref{eq:cellcomplex}.
\end{rem}

Note that in the original reference \cite{Sei08} all Fukaya categories are described in the framework of $A_{\infty}$-categories. We may turn them into $k$-linear $\infty$-categories by choosing quasi-equivalent dg-categories and passing to their $k$-linear dg-nerves, see \cite{Coh13}. We can also apply the $A_\infty$-nerve construction of \cite{Fao17}, but the $k$-linear structure of the resulting $\infty$-categories has not yet been described. In order to obtain presentable
$\infty$-categories, we then pass to $\Ind$-completions: 

\begin{defi}
	\label{def:fukaya-seidel}
	The categorical complex $\Ind(\FS(\pi_{\bullet}))$ is called the {\em ($\Ind$-completed) Fukaya--Seidel
	complex} of the family of Lefschetz fibrations $\{\pi_{\bullet}\}$.
\end{defi}

\begin{rem}
	\label{rem:wrapped} As will be seen in the examples below, the categories of cycles of the
	complex $\Ind(\FS(\pi_{\bullet}))$ correspond to wrapped variants of the Fukaya categories
	considered in \cite{Sei08} (the Fukaya categories in \cite{Sei08} are always generated by
	compact Lagrangian submanifolds). 
\end{rem}

\subsection{Mirror symmetry}
\label{sub:mirror}

In the context of homological mirror symmetry, we expect equivalences between the categorical
complexes from \S \ref{subsec:completeintersection} and \S \ref{subsec:FukayaSeidelcplx},
respectively. 

In this section, we illustrate this ``mirror symmetry for categorical complexes'' in a somewhat familiar
context and then formulate a conjecture as to what to expect in greater generality. To this end,
we start recalling a well-known example of homological mirror symmetry (see \cite{seidel:more,AKO08}
for details):

We consider the affine hypersurface
\[
	X = \{ xyz = 1 \} \subset \CC^3
\]
and note that $X \cong (\CC^*)^2$. Consider the Lefschetz fibration $\pi_1\colon X \to \CC$ given by the
restriction of the linear function $\pi_1 = x + y + z$. The general fiber of $\pi_1$ is an elliptic
curve with $3$ punctures which degenerates to a nodal cubic over the $3$ critical values of $\pi_1$.
As a regular fiber, we may take $X_1 \coloneqq  \pi_1(0) = E \setminus \{p_1,p_2,p_3\}$. A depiction
of the $3$ vanishing spheres in $X_1$ that correspond to the
critical values can be found in \cite[\S 3B]{seidel:more}:
\begin{center}
	\includegraphics[scale=.5]{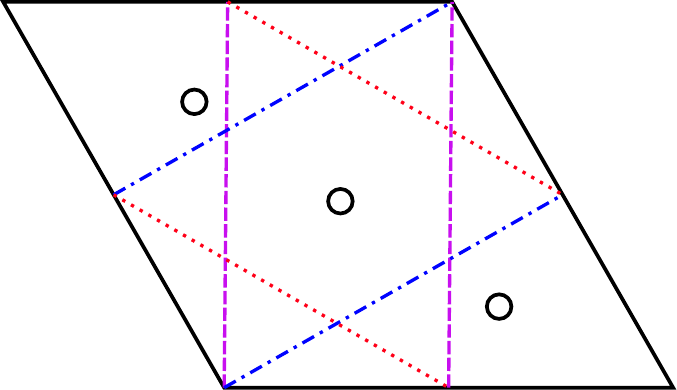}
\end{center}
This example does not directly fit into the context discussed in \cite{Sei08}, but it is explained
in \cite{seidel:more} how to associate to $\pi_1$ a $\ZZ$-graded Fukaya-Seidel category which can be
described as follows: The category $\FS(\pi_1)$ is equivalent to the derived category of the quiver
\[
	\begin{tikzcd}
		\bullet \ar[bend left=50]{r}{a_1} \ar{r}{a_2} \ar[bend right=50]{r}{a_3} & \bullet \ar[bend
		left=50]{r}{b_1} \ar{r}{b_2} \ar[bend right=50]{r}{b_3} & \bullet
	\end{tikzcd}
\]
subject to the relations $b_i a_j = b_j a_i$, for $i \ne j$, and $b_i a_i = 0$. Thus, we further
have
\[
	\FS(\pi_1) \simeq \Coh(\PP^2)
\]
by \cite{beilinson:pn}. 

As a next step, we choose the Lefschetz fibration $\pi_2 \colon X_1 \to \CC$ obtained by restricting the
linear function $\pi_2 \coloneqq  x-y$ on $\CC^3$. This is a ramified covering map of degree $3$, the fibers of
$\pi_2$ move in the pencil generated by $x-y$ and $w$ on $E$ whose general fiber consists of $3$
distinct points degenerating over the $6$ nondegenerate critical
values. Again, this setup does not directly fit into \cite{Sei08}, but the necessary modifications
are explained in \cite{seidel:more}. We may take the fiber $X_2 \coloneqq  \pi_2^{-1}(0) =
\{a,b,c\}$ as a regular fiber. Then the $6$ critical values correspond to $6$ vanishing cycles
($0$-spheres = pairs of points) in $X_2$:
\[
	\{a,b\}, \{a,b\}, \{a,c\}, \{a,c\}, \{b,c\}, \{b,c\} 
\]
The directed subcategory on these can be described by the quiver
\[
	\begin{tikzcd}
		\bullet \ar[bend left=20]{r}{a} \ar[bend right=20]{r}{b} & \bullet \ar{r}{a} \ar[bend
		left=50]{rrr}{b} & 
		\bullet \ar[bend left=20]{r}{a} \ar[bend right=20]{r}{c} & \bullet \ar{r}{c} &
		\bullet \ar[bend left=20]{r}{b} \ar[bend right=20]{r}{c} & \bullet
	\end{tikzcd}
\]
with zero relations given by the rule that the composite of composable arrows is zero iff they are
labelled by different letters. We omit the discussion of grading choices (which is explained in
\cite{seidel:more}).

\begin{figure}[p]
\[
	\includegraphics{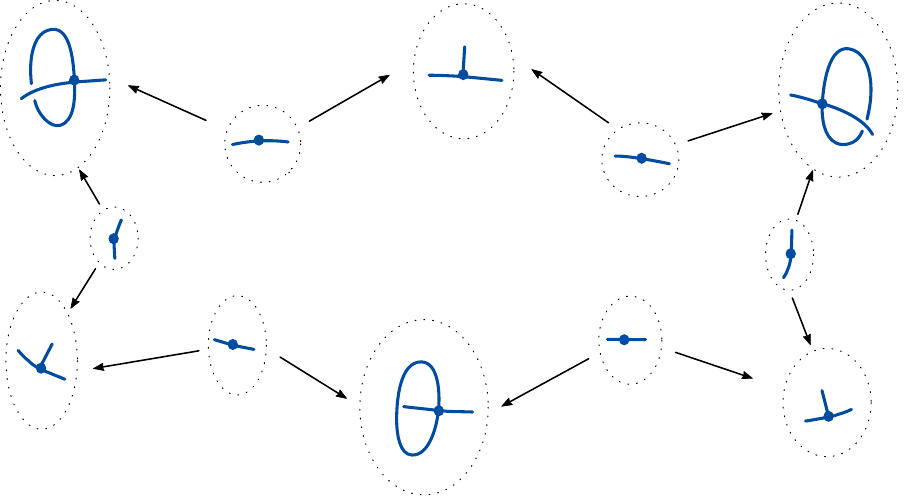}
\]
\caption{Decomposition of $\Gamma$.}
\label{fig:ribbon-split}
\end{figure}

\begin{figure}[p]
\centering
\def\svgwidth{\columnwidth}
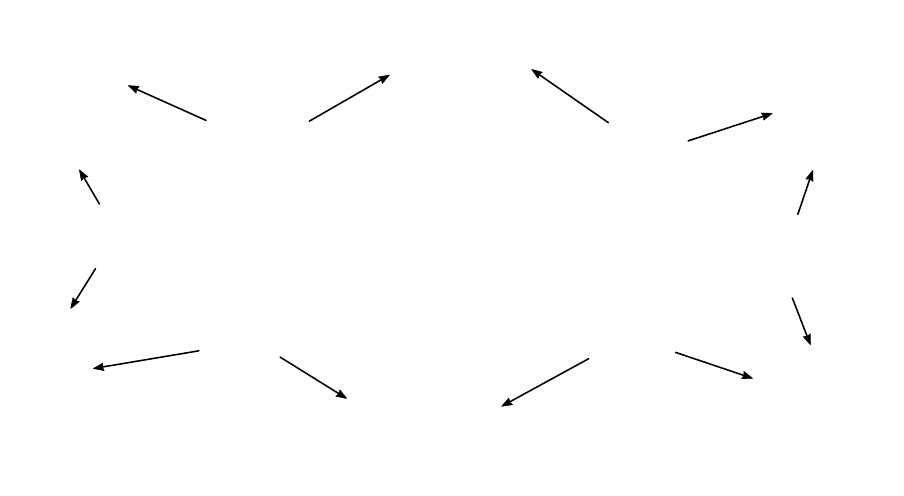
\caption{Diagram computing the topological Fukaya category of $\Gamma$.}
\label{fig:ribbon-diagram}
\end{figure}

Alternatively, we may describe this category as a topological Fukaya category of the Riemann surface
$E \setminus \{p_1,p_2,p_3\}$ with one stop at each puncture (or rather the cylindrical end
corresponding to it). The corresponding spine is given by the Ribbon graph
\[
	\Gamma =\quad \raisebox{-.4\height}{\includegraphics{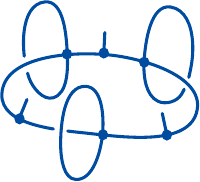}}
\]
so that the formalism of \cite{d:a1homotopy} yields the diagram of Ribbon graphs depicted in Figure
\ref{fig:ribbon-split} This latter diagram yields an equivalence between the topological Fukaya
category and the colimit of the corresponding diagram of Figure \ref{fig:ribbon-diagram}. Here,
$\Rep(A_2)$ denotes the bounded derived category of representations of the quiver
\[
	A_2 = \quad \begin{tikzcd}
		\bullet \ar{r} & \bullet,
	\end{tikzcd}
\]
the functors of the form $\Coh(\pt) \to \Rep(A_2)$ are induced by the inclusion of
source (resp. target) of the quiver. The functors of the form $\Coh(\pt) \to \Coh(\PP^1)$ are given
by pushforward along the inclusion of the points $0$ and $\infty$, respectively, into $\PP^1$. 
Again, there are auxiliary choices to be made to determine the $\ZZ$-grading, such as a
trivialization of the tangent bundle (or rather its square) of $E \setminus \{p_1,p_2,p_3\}$, which we
do not discuss here. 

Informally, the colimit of the diagram from Figure \ref{fig:ribbon-diagram} can thus be described by
starting with three disjoint copies of $\Coh(\PP^1)$ and freely adjoining three arrows connecting
the various skyscraper sheaves of the projective lines at $0$ and $\infty$, respectively. We
keep this intuition in mind by denoting the resulting category by
\begin{equation}
	\label{eq:three_lax_lines}
	\Coh(\raisebox{-.4\height}{\includegraphics[scale=.7]{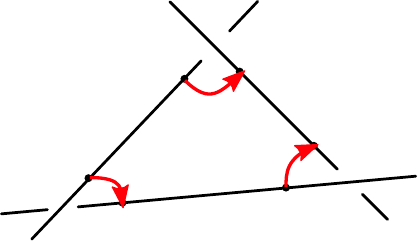}})
\end{equation}
describing it as an amalgamate of commutative geometry (the projective lines) and noncommutative
geometry (the quiver arrows). In terms of the constructions of \S \ref{subsec:laxsum}, this means that the
category is described by the \dpo{} of three copies of $\Coh(\PP^1)$ along three copies of
$\Coh(\pt)$. 


The above choices of potentials $\pi_1, \pi_2$ in total yield a Fukaya-Seidel complex of the form
\begin{equation}
	\label{eq:fs_p2}
		\FS(\pi_1) \lra \FS(\pi_2) \lra \Fuk(X_2) 
\end{equation}
with
\begin{enumerate}
	\item $\FS(\pi_1) \simeq \Coh(\PP^2)$
	\item $\FS(\pi_2) \simeq \Fuk(\Gamma) \simeq
		\Coh(\raisebox{-.4\height}{\includegraphics[scale=0.4]{three_lax_lines.pdf}})$
	\item $\Fuk(X_2) \simeq \Coh(\pt \amalg \pt \amalg \pt)$
\end{enumerate}

On the other hand, consider $\PP^2$ with the normal crossings divisor given by three lines $L_1, L_2,
L_3$ in general position: 
\begin{center}
\begingroup%
  \makeatletter%
  \providecommand\color[2][]{%
    \errmessage{(Inkscape) Color is used for the text in Inkscape, but the package 'color.sty' is not loaded}%
    \renewcommand\color[2][]{}%
  }%
  \providecommand\transparent[1]{%
    \errmessage{(Inkscape) Transparency is used (non-zero) for the text in Inkscape, but the package 'transparent.sty' is not loaded}%
    \renewcommand\transparent[1]{}%
  }%
  \providecommand\rotatebox[2]{#2}%
  \newcommand*\fsize{\dimexpr\f@size pt\relax}%
  \newcommand*\lineheight[1]{\fontsize{\fsize}{#1\fsize}\selectfont}%
  \ifx\svgwidth\undefined%
    \setlength{\unitlength}{247.30120802bp}%
    \ifx\svgscale\undefined%
      \relax%
    \else%
      \setlength{\unitlength}{\unitlength * \real{\svgscale}}%
    \fi%
  \else%
    \setlength{\unitlength}{\svgwidth}%
  \fi%
  \global\let\svgwidth\undefined%
  \global\let\svgscale\undefined%
  \makeatother%
  \begin{picture}(1,0.35608567)%
    \lineheight{1}%
    \setlength\tabcolsep{0pt}%
    \put(0,0){\includegraphics[width=\unitlength,page=1]{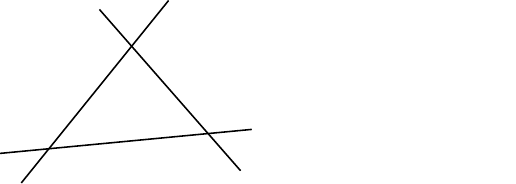}}%
    \put(0.13647325,0.19279556){\color[rgb]{0,0,0}\makebox(0,0)[t]{\lineheight{1.25}\smash{\begin{tabular}[t]{c}$L_1$\end{tabular}}}}%
    \put(0.26319841,0.02707806){\color[rgb]{0,0,0}\makebox(0,0)[t]{\lineheight{1.25}\smash{\begin{tabular}[t]{c}$L_2$\end{tabular}}}}%
    \put(0.33793377,0.2155411){\color[rgb]{0,0,0}\makebox(0,0)[t]{\lineheight{1.25}\smash{\begin{tabular}[t]{c}$L_3$\end{tabular}}}}%
    \put(0.67261819,0.16138503){\color[rgb]{0,0,0}\makebox(0,0)[t]{\lineheight{1.25}\smash{\begin{tabular}[t]{c}$\subset \PP^2$\end{tabular}}}}%
  \end{picture}%
\endgroup%

\end{center}

The corresponding categorical cubical diagram induced by the various push-forward functors takes the
following form:
\[
\begin{tikzcd}[row sep=small, column sep=small]
	\Coh(\emptyset) \arrow[dd] \arrow[rr] \arrow[rd]   &                                      &
	\Coh(L_2 \cap L_3) \arrow[dd] \arrow[rd] &                                   \\
                                     & \Coh(L_1 \cap L_3) \arrow[dd] \arrow[rr] &
				     & \Coh(L_3) \arrow[dd] \\
\Coh(L_1 \cap L_2 ) \arrow[rr] \arrow[rd] &                                      & \Coh(L_2) \arrow[rd]    &                                   \\
                                     & \Coh(L_1) \arrow[rr]    &
				     & \Coh(\PP^2)           
\end{tikzcd}
\]
Its coproduct totalization yields the categorical intersection complex 
\begin{equation}
	\label{eq:ic_p2}
	\begin{tikzcd}
		\Coh(\pt^{\amalg 3}) \ar{r} &
		\Coh(\raisebox{-.4\height}{\includegraphics[scale=0.4]{three_lax_lines.pdf}})\ar{r} & \Coh(\PP^2) 
	\end{tikzcd}
\end{equation}
where we use the above notation \eqref{eq:three_lax_lines} to denote the middle term. In particular,
all terms of the complexes \eqref{eq:ic_p2} and \eqref{eq:fs_p2} agree. We do not verify here the
expectation that the differentials in both complexes are in fact adjoint to one another (this can
probably be extracted with some effort from the existing literature) but rather formulate the more
general conjecture: 

\begin{conj}\label{conj:mirror} Let $X = (\CC^*)^n \subset \CC^{n+1}$. Then there exists a family of iterated Lefschetz
	fibrations $\pi_{\bullet}$ on $X$ such that the corresponding Fukaya--Seidel complex
	$\FS(\pi_{\bullet})$ is equivalent to the (right adjoint of the) categorical intersection complex
	of the normal crossings divisor given by $n+1$ hyperplanes in $\PP^{n}$ in general
	position.\end{conj}

The conjecture has its natural generality in some context of normal crossings divisors in
toric Fano varieties, where Hori-Vafa mirror symmetry provides a prediction of the mirror potential \cite{hori2000mirror}, but for the sake of concreteness, we leave it as stated.

While presenting some of this material during a talk in the Edinburgh Hodge Seminar, we learned
about some work in progress on wrapped Fukaya categories of ``multi-potentials'' which seems to be
closely related to the perspective on homological mirror symmetry via categorical complexes. Indeed,
once the multi-potential approach is implemented, one expects that this directly yields cubical
diagrams of wrapped Fukaya categories which are mirror to the cubical diagrams of coherent sheaves
from which we build the categorical intersection complex via totalization.

This ``mirror symmetry conjecture for cubes'' has for example been described in the recent article
\cite{Lee22}. Once this type of cubical mirror symmetry is established, our Conjecture
\ref{conj:mirror} then essentially reduces to a statement that our Fukaya-Seidel complexes from
Definition \ref{def:fukaya-seidel} arise as totalization of the cubical wrapped Fukaya-category
diagrams that one expects to associate to a multi-potential Landau--Ginzburg model (cf.
\cite{auroux:multi}). This statement, as well as its relevance for higher--dimensional perverse
schobers (cf.\S \ref{subsec:sphcub} and \S \ref{subsec:spherical}), seems interesting in its own
right and will be investigated in future work. 

\subsection{Manifolds with corners and complexes of \texorpdfstring{$\infty$}{infinity}-local systems} 
\label{subsec:localsystems}

We fix a field $k$. Given a space $X$, we denote by $\on{Loc}(X)=\on{Fun}(X,\mathcal{D}(k))$ the stable $\infty$-category of $\mathcal{D}(k)$-valued local systems on $X$. Given a morphism $f\colon X\rightarrow Y$ between spaces, we denote by $f^*\colon \on{Loc}(Y)\rightarrow \on{Loc}(X)$ the pullback functor. It admits left and right adjoints $f_!$ and $f_*$, given by left Kan extension and right Kan extension, respectively, see \cite[4.3.3.7]{lurie:htt}. Recall the following Theorem from \cite{BD19}.

\begin{thm}\label{thm:mfd}
Let $X$ be a compact oriented manifold of dimension $n$ with boundary $f\colon \partial X\subset X$. Then the functor 
\[ f_!\colon \on{Loc}(\partial X)\longrightarrow \on{Loc}(X)\]
admits a canonical left $n$-Calabi--Yau structure.
\end{thm}

The goal of this section is to show that \Cref{thm:mfd} admits an extension to categorical complexes arising from oriented bordisms. In the simplest case, an oriented $n$-dimensional bordisms between two closed, oriented $(n-1)$-dimensional manifolds $M,M'$ consists of an oriented $n$-dimensional manifold $N$ with boundary $\partial N=M\amalg M'$, with $M$ and $M'$ carrying the induced orientation. If instead $M$ and $M'$ are not closed and not disjoint in $N$, but instead overlap on their boundary, one can further ask that $M$ and $M'$ again define two $(n-1)$-dimensional bordisms between $(n-2)$-dimensional manifolds, and so on. This perspective can be formalized by organizing oriented bordisms up to dimension $n$ in an $(\infty,n)$-category of oriented bordisms $\on{Bord}_n^{\on{or}}$, whose $m$-cells are $m$-dimensional oriented bordisms for $m\leq n$, see for instance \cite{Lur:tft,CS19}. In the following, we associate to a functor from the $n$-simplex to the $(\infty,n)$-category $\on{Bord}_n^{\on{or}}$ a categorical $n$-cube whose totalization is left Calabi--Yau.

We fix an integer $n\geq 1$. Given an $n$-simplex $\Delta^n$, we denote by $\on{Sd}\Delta^n$ its barycentric subdivision, considered as a poset. We denote by $\on{Mfd}$ the $1$-category with objects compact oriented manifolds of any dimension with boundary and morphisms oriented inclusions into the boundary. A functor from the $n$-simplex to $\on{Bord}_n^{\on{or}}$ amounts to a functor $A_\ast\colon \on{Sd}\Delta^n\rightarrow \on{Mfd}$, mapping each $m$-simplex $x\in (\Delta^n)_m$ to an $m$-dimensional manifold $A_x$ (possibly empty), such that the boundary of $A_x$ is given by
\[ \partial A_{x}=\bigcup_{i=0}^m A_{d_i(x)}
\]
for any $x\in (\Delta^n)_m$. We further ask that the intersection of $A_{d_i(x)}$ and $A_{d_j(x)}$ is given by $A_{d_jd_i(x)}$ for any $j<i$.

There is an equivalence of posets $\phi\colon \on{Sd}\Delta^n\simeq  \I^n\backslash \{(1,\dots,1)\}$, i.e.~the barycentric subdivision of $\Delta^n$ is a cube with the initial vertex removed. Given an $n$-simplex $A_\ast$ in $\on{Bord}_n^{\on{or}}$, we hence find a categorical cube 
\[ \on{Loc}(A_\ast)\colon \I^n\rightarrow \Stk\,,\]
with $\on{Loc}(A_\ast)_{\phi(x)}=\on{Loc}(A_{\phi(x)})$ and $\on{Loc}(A)_{((1,\dots,1))}=0$. We will omit the usage of the equivalence $\phi$ in the following, and write for instance $A_J$ and $\on{Loc}(A_J)$ for $A_{\phi^{-1}(J)}$ and $\on{Loc}(A_{\phi^{-1}(J)})$, with $J\subset [n]$.  

\begin{thm}\label{thm:bord}
Let $A_\ast$ be an $n$-simplex in $\on{Bord}_n^{\on{or}}$. Then the categorical $n$-cube  $\on{Loc}(A_\ast)$ admits a canonical left $n$-Calabi--Yau structure. 
\end{thm}

\begin{cor}\label{cor:bord}
Let $A_\ast$ be an $n$-simplex in $\on{Bord}_n^{\on{or}}$. Then the coproduct totalization $\on{tot}^{\amalg}(\on{Loc}(A_\ast))$ admits a canonical left $n$-Calabi--Yau structure.
\end{cor}

\begin{proof}
Combine \Cref{thm:bord} and \Cref{thm:CYcubetot}.
\end{proof}

Before proving \Cref{thm:bord}, we recall some details from the proof of \Cref{thm:mfd} from \cite{BD19}. Given a space $X$, we denote by $C_\bullet(X)\in \mathcal{D}(k)$ its singular chain complex and by $LX$ its free loop space. There are natural morphisms in $\mathcal{D}(k)$
\begin{equation}\label{eq:SingulartoHochschild} C_\bullet(X)\rightarrow C_\bullet(LX)\simeq  \on{HH}^{S^1}(\on{Loc}(X))\,,\end{equation}
where the first morphisms arises from the inclusion of the constant loops $X\rightarrow LX$. If $X$ is a closed $n$-dimensional manifold, any orientation, considered as an element of $H_nC_\bullet(X)$, gives via the above map rise to a left $n$-Calabi--Yau on $\on{Loc}(X)$. If $Y$ is an $n$-dimensional manifold with boundary $X$, then a relative orientation of $Y$ is a suitable $n$-th homology class of $C_\bullet(Y,X)\coloneqq \on{cof}(C_\bullet(X)\rightarrow C_\bullet(Y))$ and its image in $\HH^{S^1}(\on{Loc}(Y),\on{Loc}(X))$ gives rise to a relative Calabi--Yau structure on the functor $\on{Loc}(X)\xrightarrow{(X\subset Y)_!} \on{Loc}(Y)$. 

\begin{lem}\label{lem:totalSingulartoHochschild}
There exists a natural map 
\[ C_\bullet(A_{(0,\dots,0)},\partial A_{(0,\dots,0)})\longrightarrow \HH^{S^1,\tot}(\on{Loc}(A_\ast))\,.\]
\end{lem}

\begin{proof}
Consider the diagram $C_\bullet(A_\ast)\colon \I^n\to \Mod_k$. The colimit of its restriction to $\I^n_{<(0,\dots,0)}=\I^n\backslash \{(0,\dots,0)\}$ is by Mayer-Vietoris equivalent to $C_\bullet(\partial A_{(0,\dots,0)})$. Consider the following diagram in $\mathcal{D}(k)$:
\[
\begin{tikzcd}
{C_\bullet(\partial A_{(0,\dots,0)})} \arrow[r] \arrow[d]                         & {C_\bullet(A_{(0,\dots,0)})} \arrow[d] \arrow[r] & {C_\bullet(A_{(0,\dots,0)},\partial A_{0,\dots,0)})} \arrow[d, dotted] \\
{\colim\nolimits_{J<(0,\dots,0)} \HH^{S^1,\tot}(A_{J})} \arrow[r] & {\HH^{S^1}(A_{(0,\dots,0)})} \arrow[r]     & {\HH^{S^1,\tot}(A)}                                     
\end{tikzcd}
\]
The left square commutes by the naturality of the morphisms \eqref{eq:SingulartoHochschild}. Both the upper and lower sequence in this diagram are cofiber sequences. It follows that there exists a dotted arrow as indicated, making the diagram commutative. 
\end{proof}

\begin{proof}[Proof of \Cref{thm:bord}.]
Let $x\in \Delta^n$ be the top cell. Any choice of orientation of $A_{x}$ relative to $\partial A_x$ determines by \Cref{lem:totalSingulartoHochschild} a total negative cyclic homology class $\eta\colon k[n]\to \HH^{S^1,\tot}(\on{Loc}(A_\ast))$. Its restriction to $\HH^{S^1,\tot}(\on{Loc}(A_{\partial_i \dots \partial_j \ast}))$ arises via an analogue of \Cref{lem:totalSingulartoHochschild} for the cube $A_{\partial_i \dots \partial_j \ast}$ from an induced orientation of $A_{d_i\dots d_j(x)}$ relative to $\partial A_{d_i\dots d_j(x)}$. 

Using that passing to local systems preserves colimits, we find for any $J\subset [n]$ a commutative diagram
\[
\begin{tikzcd}[row sep=tiny]
{\on{Loc}(\partial A_{J}))} \arrow[r, "=", no head] & \on{Loc}(\bigcup_{j\not\in J}^m A_{J\cup \{j\}}) \arrow[rd] \arrow[dd, "\simeq"', no head] &                            \\
                                                                                  &                                                                                                  & {\on{Loc}(A_{J})\,.} \\
                                                                                  & {\colim_{\I^n_{<J}}\on{Loc}(A_\ast)} \arrow[ru]                                                    &                           
\end{tikzcd}
\]
The total Hochschild class of $\on{Loc}(A_\ast)$ thus restricts by \Cref{thm:mfd} to a Calabi-Yau structure on the lower horizontal morphism above. This shows that $\eta$ defines a left $n$-Calabi-Yau structure for $\on{Loc}(A_\ast)$, as desired. 
\end{proof}

\begin{exa}
Consider the $(n+m)$-simplex $A$ in $\on{Bord}_{n+m}$, with 
\[ A_{[n+m]}=\mathbb{D}^n\times \mathbb{D}^m\] the product of the $n$- and $m$-dimensional real unit discs and  
\[ A_{d_0[n+m]}=(\partial \mathbb{D}^n)\times \mathbb{D}^m=S^{n-1}\times \mathbb{D}^m\,,\quad A_{d_1[n+m]}=\mathbb{D}^m\times S^{m}\,,\quad A_{d_0^2[n+m]}=S^{n-1}\times S^{m-1}\]
and $A_{x}=\emptyset$ otherwise. The totalization $\on{tot}^{\amalg}(\on{Loc}(A))$ is given by the complex
\[
\on{Loc}(S^{n-1}\times S^{m-1})\longrightarrow \laxcolima{\on{Loc}(\mathbb{D}^n\times S^{m-1})}{\on{Loc}(S^{n-1}\times S^{m-1})}{\on{Loc}( S^{n-1}\times \mathbb{D}^m)}\longrightarrow \on{Loc}(\mathbb{D}^n\times \mathbb{D}^m)\]
concentrated in degrees $2$ to $0$. Except for being left $(n+m)$-Calabi--Yau, this complex is furthermore spherical as it arises from totalizing the following spherical square:
\[
\begin{tikzcd}
\on{Loc}(S^{n-1}\times S^{m-1}) \arrow[d] \arrow[r] & \on{Loc}(\mathbb{D}^n\times S^{m-1}) \arrow[d] \\
\on{Loc}( S^{n-1}\times \mathbb{D}^m) \arrow[r]     & \on{Loc}(\mathbb{D}^n\times \mathbb{D}^m)     
\end{tikzcd}
\]
The above square arises from a pullback square of Kan fibrations and is hence Beck-Chevalley, see \cite[Prop.~4.4.11,~Thm.~6.4.13]{Cis}. These Kan fibrations are furthermore spherical fibrations, meaning that their fibres are spheres, the sphericalness of the functors was thus shown in \cite{Chr20}. We note that in general, $\on{tot}^{\amalg}(\on{Loc}(A))$ is not a spherical complex. 

A version of this example producing an $n$-term complex arises by starting with the product of $n+1$ discs. 
\end{exa}





\bibliographystyle{alpha} 
\bibliography{refs} 

\end{document}


\subsection{}
\newcommand\lQ{Q}
\newcommand\lI{I}
\newcommand\loneC{C}
\newcommand\lC{\C}
\newcommand\Funp{\Fun_{*}}
\newcommand\Funz[1][0]{\Fun_{#1}}
\newcommand\lquot{\lQ/\lI}
\newcommand\lquotred{
  \left(
    \lQ/\lI
  \right)'
}

Let \(\lQ\) be a small category and \(\lI\subset \lQ\) an ideal,
i.e., a set of morphisms such that \(f\circ g \in \lI\)
whenever \(f\in\lI\) or \(g\in\lI\).

We denote by \(\lquotred\) the category with the same objects
as \(\lQ\) and with hom-sets
\begin{equation}
  \lquotred(x,y)=\set{0}\amalg \lQ(x,y)\setminus \lI
\end{equation}
with the unique composition law which makes the assignment
\begin{equation}
  \lQ\to\lquotred;
  \quad\quad f\mapsto
  \begin{cases}
    f,\text{ if } f\notin\lI,\\
    0,\text{ if } f\in\lI
  \end{cases}
\end{equation}
into a functor. This is well defined because \(\lI\) is an ideal.
Furthermore we denote by \(\lquot\) the
pointed category obtained from \(\lquotred\) by adding a new zero object \(0\)
and defining the composition
\begin{equation}
  x\to 0 \to y
\end{equation}
to be the morphism \(0\in\lquotred(x,y)\) for all objects \(x,y\).
Note that the inclusion \(\lquotred\hookrightarrow\lquot\) is fully faithful.

The following is a straightforward observation.

\begin{lem}
  \label{lem:1-categorical-quot-extension}
  Let \(\loneC\) be a pointed category.
  The restriction functors along
  \(\lQ\to\lquotred\to\lquot\)
  induce equivalences of
  between the full subcategories
  \begin{align}
    \Funp\left(\lquot,\loneC\right)
    &
      \hookrightarrow
      \Fun\left( \lquot,\loneC \right)
    \\
    \Funz\left( \lquotred,C \right)
    &
      \hookrightarrow
      \Fun\left( \lquotred,C \right)
    \\
    \Funz[I]\left( Q \right)
    &
      \hookrightarrow
      \Fun(Q,C)
  \end{align}
  \begin{itemize}
  \item
    of functors \(\lquot\to C\) which are pointed, i.e.,
    send the zero object to a zero object;
  \item
    of functors \(\lquotred \to C\) which send zero arrows to zero morphisms;
  \item
    of functors \(\lQ\to C\) which send arrows in \(\lI\) to zero morphisms.
  \end{itemize}
\end{lem}

Now let \(\lC\) be an \(\infty\)-category.
We still have the full subcategory
\begin{equation}
  \Funp(\lquot,\lC)
  \subset\Fun(\lquot,\lC)
\end{equation}
of pointed functors.
Additionally, we define the full subcategories
\begin{align}
  \Funz(\lquotred,\lC)
  &\coloneqq \Fun(\lquotred,\lC)
  \times_{\Fun(\lquotred,\ho\lC)}\Funz(\lquotred,\ho\lC)
  \\
  &=\setP{F}
  {F(0\colon x\to y) = 0 \text{ in } \pi_0\loneC(F(x),F(y))}
  \\
  \Funz[I](\lQ,\lC)
  &\coloneqq
  \Fun(\lQ,\lC)
  \times_{\Fun(\lQ,\ho\lC)}
  \Funz[I](\lQ,\ho\lC)
  \\
  &= \setP{F}
    {\forall f\in\lI : F(f\colon x\to y)= 0 \text{ in } \pi_0\loneC(F(x),F(y))}
\end{align}
of \(\Fun(\lquotred,\lC)\) and \(\Fun(\lQ,\lC)\), respectively,
where we impose the zero condition on the morphisms up to homotopy.

In general, the analog of \Cref{lem:1-categorical-quot-extension}
does not hold:
For example, let
\(\lQ=\set{1\xrightarrow{a} 2\xrightarrow{b} 3}\)
and \(\lI=\set{ba}\).
Then extending a diagram \(F\colon \lQ\to\lC\) along \(\lQ\to \lquot\),
involves the choice of a nullhomotopy \(F(ba)\simeq 0\)
in the mapping space \(\lC(F(1),F(3))\);
while the condition on \(F\) guarantees its existence,
it will in general not be unique if the component of the zero map
is not contractible in \(\lC(F(1),F(3))\).
More generally, extending a diagram along \(\lQ\to\lquot\)
involves choosing coherent trivializations of all morphisms \(F(f)\),
for \(f\in I\).
Each choice of trivialization and each choice of coherence datum
takes place in the mapping spaces of \(\lC\);
more specifically in the connected components of the zero morphisms
\(0\colon x\to 0\to y\).
Thus it is not surprising that all these choices can be made uniquely
when all these spaces are contractible.

\begin{lem}
  \label{lem:infty-categorical-quot-extension}
  Let \(\lC\) be a pointed \(\infty\)-category
  satisfying the following assumption:
  \begin{itemize}
  \item
    For each \(x,y : \lC\), the connected component of the zero morphisms
    \(x\to 0\to y\) in \(\lC(x,y)\) is contractible.
  \end{itemize}
  Then restriction along \(Q\to \lquotred\to\lquot\)
  induce equivalences of \(\infty\)-categories
  \begin{equation}
    \Funp(\lquot,\lC)
    \xrightarrow{\simeq}
    \Funz(\lquotred,\lC)
    \xrightarrow{\simeq}
    \Funz[I](\lQ,\lC).
  \end{equation}
\end{lem}


\todo{finish proof}